\documentclass{article}
\usepackage{amsmath,latexsym,amssymb,amsthm,enumerate,amscd}
\theoremstyle{plain}
\newtheorem{theorem}{Theorem}[section]
\newtheorem{proposition}[theorem]{Proposition}
\newtheorem{corollary}[theorem]{Corollary}
\newtheorem{lemma}[theorem]{Lemma}
\theoremstyle{definition}
\newtheorem{definition}[theorem]{Definition}

\newtheorem{remark}[theorem]{Remark}
\newtheorem{example}[theorem]{Example}

\newtheorem{thm}{Theorem}

\theoremstyle{plain}

\theoremstyle{definition}

\newtheorem{rem}[thm]{Remark}

\newtheorem*{urem}{Remark}

\newtheorem*{uthm}{Theorem}

\newtheorem*{ack}{Acknowledgment}

\numberwithin{equation}{section}
\numberwithin{table}{section} 
\DeclareMathOperator{\Grass}{\rm{Grass}}

\def\GA{\mathrm{GA}}
\def\LA{\mathrm{LA}}

\def\cha{\mathrm{char}\ }

\def\Ann{\mathrm{Ann}\ }
\def\Proj{\mathrm{Proj}\ }

\def\Hilb{\mathrm{Hilb}}
\def\Sat{\mathrm{Sat}}
\def\Soc{\mathrm{Soc}}

\def\GCD{\mathrm{GCD}}
\def\Grass{\mathrm{Grass}}

\providecommand{\bysame}{\makebox[3em]{\hrulefill}\thinspace}
\def\cod{\mathrm{cod}\ }
\def\<{\left<}
\def\>{\right>}

\def\GAD{\mathrm{GAD}}

\def\ns{\footnotesize \it}

\def\Sat{\mathrm{Sat}}

\def\Z{\mathfrak{Z}}

\def\max{\mathrm{max}}
\def\rem{\mathrm{rem}}
\def\Anc{\mathrm{Anc}}
\title{Ancestor ideals of vector spaces of forms, and level
algebras}
\author{Anthony Iarrobino\\[.05in]
{\ns Department of Mathematics, Northeastern University, Boston, MA 02115, USA.
}\\[.2in]}

\date{May 22, 2003}

\begin{document}

\maketitle
\begin{abstract} Let $R=k[x_1,\ldots ,x_r]$ denote the polynomial ring
in $r$ variables over a field $k$, with maximal ideal $M=(x_1,\ldots
,x_r)$, and let
$V\subset R_j$ denote a vector subspace of the space $R_j$ of degree-$j$
homogeneous elements of $R$. We study three related algebras determined
by $V$. The first is the \emph{ancestor algebra} $\Anc
(V)=R/\overline{V}$ whose defining \emph{ancestor ideal} $\overline{V}$
is the largest graded ideal of
$R$ such that $\overline{V}\cap M^j=(V)$, the ideal generated by $V$. The
second is the level algebra $\LA (V)=R/L(V)$ whose defining ideal $L(V)$,
is the largest graded ideal of $R$ such that the degree-j component
$L(V)\cap R_j$ is $V$; and third is the algebra $R/(V)$. We have that 
$L(V)=\overline{V}+M^{j+1}$. When $r=2$ we
determine the possible Hilbert functions $H$ for each of these algebras,
and as well the dimension of each Hilbert function stratum (Theorem
\ref{dimstrata}). We characterize the graded Betti numbers of these
algebras in terms of certain partitions depending only on $H$, and give
the codimension of each stratum in terms of invariants of the partitions
(Theorem \ref{codpartition}). We show that when
$r=2$ and $k$ is algebraically closed the Hilbert function strata for each
of the three algebras attached to $V$ satisfy a frontier property that the
closure of a stratum is the union of more special strata. In each case
the family
$G(H)$ of all graded quotients of $R$ having the given Hilbert function is
a natural desingularization of this closure (Theorem
\ref{closureofstrata}). \par
 We then solve a
refinement of the simultaneous Waring problem for sets of degree-j binary
forms. Key tools throughout include properties of an invariant
$\tau (V)$, the number of generators of
$\overline{V}\subset k[x_1,x_2]$, and
previous results concerning the projective variety $G(H)$ in \cite{I2}. 
\end{abstract}
\tableofcontents
\section{Introduction}
In Section \ref{defnintro} we first define what we term the \emph{ancestor
ideal} $\overline{V}$ and \emph{ancestor algebra} $\Anc (V)$ and also the
\emph{level algebra}
$\LA (V)$ of a vector space 
$V\subset R_j$ of degree-$j$ forms in the polynomial ring
$R=k[x_1,\ldots ,x_r] $ in
$r$ variables over a field $k$. We then show some initial results about
the three algebras
$\Anc (V), \LA (V)$ and $R/(V)$ determined by $V$. In Section
\ref{main} we state our main results about these three algebras for $r=2$,
and we give  context in the literature. In Section \ref{HFstrata} we show
some general results about the Hilbert function strata of ancestor
ideals. In Section \ref{anc2} we show our main results about the three
algebras of $V$ for
$r=2$ variables. In Section \ref{HF2} we determine the dimensions of the
Hilbert function strata (Theorem \ref{dimstrata}); in
Section~\ref{minressub} we express the codimensions of these strata in
terms of partitions given by the graded Betti numbers of the three
algebras attached to
$V$ (Theorem \ref{codpartition}); and in Section \ref{closuresub} we
determine the Zariski closure of each Hilbert function stratum when $k$
is algebraically closed. We show that the strata for each of the three
algebras satisfy the frontier property, that the closure is a union of
more special strata in a natural partial order (Theorem
\ref{closureofstrata}).  In Section
\ref{subsimultWaring} we study a refinement of the simultaneous Waring
problem for vector spaces of degree-$j$ forms when
$r=2$. In Section
\ref{related} we develop a concept of \emph{related} vector spaces of
forms, then we state some open problems.
\par
\subsection{Three algebras attached to the vector
space
$V\subset R_j$}\label{defnintro}
We let
$k$ be an arbitrary field, and we denote by
$R=k[x_1,\ldots ,x_r]$ the polynomial ring over $k$, with maximal ideal
$M=(x_1,\ldots ,x_r)$, and the standard grading. For an integer
$j\ge 0$ we denote by
$R_j$ the vector space of degree-$j$ homogeneous elements of $R$. Let
$j>0$ and suppose that 
$V\subset R_j$ is a vector subspace of the space of degree-$j$ homogeneous
forms of $R_j$. We denote by
$(V)$ the ideal generated by $V$, and by $\overline{V}$ the largest ideal
of $R$ such that $\overline{V}\cap M^j=(V)$ (see Definition
\ref{basics}). For a form
$f\in R_j$ and an integer
$i\ge 0$ we denote by
$R_i\cdot f$ the vector space 
\begin{equation*}R_if=\langle hf\mid h\in R_i\rangle \subset R_{i+j};
\end{equation*}
For a vector space $V\subset R_j$ and an integer $i\ge 0$ we
denote by
$R_iV$ the vector space span 
\begin{equation}R_iV=\{ hf\mid h\in R_i, f\in V\} .
\end{equation}
 For $0\le
i\le j$ we denote by $R_{-i}V$ the vector space
 satisfying
\begin{equation}R_{-i}V=\{ f\in R_{j-i}\mid f\cdot R_i\subset
V \} .
\end{equation}
We now define the three algebras determined by $V$ that we study.
\begin{definition}\label{basics} Let $V\subset R_j$ be a vector space of
forms. The \emph{level ideal} $L(V)$
determined by
$V$ is
\begin{equation}L(V)=M^{j+1}\oplus V\oplus R_{-1}V\oplus \cdots \oplus
R_{-j}V,
\end{equation}
and the {\em level algebra} determined by $V$ is $\LA (V)=R/L(V)$. The 
\emph{ancestor ideal} $\overline{V}$ of $V$ is the ideal
\begin{equation}\overline{V}=(V)\oplus R_{-1}V\oplus \cdots
\oplus R_{-j}V,
\end{equation}
and the \emph{ancestor algebra} determined by $V$ is $\Anc
(V)=R/\overline{V}$. The usual ideal determined by $V$ is $(V)\subset
R_j$, and we denote by $\GA (V)=R/(V)$ the graded algebra quotient. \par
Recall that the {\emph socle} of an Artinian algebra $A=R/I$ is
\begin{equation*}
\Soc (A)=(0:M)_A =\langle f\in A\mid M\cdot f=0\rangle .
\end{equation*}
The {\em{type}} of $A$ is the vector space dimension $\dim_k(\Soc (A))$
of the socle.
\end{definition}
\begin{remark}\label{ancid} The ancestor ideal $\overline{ V}$ is the
largest graded ideal of
$R$ such that $\overline{V}\cap M^j=(V)$, the ideal of $R$ generated by
$V$. The level ideal
$L(V)$ is the largest graded ideal of $R$ such that $L(V)\cap R_j=V$: it
satisfies 
$L(V)=\overline{V}+M^{j+1}$; and the socle of the level algebra $\LA
(V)=R/L(V)$ satisfies
$\Soc (\LA (V))\cong R_j/V$.  The
ideal
$(V)$ satisfies
$(V)=\overline{V}\cap M^j$.
{\it Note:} The maximality statements
for the ancestor ideal
$\overline{V}$ and for the level ideal $L(V)$ may appear similar, but they
are quite different. The two ideals are equal only when
$R_1\cdot V=R_{j+1}$. 
\end{remark}
\begin{proof}[Proof of Remark] For $i>0, R_{-i}V\subset R_{i-j}$
is the largest subset of $R_{i-j}$ satisfying $R_i(R_{-i}V)\subset V$; and
evidently
$\overline{V}$ of Definition \ref{basics} is the largest graded ideal such
that
$\overline{V}\cap M^j=(V)$, the ideal generated by $V$. The other
statements are also immediate from the relevant definitions.
\end{proof}
\begin{lemma}
There are exact sequences
\begin{align}
0&\to \overline{V}/(V)\to R/(V)\to R/\overline{V}\to 0, \,\text {  and 
}\notag\\
 0&\to M^j/(V)\to R/\overline{V} \to R/L(R_{-1}V)\to 0.\label{ebasicexact}
\end{align}
\end{lemma}
\begin{proof} Immediate from the definitions.
\end{proof}
\begin{example}(See \cite[\S 60ff]{Mac1},\cite[Lemma 2.14]{IK}). When the codimension
of
$V$ as a vector subspace of
$R_j$ is one, then $\LA (V)=R/L(V)$ is a graded Artinian Gorenstein algebra,
and all standard graded Artinian Gorenstein algebras quotients of $R$
having socle degree
$j$  arise in this way. When $V=\langle xy^2+yx^2,x^3,y^3\rangle\subset
R=k[x,y]$ then $L(V)=(x^2+xy+y^2,x^3)$ and $\LA (V)$ is a complete
intersection of Hilbert function $H(A)=(1,2,2,1)$. Here, as usual in the Gorenstein Artinian case,
$\overline{V}=L(V)$; the exception is when $V=(m_p)\cap R_j$ for the
maximal ideal of a point $p\in {{\mathbb P}}^{r-1}$, then
$\overline{V}=m_p$.
\end{example}
\begin{example}
Let $I_\Z$ be the defining ideal of a closed subscheme $\Z\subset
{{\mathbb P}}^{r-1}$, and let
$V=I_\Z\cap R_j$. Then $\overline{V}\subset I_\Z$. If also $j\ge \sigma (\
Z)$, the regularity degree of $\Z$, then $\overline{V}=I_\Z$.
\end{example}
\par Recall that the saturation $\Sat (I)$ of a graded ideal $I\subset R$
is the ideal 
\begin{equation}
\Sat(I)=I:M^{\infty}=\{f\mid \exists i\text{ with } R_if\subset I\}.
\end{equation} 
Denote by $\sigma (V)$ the Castelnuovo-Mumford regularity degree of the
projective scheme
$\Z_V=\Proj (R/(V))\subset {{\mathbb P}}^{r-1}$. In case $(V)\supset
M^\sigma $ but $(V)\nsupseteq M^{\sigma -1}$, when $\Z (V)$ is empty, we
set
$\sigma (V)=\sigma$. We denote this same integer $\sigma(V)$ also by
$\sigma(\Anc (V))$ and $\sigma (\overline{V})$. 
\begin{lemma}\label{satV} Let $V\subset R_j$ be a vector subspace. For
$i\ge 0$, 
\begin{equation}\label{einclusion}
R_i\cdot R_{-i}\cdot V\subset V, \text{ and
}R_{-i}\cdot R_i V\supset V.
\end{equation}
  When  $V\ne R_j$ we have 
\begin{equation}\label{eqancseqi}
0=\overline{R_{-j}V}\subset \cdots
\subset\overline{R_{-1}V}\subset\overline{ V},
\end{equation}
and
\begin{equation}\label{eqancseqii}
 \overline{V}\subset \overline{R_1V}\subset \overline
{R_2V}\subset \cdots \subset Sat((V)).
\end{equation}
Also, for $i\ge \sigma (\Anc (V))-j,\,$ we have $ \overline{R_iV}=\Sat
((V))$.
\end{lemma}
\begin{proof} The inclusions of equation \eqref{einclusion} are immediate
from the definitions, and they imply equations
\eqref{eqancseqi} and \eqref{eqancseqii} (see also Lemma \ref{compare}).
The increasing sequence of ideals of equation \eqref{eqancseqii} evidently
terminates in
$\Sat ((V))$. Concerning the last claim, that $\overline{R_iV}=\Sat (V)$
for $i\ge \sigma (V)-j$ we first note that, taking $ W=R_{\sigma
-j}V$; that $\sigma (V)=\sigma$ implies $\sigma(W)=\sigma$. When
$R_1W=R_{\sigma +1}$ the claim is trivially satisfied; otherwise the
regularity degree of $\Proj (R/(W)$ is $\sigma$. It follows that $W=\Sat
((W))_\sigma$, and $\overline {W}=\Sat ((W))$. This completes
the proof.
\end{proof}
\begin{lemma}\label{AncandI} Let
$I$ be a graded ideal of $R$ satisfying $H(R/I)=H$, and let $V=I_j$. Then we have
\begin{equation}\label{eancineq}
I+M^{j+1}\subset \overline{V}+M^{j+1} \text { and } I\cap M^j \supset (V).
\end{equation}
\end{lemma}
\begin{proof} Let $a>0$ and $ i=j-a$, then we have  $V=I_j\supset R_a
I_i$, hence
\begin{equation*}
\overline{V}_i=R_{-a}\cdot V\supset R_{-a}R_{a}I_i\supset I_i
\end{equation*}
by \eqref{einclusion} of Lemma \ref{satV}. This shows $I+M^{j+1}\subset
\overline{V}+M^{j+1}$. Now let $a>0$ and $i=j+a$. We have
$R_aV=R_aI_j\subset I_i$, hence $I\cap M^j \supset (V)$. 
\end{proof}\par
\begin{definition}\label{def1.6} Let $V\subset R_j$ and $W\subset
R_i$. We say that
$V$ is \emph{equivalent} to $W$ ($V\equiv W$) if $\overline{V}=\overline
{W}$. We will say that 
$W$ is \emph{simpler than} $V$ if $W=R_{i-j}V$ and
$\overline{W}\ne\overline{V}$.
\end{definition}
The first principle behind this article is that each vector space in one
of the sequences  
\begin{equation*} 
V,R_{-1}V,R_{-2}V,\ldots \, \text { or }\, V,R_1V, R_2V \ldots
\end{equation*}
 should be either equivalent to or simpler than the preceding space. The
complexity of a vector space $V\subset R_j$ should be measured by an
invariant
$\tau (V)$ that is nonincreasing along each sequence
above, and where equality
$\tau (V)=\tau (R_iV)$ implies $V\equiv W$. We succeed in this enterprise
of measuring the complexity of $V$ only when $r=2$. In this case, we take
$\tau (V)=\dim_k R_1V-\dim_k V$, and show that $\tau (V)=\nu
(\overline{V})$, the number of generators of the ancestor ideal of $V$
(Lemma \ref{taugen}). We show that this $\tau$ has the needed
properties (Theorem \ref{taubasic}). When $r\ge 3$ an analogous
invariant with such strong properties is not possible due to an
example of D. Berman (Example
\ref{berman}).
\par The second principle is that, fixing a degree $j$ and
vector space dimension $d$, the Grassmanian 
$\Grass (d,R_j)$ parametrizing $d$-dimensional subspaces of $V\subset R_j$
is stratified by locally closed subschemes $\Grass (H)=\Grass_H(d,j)$,
parametrizing the vector spaces $V$ for which the Hilbert function
$H(R/\overline{V})=H$ is fixed. Letting $G(H)$ be the scheme
parametrizing all the graded ideals
$I\subset R$ with $H(R/I)=H$, we have that $\Grass(H)$ is an open
subscheme of $G(H)$ (Theorem \ref{loclostrata}). Natural questions are,
when is $\Grass (H)$ nonempty? Is $\Grass (H)$ irreducible? What are the
dimensions of its components? Is $\Grass (H)$ smooth? Describe the
Zariski closure
$\overline {\Grass (H)}\subset \Grass (d,R_j)$. 
\subsection{Background and main
results}\label{main} We first give the immediate background of the paper,
and outline our main results, then we discuss related work of others.\par
 Our main results are
for the case
$r=2$, where we answer the above questions. We further show that $G(H)$
is a natural desingularization of
$\overline{\Grass (H)}$ when $r=2$, and we determine the fibre of $G(H)$
over a point in the closure of
$\Grass (H)$.  \par
When
$r=2$ we denote by $\Grass_\tau (d, R_j)$ the locally closed subscheme of
$\Grass (d,R_j)$ parametrizing vector spaces $V$ with $\tau (V)=\tau$.
Recall that here,
$\tau (V)$ is the number of generators of $\overline
{V}$.  Given a sequence $H=(H_0,H_1,\ldots )$ of non-negative
integers, we define the first difference sequence $E(H)=\Delta
H$ by
\begin{equation}
E(H)=(e_1,\ldots ,e_i,\ldots ) ,\text{ where } e_i=H_{i-1}-H_i.
\end{equation}
We let $e_0=-1$. When $H=H(R/\overline{V}), \text{ then } e_i=\tau
(R_{i-j}V)-1$ for $i<j$, and $e_i=\tau
(R_{i-j-1}V)-1$ for $i>j$ (Proposition \ref{etau}). For $H',H$ two
sequences of integers that occur as Hilbert functions of ancestor
algebras $\Anc (V), V\subset R_j$, $\dim V=d$ we let (see Definition \ref{defpo})
\begin{equation}\label{partial}
H'\ge_{\mathcal P} H \text{ if for each
 } i\le j\text{ we have } H'_i\le H_i,\text{ and for each } i\ge j 
\text{ we
have }H'_i\ge H_i.
\end{equation}
 We denote by $a^+$ the number $a$ if $a\ge 0$ and $0$
otherwise. It is well known that in two variables, the Hilbert function
$H$ of a quotient
$A=R/I$ by a proper non-zero ideal (so $H$ is a \emph{proper}
$O$-sequence) satisfies, for some positive integer $\mu$, the
\emph{order} of $H$ (so $M^\mu\supset I,I_\mu\ne 0$)
\begin{align}
H&=(1,2,\ldots
,\mu,H_\mu,H_{\mu +1},\ldots ,H_i,\ldots ) \text{ with }\mu=\min\{i\mid
H_i<i+1\},\text{ and }\notag\\
&\mu\ge H_\mu\ge H_{\mu +1}\ge \cdots \ge c_H \text{ and }\lim_{i\to
\infty}H_i=c_H\ge 0.\label{Hcond}
\end{align}
\begin{definition}\label{defGH}
Given a sequence
$H$ satisfying \eqref{Hcond} with $c_H=0$, let $\sigma =\sigma_H$ satisfy
$H_{\sigma-1}\not=0,H_\sigma =0$. We denote by
$G(H)$ the closed subscheme 
\begin{equation}\label{eqGH}
G(H)\subset \prod_{ \mu\le i\le
\sigma-1}\Grass(i+1-H_i,R_i)
\end{equation}
 parametrizing graded ideals of $R$ having Hilbert
function $H$:  here $\prod_{ \mu\le i\le
\sigma-1}\Grass(i+1-H_i,R_i)$ parametrizes sequences $V_\mu,V_{\mu
+1},\ldots V_{\sigma -1}$ of vector spaces with each $V_i\subset R_i$ and
$\dim V_i=i+1-H_i$;
we assume
$V_i=0$ for
$i<\mu$ and
$V_i=R_i$ for
$i>j$. The subscheme $G(H)$ is defined by the conditions $xV_i\subset
V_{i+1}$ and $yV_i\subset V_{i+1}$ for $\mu\le i<j$. \par
When $c_H>0$, let $\sigma_H=\min\{i\mid H_{i-1}>c_H\}$. It is not hard to
show that each ideal
$I$ with
$H(R/I)=H$,  satisfies 
\begin{equation}
\exists f\in R_{c_H} \mid i>\sigma_H\Rightarrow I_i=(f)\cap R_i.
\end{equation}
Thus, when $c_H>0$ we may regard $G(H)\subset \prod_{\mu\le i\le
\sigma}\Grass(i+1-H_i,R_i)$, in a manner similar to that above in
\eqref{eqGH} for the case $c_H=0$.
\end{definition}
 We will use the
following result, essentially from
\cite{I2}, valid over a field
$k$ of
arbitrary characteristic.
\begin{theorem}\label{basic}\cite[Theorems 2.9,2.12,3.13,4.3, Proposition
4.4, Equation 4.7]{I2} Let
$r=2$, and for \eqref{basici} let the field $k$ be
algebraically closed. Let 
$H$ be an
$O$-sequence that is eventually constant, so $H$ is a sequence satisfying
\eqref{Hcond}, let $c=c_H$ and let $H_s=c_H, H_{s-1}\not=c_H$.
\begin{enumerate}[i.] 
\item\label{basici} Then $G(H)$ is a smooth projective variety of
dimension $c+\sum_{i\ge \mu} (e_i+1)(e_{i+1})$. $G(H)$ has a finite
cover by opens in an affine space of this dimension. If $\cha k=0$ or 
$\cha k>s$ then $G(H)$ has a finite cover by opens that are affine spaces.
\item\label{basicii} \cite[Theorem 4.3]{I2} The number of
generators
$\nu (I)$ of a graded ideal
$I$ for which
$H(R/I)=H$, satisfies $\nu (I)\ge \nu (H)=1+e_\mu +\sum_{i\ge
\mu}(e_{i+1}-e_i)^+$.
\item\label{basiciii}\cite[Proposition 4.4]{I2} Assume
that
$k$ is an infinite field. 
The graded ideals
$I$ with
$H(R/I)=H$ and having the minimal number $\nu (H)$ of generators given by
equality in
\eqref{basicii} form an open subscheme of $G(H)$ having the dimension
specified in \eqref{basici}, that is dense in $G(H)$ when $k$ is
algebraically closed.
\end{enumerate}
\end{theorem}
\begin{proof}[Remark on the Proof] The proof of \eqref{basici} in the
case $R/I$ Artinian, so 
$c=0$ is one of the main results of \cite{I2}. The characteristic 0 case
is handled in Theorems 2.9, 2.12, and the characteristic p case in
Theorem 3.13 of \cite{I2}. The proof of \eqref{basici} when 
$c>0$ relies on the fact that $t_s=t_{s+1}=c$ implies there is a form $f$
of degree $c$ such that
$I_s=(f)\cap R_s,I_{s+1}=(f)\cap R_{s+1}$ (for a proof, see Proposition
\ref{taubasic} \eqref{taubasicvi} below). This implies that
$f\mid I_i$ for $i\le s$. Thus, when $c>0,\, I=fI'$ where $I'$ is a
graded ideal such that $H(R/I')=H'$,  where
$H'$ is defined by $H'_i=H_{i+c}-c$. It follows that
$G(H)\cong
{\mathbb P}^c\times G(H')$. Here $H'$ is eventually zero, so the
dimension and structure of $G(H')$ is given by Theorems 2.9, 2.12, and
3.13 (see also Equation 4.7) of
\cite{I2}. In \cite{I2} we defined certain subfamilies $U_H\subset G(H)$
parametrizing ideals $I$ having ``normal patterns'': such that $I$ has a
Gr\"{o}bner basis with leading terms the first
$i+1-H_i$ degree-$i$ monomials in lexicographic order for
each
$i$. We showed that these subfamilies are affine spaces of dimension
specified in
\eqref{basici}; this result in fact requires only that $k$ be an
infinite field. However, that $U_H$ be dense in $G(H)$ requires that $k$
be algebraically closed.
\end{proof}\par
We will show the following main results for ancestor ideals of a vector
space
$V\subset R_j$ of homogeneous polynomials  when
$r=2$. Analogous results for level algebras and the algebras
$R/(V)$ follow, and are stated in the appropriate section. Recall that we
denote
$\Grass_H(d, R_j)$ by
$\Grass(H)$, and that we have $e_i=E(H)_i=H_{i-1}-H_i$. We denote by
$c_H=\lim_{i\to
\infty} H_i$. Theorem A is Theorem \ref{allH}\eqref{allHii}. Theorem B is
\eqref{edimH} of Theorem \ref{dimstrata}\eqref{dimstrataiii}; other
dimension results are in Theorems
\ref{dimstrata} and \ref{codpartition}. Theorems C,D are 
the two parts of Theorem
\ref{closureofstrata}, Theorem E is Theorem \ref{equivpos}.  For Theorems
B-E we assume that the field $k$ is infinite, and the
$O$-sequences
$H,H'$ belong to the set
$\mathcal H(d,j)$ of
\emph{acceptable} sequences (Definition \ref{defaccept}), which by
Corollary \ref{accept} are those $O$-sequences $H$ with $d$ fixed
satisfying the conditions of Theorem A; the partial order is that of
\eqref{partial}. We denote by $\LA (N)=\LA_N(d,j)\subset
\Grass(d,R_j)$ the scheme parametrizing those  vector spaces $V\subset
R_j$ whose level algebra $\LA (V)$ satisfies $H(\LA (V))=N$; and we let 
$\GA(T)=\GA_T(d,j)\subset \Grass(d,R_j)$ parametrize graded algebras
$R/(V), V\subset R_j$ satisfying
$H(R/(V))=T$. For Theorem E the set
$\mathcal PA(d,j)$ is a certain partially ordered set of pairs of
partitions (Definition \ref{defpos}).
\begin{uthm}[A]
 The proper $O$-sequence $H=(H_0,H_1,\ldots
,H_j,H_{j+1},\ldots )$ as in \eqref{Hcond} occurs as the Hilbert function
of the ancestor algebra of a proper vector subspace of $R_j$ if and only
if the first difference
$E=\Delta (H)$ satisfies the conditions
\begin{align}\label{efirst} e_j&=e_{j+1}\ge e_{j+2}\ge \cdots \ge
e_{\sigma (V)}=0 
\\ \label{esecond} e_j&\ge e_{j-1}\ge e_{j-2}\ge \cdots \ge e_1 \ge e_0=
-1 \text { and }\\
\label{ethird}
\sum_{i\le j}&(e_i+1)+\sum_{i>j}e_i+c_H=j+1.
\end{align}
Each such sequence $E$ satisfying the three conditions occurs, and for a
vector space of dimension $d=\sum_{i\le j}(e_i+1)$.
\end{uthm}
\begin{uthm}[B] Let $d\le j$ be positive integers, and let $H$ be an
acceptable
$O$-sequence. The dimension of $\Grass (H)$ is $c_H+\sum_{i\ge \mu (H)}
(e_i+1)(e_{i+1})$.
\end{uthm}
\begin{uthm}[C]{\sc Frontier property} Assume that $k$ is algebraically
closed. The Zariski closure
$\overline {\Grass (H)}$ is
$\bigcup_{H'\ge_{\mathcal{P}} H} \Grass(H')$.
\end{uthm}
\begin{uthm}[D] Assume that $k$ is algebraically closed. Let $d,j$ be
positive integers satisfying
$d\le j$, and suppose that $H$ is an acceptable $O$-sequence (Definition
\ref{defaccept}).  There is a surjective
morphism $\pi :G(H)\to
\overline{Grass(H)}$ from the nonsingular variety $G(H)$, given by
$I\to I_j$. The inclusion $\iota :\Grass_H(d,j)\subset G(H), \iota :
V\to \overline{V}$ is a dense open immersion. For $H'\in \mathcal
H(d,j), H'\ge_\mathcal P H$, the fibre of
$\pi$ over
$V'\in
\overline{\Grass_H(d,j)}\cap
\Grass_{H'}(d,j)$ parametrizes the family of graded ideals 
\begin{equation*}
\{ I\mid H(R/I)=H \text { and } I_j=V'\}.
\end{equation*}
The schemes $\overline{\LA_N(d,j)}$ and $\overline{\GA_T(d,j)}$
 have desingularizations
$G(N)$ and
$G(T)$, respectively, with analogous properties. 
\end{uthm}
\begin{uthm}[E]   There is an isomorphism $\beta$ from the partially
ordered set $\mathcal H(d,j)$ under the partial order $\mathcal P=\mathcal P(d,j)$,
and the partially ordered set $\mathcal PA(d,j)$ under the product of the
majorization partial orders (see Definition \ref{defpos}). The isomorphism is given
by
$\beta (H)=(P,Q),P=P(H)=A(H)^\ast, Q=Q(H)=B(H)^\ast$ (see Definitions \ref{defP} and
\ref{defA-D}). This is the same order as is
induced by specialization (closure) of the strata $\Grass(H)$.
\end{uthm}
We show similar results to Theorems A-E for the Hilbert function strata
$\LA_N(d,j)$ and $\GA_T(d,j)$. Of these results Theorems C,D --- Theorem
\ref{closureofstrata} in Section \ref{closuresub} --- are the
deepest of the paper. The kind of frontier property shown is rare in this
context of Hilbert schemes of families of ideals. The key step in the
case of $R/(V)$ is the construction of an ideal
$I$ of a given Hilbert function $T=H(R/I)$ such that $I$ contains a given
ideal
$ I'$ of Hilbert function $T'=H(R/I')$, where $T'\ge T$
termwise, and $T,T'$ are permissible Hilbert functions $T=H(R/(V),
T'=H(R/(V')$ for algebras $R/(V)$. This key step is made in
Lemma~\ref{build}, and involves constructing a sequence of intermediate
ideals.\par
 Many of the main results here, including Theorems A-D 
are rewritten from a youthful preprint \cite{I1} of 1975, that was
circulated then, even submitted, but not published, and is hereby
retired! We have chosen to restrict the focus of the present paper to
ancestor algebras, level algebras, and also the algebra
$R/(V)$ determined by $V$, and several applications. We omit the
developing of basic facts about apolarity/Macaulay's inverse systems that
comprised an important part of \cite{I1}, but was both classically
known, and is now well-known in recent literature in the form that we use
in section~\ref{subsimultWaring} (see, for example
\cite{I4,EmI,IK, Ge}). We give here a much-changed
and clearer exposition of Theorems A-D above, and their analogues for
level algebras and the algebras $R/(V)$; the latter case $R/(V)$ was
treated in
\cite[\S 4B]{I2}, but the exposition here is improved.\par
 Several
advances since 1975 have modified our exposition and influenced our
results. The Persistence theorem of Gotzmann, which appeared in 1978,
resolved a natural question that was open at the time of our original
preprint and is a result that had been conjectured by D. Berman
\cite{Be,Go}: see also
\cite{BH,IKl} for further exposition of the persistence and Hilbert
scheme result of G. Gotzmann, a refinement of Grothendieck's construction
of the Hilbert scheme \cite{Gro}. New here is the use of the Gotzmann
results in Section
\ref{HFstrata} to help parametrize the Hilbert function strata of ancestor
ideals, when $r>2$ and  $H$ is not eventually zero.\par
Several authors have written about the restricted tangent bundle to a
rational curve \cite{GhISa,Ra,Ve}, closely related to the Hilbert function
stata $\GA_H(d,j)$. The form of the codimension results there have
inspired  an entirely new Section
\ref{minressub} on the minimal resolutions of the three algebras attached
to $V$. We define partitions $A,B$ giving the generator and
relation degrees of the ancestor ideal $\overline{V}$, and depending only
on the Hilbert function $H(R/\overline{V})$ (Lemma \ref{writeA}); and we
find compact formulas for
the codimensions of $
\Grass_H(d,j), \LA_N(d,j)$ and $ \GA_T(d,j)$ in terms of natural
invariants of these partitions (Theorem \ref{codpartition}). We also count
level algebra and ancestor algebra Hilbert functions using the partitions
(Theorem \ref{allH}, Corollary \ref{count}) and as well we describe the
closures of strata using them (Lemma \ref{comparepart}, Theorem E).  The
Betti strata for more general $O$-sequences $H$ --- not arising from
ancestor algebras --- are studied in a sequel \cite{I6}.\par
The methods of this paper, in particular the proof of the frontier
property of Theorem C for the parameter spaces $\GA_T(d,j)$ of the ideal
$(V)$, can be applied to show a similar frontier property for the
stratification of the family of rational normal curves in ${\mathbb P}^r$
according to the decomposition of the restricted tangent bundle into a
direct sum of line bundles
(see \cite{GhISa}, also \cite{Ra}. The analogous result for
$\LA_N(d,j)$ has a similar interpretation for the stratification of such
a family by the minimal rational scroll upon which they lie \cite{I5}.\par
 In
Section
\ref{subsimultWaring} we apply our results to solve a refined version of
the simultaneous Waring problem for a vector space
$\mathcal W$ of degree-$j$ forms in
$\mathcal R=k[X,Y]$, using apolarity or Macaulay inverse systems. 
The simultaneous Waring problem for a set of $c$ general forms of
specified degrees is to find a smallest integer $\mu$ such that $c$
generic forms of these degrees may be written as linear combinations of
powers of
$\mu$ linear forms. It was studied classically by A. Terracini, whose
approach is generalized and modernized in
\cite{DiFo}. Recently E. Carlini has interpreted the result concerning the
generic (largest) Hilbert function for a level algebra, in
terms of the simultaneous Waring problem, while making explicit the
connection with secant varieties to the rational normal curve \cite{Ca}.
This well known connection of ideals in $k[x,y]$ to secant bundles is
explained in the complete intersection case related to the Waring problem
for a single form in
\cite[Section 1.3]{IK}. Another recent solution of the Waring problem for
forms in two variables occurs in a unpublished preprint with Jacques
Emsalem, a result that can be readily derived from the theory of
compressed algebras 
\cite[Theorem 4.6C]{I4}. In the special case of equal degrees, so one
considers $f\in \mathcal W$, for a general vector space
$\mathcal W\subset \mathcal R_j, r=2$ solutions are given in
\cite[Theorem 3.1]{CaCh},\cite[Theorem 3.3]{Ca}, and \cite[Theorem
3.16]{ChGe}; the latter result also determines the dimension of the
subscheme of $\Grass(c,\mathcal R_j)$ parametrizing vector spaces
$\mathcal W$ having a length $\mu$ simultaneous decomposition. Our
refinement here is two-fold, first to consider vector spaces of degree-$j$
forms $\mathcal W$ having a given \emph{differential
$\tau$ invariant}, and second, we use Theorem
\ref{closureofstrata} to determine the closure of the relevant
$\LA_N(d,j)$ strata (Theorem
\ref{sWaring}).
\par Section \ref{related} has results from the original
preprint
\cite{I1} concerning related vector spaces $V,W$, where $W=R_{i_k}\cdot
R_{i_{k-1}}
\cdots  R_{i_1}V$.
David Berman's article
\cite{Be}  showed that a \emph{complete Hilbert function}
associated to a vector subspace of
$R_j$, ostensibly a function from a countable set of sequences to
$\mathbb{N}$, the non-negative integers giving the dimension of each
space $W$ related to $V$, is determined by its restriction to a finite
subset of the sequences. Here we study primarily the case $r=2$ and we
bound the number of classes $\overline{W}$ related to $V$ (Proposition
\ref{numbrelated}).
\par The results of Sections
\ref{HF2} and \ref{closuresub} in the special case of the algebras
$R/I$ where 
$I=(V)$ when $r=2$ were stated and shown in Proposition 4.7-4.9 and
Theorem 4.10 of
\cite[Section 4B]{I2}. Our exposition here is rather more detailed and
careful even in this special case.
Other results of this article for the case $r=2$ were announced in
\cite[Appendix B]{I3} (the case (V), with an allusion to the ancestor
ideal case), in \cite[Proposition 4.6A,B,C]{I4} (level algebras), in
\cite[Theorem 8.1]{IK} (Gorenstein Artinian algebras), and in a
note on level algebras when $r=2$ at the end of
\cite{ChoI}. But proofs of the results of Sections
\ref{HF2} and \ref{closuresub} for ancestor ideals and level
algebras, when
$r=2$ were in the original preprint
\cite{I1} and appear here for the first time. 
\par
\par Several authors have recently studied level algebras, but from a
rather different viewpoint than taken here
\cite{ChoI,BiGe,Bj1,St}. In addition E. Carlini, and J. Chipalkatti with
Tony Geramita have written about the two variable case, each determining
the possible Hilbert functions for level algebras \cite{Ca,ChGe}.
 E. Carlini and J.~Chipalkatti have made some remarkable progress in
the simultaneous resolution problem in certain other cases
for $r\ge 3$ variables \cite{CaCh}.  J.~Chilpakatti and
A.~Geramita give a geometric description of Hilbert function stratum
$\LA_N(d,j)$ for level algebras in \cite[Propositions 3.7,3.10]{ChGe};
and they draw conclusions for the simultaneous Waring problem for
binary forms (ibid, Theorme 3.16). They also show that certain quite
special unions of these strata are projectively normal, or arithmetically
Cohen-Macaulay (ibid., Theorem 4.4): these unions are different from the
closures $\overline{\LA_N(d,j)}$ studied here.\par
 In higher dimensions $r>2$, until recently only the Gorenstein case $\cod
V = 1$ of level algebras had been extensively studied (see
\cite{IK} for results and references); also a compressed algebra case
where $H$ is maximum given the codimension of $V$ and $r$ had been studied
\cite{I4,FL,Bj1}. The analogue for
$r>2$ of the frontier property of Theorem C does not usually hold even in
the Gorenstein height three case
\cite[Example 7.13]{IK}, nor is $G(H)$ a desingularization of $\Grass(H)$
\cite[Lemma 8.3 with J. Yam\'{e}ogo]{IK}. The sequences $H$ that occur as
Hilbert functions $H=H(R/\overline{V})$ are known when $r=3$ in the
Gorenstein case \cite{BE,St,D}, (see \cite[\S 5.3.1]{IK}); also in
this Gorenstein case the family $\Grass(H)$ is irreducible and
nonsingular
\cite{D,Klp}. The question of which sequences $H$ occur as Hilbert
functions of level algebras $LA(V)$ is studied by A. Geramita, T. Harima,
and Y. Shin in
\cite{GHS1} using \emph{skew configurations} of points in $\mathbb P^n $.
 With J. Migliore they develop further results, including
necessary conditions  and new techniques and constructions for arbitrary
socle degree and type; they also include a complete list
of level Hilbert functions for $r=3$, socle degree at most 5, of of socle
degree 6 and type
$\cod V=2$ \cite{GHMS}. When
$r\ge 4$ even the set of Gorenstein sequences are unknown. However,
several authors have established both minimum and maximum Hilbert
functions for level algebrs
$\LA(d,j)$ in any codimension $r$ (see
\cite{BiGe,ChoI}).\par
\subsection{The Hilbert function
strata}\label{HFstrata} Fix
$r$ and the polynomial ring $R=k[x_1,\ldots ,x_r]$. Recall that we denote
by $\Grass (d,R_j)$ the Grassmanian parametrizing $d$-dimensional vector
subspaces of $R_j$. A reader primarily interested in $r=2$ may wish to
skip over or skim this section and consult Proposition \ref{GCD} in its
place. 
\begin{definition} Let $H$ be a sequence of non-negative
integers that occurs as the Hilbert function $H=H(R/\overline{V})$ where
$ V$ is a
$d$-dimensional vector subspace of
$R_j$. We denote
by
$\Grass_H(d,j)\subset \Grass (d,R_j)$ the subscheme of the Grassmanian
parametrizing vector spaces $V$ satisfying the rank
conditions
\begin{equation}\label{estrata} \cod R_iV=H_{i+j}
\text { in } R_{i+j}, 
\text { for } i=-j,-j+1,\ldots .
\end{equation}
\end{definition}
When $H$ is eventually zero, evidently equation \eqref{estrata} imposes a
finite number of algebraic conditions on $V$ (which we study shortly). 
When $H$ is not eventually zero, we will use Gotzmann's Persistence and
Hilbert scheme theorems, a refinement of the Grothendieck Hilbert scheme
theorem, to show that the number of algebraic conditions imposed by 
\eqref{estrata} is
finite. \par Recall that every sequence $H=(H_0,\ldots )$ occurring as the
Hilbert function $H=H(A)$ of a quotient algebra $A=R/I$ is eventually
polynomial: there exists a pair $(p_H\in \mathbb{Q}[t],s=s(H) \in {\mathbb
N} )\mid H_i=p_H(i)\text{ for }i\ge s(H)$. We denote by $\sigma=\sigma
(p_H)$ the Gotzmann regularity degree of $p_H$ (see
\cite{Go,IKl}). It is easy to see that $\sigma \ge s(H)$.
Recall that the Grothendieck Hilbert scheme
$\Hilb^p({\mathbb P}^{r-1})$ parametrizes subschemes of ${\mathbb P}^{r-1}$
having Hilbert polynomial $p$ \cite{Gro}. We denote by $r_i$ the integer
$r_i=\dim_k R_i=\binom{r+i-1}{i}$, and define $q=q_H$ by
$q(i)=r_i-p_H(i)$. We denote by $M(d,j)$ the vector space span of the
first
$d$ monomials of degree $j$ in $R$, in lexicographic order.
\begin{theorem}{\sc Macaulay Growth Theorem}\cite{Mac2}\label{Macgrowth} A
vector space
$V\in
\Grass (d,R_j)$ satisfies
\begin{equation}\label{Maceq}
\dim R_1\cdot V\ge \dim R_1\cdot M(d,j).
\end{equation}
\end{theorem}
\begin{theorem}{\sc Gotzmann 
Hilbert scheme and Persistence Theorem}\label{Gotzper}\cite{Go} Let $p$
be a Hilbert polynomial, and $\sigma=\sigma(p)$. The Hilbert scheme
$\Hilb^p({\mathbb P}^{r-1})$ is the locus of pairs of vector spaces
\begin{equation}\label{eGotzper}
(V,V')\in\Grass(q(\sigma ),R_\sigma )\times\Grass(q(\sigma
+1),R_{\sigma +1})
\end{equation}
satisfying $R_1\cdot V=V'$, or, equivalently $R_1\cdot V\subset V'$. Such
vector spaces $V$ satisfy equality in \eqref{Maceq}.
\par \noindent
{\sc (Persistence)} A vector space $V$ occurring in such an extremal
growth pair $(V,V')$ satisfies 
\begin{equation}
\dim (R_{\sigma +i}/R_iV)=p(\sigma +i) \,\, \forall i\ge 0;
\end{equation}
the space $R_iV$ has dimension $q(\sigma +i)$, and also satisfies equality
in
\eqref{Maceq}.
\end{theorem}
For an exposition of the persistence result over $k$, see \cite[\S
4.3]{BH}; for an exposition of the Gotzmann-Grothendieck Hilbert scheme
results and further references see
\cite{IKl}. One consequence of Theorem \ref{Gotzper} for us is that
one may suppose that $i\le \max\{1,\sigma p_H+1-j\}$ in equation
\eqref{estrata}. Thus \eqref{estrata} defines a scheme structure on
$\Grass_H(d,j)$ as locally closed subscheme of $\Grass (d,R_j)$, for all
occurring sequences $H$.\par
 Given such a sequence
$H$ we define a projective scheme
$G(H)$ parametrizing the graded ideals
$I\subset R$ that determine a quotient algebra $A=R/I$ having Hilbert
function
$H(A)=H$. When
$H$ is eventually zero, so $H_s=0$, the parametrization of $G(H)$ is as a
subset of 
$\prod_{i\le s}\Grass(r_j-h_j, R_j)$, where $r_j=\dim_k R_j$. When $H$ is
not eventually zero, then $H$ is eventually polynomial $H_i=p_H(i)$ for
$i\ge s(H)$ for some polynomial $p=p_H$. As before, we take
$\sigma (H)$ the regularity degree of the polynomial, and parametrize\par
\begin{equation}\label{eGH} G(H)\subset \left(\prod_{i<\sigma}
\Grass(r_j-h_j,R_j)\right)\times
\Hilb^p({\mathbb P}^{r-1}).
\end{equation}
By Theorem \ref{Gotzper}, we may replace the product in equation
\eqref{eGH} by  $\prod_{i\le \sigma +1}\Grass(r_j-h_j, R)$.\par
Results of D. Mall (when $\mathrm{char} k=0$ or
$\mathrm{char} k>
\sigma (p_H)$
and K. Pardue (for arbitrary characteristic) show that when the base
field $k$ is algebraically closed, the scheme $G(H)$ is connected
\cite{Mall,Par}.
\begin{definition}\label{defpo}
We define a partial order $\mathcal P=\mathcal { P}(d,j,r)$ on the
set
$\mathcal {H}(d,j,r)$ of Hilbert functions possible for
$H(A),A=R/\overline{V}$, as follows:
\begin{equation}\label{po}
H'\ge_{\mathcal{P}(d,j,r)}H\Leftrightarrow H'_i\le H_i \text { for } i\le
j
\text { and } H'_i\ge H_i \text { for } i\ge j.
\end{equation}
\end{definition}
When the triple $(d,j,r)$ is obvious from context we write $H'\ge
_{\mathcal P}H$ for
$H'\ge_{\mathcal P(d,j,r)}H$. Recall that $H$ \emph{occurs} or is
possible for us if it occurs as the Hilbert function of an ancestor
algebra $\Anc (V)$ for some $d-$ dimensional vector subspace of $R_j$.
\begin{theorem}\label{loclostrata} Let $H$ be a sequence that occurs
as the Hilbert function of an ancestor algebra.  The scheme
$\Grass_H(d,j)$ is a locally closed subscheme of $\Grass (d,R_j)$.
The condition $H'=H(R/\overline{V})\ge_{\mathcal P} H$ is a closed
condition on
$V\in\Grass (d,R_j)$. Also the inclusion
$\iota: \Grass_H(d,j)
\to G(H)$ given by $ \iota : V\to \overline{V}$ is an open immersion.
\end{theorem}
\begin{proof} Let $I=I_V=\overline{V}$. It is not hard to show that $\dim
I_i\ge r_i-H_i$ is a closed condition, and $\dim I_i< r_i-H_i+1$ is an
open condition on $V\in \Grass (d,R_j)$, when $i\le j$. Likewise, it is
not hard to show that for each $i\ge j$ then $\dim I_i\le r_i-H_i$ is a
closed condition, while $\dim I_i>r_i-H_i-1$ is an open condition. By the
Gotzmann persistence and regularity theorems, if $V$ satisfies each of
these conditions for all positive integers $i\le \sigma (p_H)+1$ (which we
may suppose greater than $j$), then
$H(R/\overline{V})=H$. Thus, we have shown that $\Grass_H(d,j)\subset
\Grass (d,R_j)$ is defined by the intersection of a finite number of open
and closed conditions, so it is locally closed, as claimed. \par
That the inclusion $\iota$  is an open immersion, follows from $I_{\ge j}$
being generated by
$I_j$, and $I_i, i<j$ being  $R_{i-j}I_j$. For $a>0$ the condition that
$V=I_j$ generates
$I_{j+a}$ is equivalent to the rank of the multiplication map:
$R_a\otimes V\to R_i$ 
 being greater than $\dim I_i-1=r_i-H_i-1$ on $G(H)$  --- an open
condition. Let $W=V^\perp \subset \mathcal R_j$ in the Macaulay
duality. For $a>0$ the condition that
$I_{j-a}=R_{-a}V$ is equivalent to the rank of the contraction map
$R_a\times W\to R_a\circ W \subset \mathcal R_{j-a}$ being greater than
$H_i-1$, on
$G(H)$, also an open condition. This completes the proof.
\end{proof}
\begin{corollary}\label{closure}
The Zariski closure
$\overline{\Grass_H(d,j)}\subset
\bigcup_{H'\ge_{\mathcal P} H}\Grass_{H'}(d,j)$. Similar inclusions hold
for
$\overline{\LA_N(d,j)}$ and for $\overline {\GA_T(d,j)}$.
\end{corollary}
\begin{remark} The partial order $\mathcal P(d,j,r)$ for $r\ge 2$ is not
in general subordinate to or equal to a simple order. For $r=2$ a simply
ordered exception are the complete intersection cases
$(d,j)=(d,d+1)$,  where $V$ has
codimension one: see \cite[\S
1.3]{IK}. Also for $r=2$, Example
\ref{exsimple} gives a different simply ordered case, $(d,j)=(4,5)$,
while Example
\ref{nsimple}(A) below $(d,j)=(3,5)$ and Example
\ref{nsimple}(B) $(d,j)=(10,12)$ illustrate the more
general situation $\mathcal P(d,j,2)$ not a simple order, for ancestor
algebras and level algebras, respectively.
\end{remark}

\section{The ancestor ideal in two variables}\label{anc2}
 Throughout this section, $R$ is the polynomial ring $R=k[x,y]$ over
an arbitrary field $k$, and we denote by $M=(x,y)$ the homogeneous
 maximal ideal. The vector space
$R_j$ of degree-$j$ forms in
$R$ satisfies, $R_j=\langle x^j,x^{j-1}y,\ldots ,y^j\rangle $, of
dimension
$j+1$, and $V\subset R_j$ will be a vector subspace having dimension $\dim
V=d$. In Section \ref{HF2} we give our main results concerning the
individual Hilbert function strata of the three algebras related to $V$
when $r=2$. These include a characterization of ancestor ideals
(Proposition \ref{ancideals}) and the dimension/structure Theorem
\ref{dimstrata}. In Section \ref{minressub} we give our results relating
the graded Betti numbers of these three algebras to certain partitions
$A,B,C,D$ (Lemma \ref{writeA}); also we give the codimension of the
Hilbert function strata in terms of the partitions $A,B$ or $C,D$
(Theorem \ref{codpartition}). In Section \ref{closuresub} we determine
the closures of the Hilbert function strata (Theorem
\ref{closureofstrata}).
\subsection{The Hilbert function strata when $r=2$}\label{HF2}
We first present the main tool we need, the simplicity
$\tau (V)$, and a key exact sequence.
\begin{definition}\label{taudef} For $V\subset R_j$ we define 
\begin{equation}\label{etaudef}
\tau (V)=\dim_k R_1
V-\dim_k V.
\end{equation}
 We define the sequence
\begin{equation}\label{tauseq}
   0\to R_{-1}V\xrightarrow{\phi} R_1\otimes V\xrightarrow{\theta} 
R_1\cdot V
\to 0,
\end{equation}
where $\phi : f\to y\otimes xf-x\otimes yf$, and $\theta
:\sum_i\ell_i\otimes v_i\to \sum_i \ell_iv_i$, where the $\ell_i$ are
elements of $R_1$ (linear forms). 
\par
\end{definition}
For $I$ a graded ideal of $R$, we denote by $\nu (I)$ the
number of minimal generators for $I$. For a vector subspace $W\subset
R_i$ we denote by $\cod W=i+1-\dim W$, the codimension of $W$ in $R_i$.
\begin{lemma}\label{taugen} The sequence \eqref{tauseq} is exact. We have
\begin{align} \tau (V)&=\dim V-\dim
R_{-1}V\label{tauseq1}\\
&=1+\cod R_{-1}V-\cod V=1+\cod V-\cod
R_1V\label{tauseq2}\\
 &=\nu (\overline {V}).\label{tauseq3}
\end{align}
Also, $\tau (V)\le \min\{ d,j+2-d\}$.
\end{lemma}
\begin{proof} Clearly $\phi$ is a monomorphism, and $\theta$ is
surjective, so we need only show the exactness of \eqref{tauseq} in the
middle. Suppose that $U\in R_1\otimes V$ and $\theta (U)=0$. We may
suppose $U=x\otimes v_1+y\otimes v_2$, thus $xv_1+yv_2=0$, implying
$y$ divides $v_1$ and $x$ divides $v_2$. Thus $w=v_2/x=-v_1/y\in R_{-1}V$
satisfies 
\begin{equation}\label{etaudim}
\phi (w)=y\otimes xw-x\otimes yw=y\otimes v_2-x\otimes (-v_1)=U.
\end{equation}
This completes the proof of the exactness of \eqref{tauseq}. Thus,
counting dimensions in \eqref{tauseq} we have
\begin{equation}\label{etaudefi}
2\dim V=\dim R_1\otimes V=\dim R_{-1}V+\dim R_1V.
\end{equation}
Noting the definition of
$\tau$ in \eqref{etaudef}, we have shown
\eqref{tauseq1}. The equations
\eqref{tauseq2} follow immediately. To show that $\tau (V)=\nu
(\overline{V})$, we first note that applying \eqref{etaudefi} to $R_iV$ we
have for any integer $i$ satisfying $ -j\le i$,
\begin{equation}
\dim R_{-1}R_iV+\dim R_1R_iV=2\dim R_iV.
\end{equation}
When $i\le 0$ we have $R_{-1}R_iV=R_{i-1}V$, so we have
\begin{equation}\label{taus}
\text{ for  } i\le 0 \quad \dim R_1R_iV=2\dim R_iV-\dim R_{i-1}V.
\end{equation}
The number of generators $\nu (\overline{V})$ of the ancestor ideal
of $V$ satisfies, $\nu (\overline{V})=\dim_k
(\overline{V}/M \overline{V})$, where $M\overline{V}=R_1\overline{V}$,
since
$\overline{V}$ is graded.
 We
have
\begin{align}
\overline{V}/R_1\overline{V}&=\oplus_{i=-j}^{+\infty}
(R_iV/R_1R_{i-1}V)\notag
\\ &=\oplus_{i=-j}^{0} (R_iV/R_1R_{i-1}V),
\label{taustep}
\end{align}
since for $i\ge 0$ we have $ R_1R_{i-1}V=R_iV$. Let $d_i=\dim R_i V$.
From \eqref{taustep} we have
\begin{align*}
\nu (\overline{V})&=\sum_{i=-j}^0\dim R_iV-\sum_{i=-j}^0\dim R_1R_{i-1}V\\
&=\sum_{i=-j}^0d_i-\left( 2\sum_{i=-j}^0
d_{i-1}-\sum_{i=-j}^0d_{i-2}\right)
\text{ by \eqref{taus} }\\
&=d_0-d_{-1}\\
&=\tau({V})  \text{   \quad by \eqref{tauseq1} }.
\end{align*}
This completes the proof of \eqref{tauseq3}. The upper bound on $\tau
(V)$ is immediate from \eqref{tauseq1} and \eqref{tauseq2}. 
\end{proof}\par
Recall from Definition \ref{def1.6} that the subspace $V\subset
R_j$ is
\emph{equivalent} to $W\subset R_i$ if $\overline{V}=\overline{W}$. A
generalization of (iii) below is shown in Corollary \ref{reltauequal}.
\begin{proposition}{\sc{Equivalence.}}\label{taubasic} We assume that
$V\subset R_j$; here $R=k[x,y]$. 
\begin{enumerate}[{i.}]
\item\label{taubasici} For $s\ge -j$ we have
$\tau
(R_sV)\le
\tau (V)$, with equality if and only if $\, \overline{R_sV}=\overline{V}$.
\item\label{taubasicii} In the
sequence
\begin{equation*}
\tau (R_{-j}V), \ldots , \tau (R_{-1}V),\tau (V),\tau(R_1V),\ldots
\end{equation*}
the values of $\tau(R_iV) $ are monotone non-decreasing for $i\le 0$, and
monotone non-increasing for $i\ge 0$. 
\item\label{taubasiciii} For two vector spaces $R_sV, R_tV$, we have
\begin{align*}
\overline{R_sV}=\overline {R_tV}&\Leftrightarrow
R_sV=R_{s-t}R_tV\text{ and } R_tV=R_{t-s}R_sV.\\
&\Leftrightarrow 
\begin{cases} \text { either }\tau (R_sV)=\tau (R_tV)&=\tau(V),\\\text {
or }
\mathrm {sign} (s)=\mathrm {sign} (t) \text { and }&\tau (R_sV)=\tau
(R_tV) .
\end{cases}
\end{align*}
\item\label{taubasiciv}
\begin{equation*}
 \overline{R_sV}=\overline{V}\Leftrightarrow
\begin{cases} \text{ If $s>0$},\, & \dim R_{s+1}V=\dim V+(1+s)\tau (V);\\
 \text{ If $s\le 0$},\, & \dim R_{s-1}V=\dim V -(1-s)\tau (V).
\end{cases}
\end{equation*}
\item\label{taubasicv} For any two vector spaces $V\subset R_j, W\subset
R_i$,
\begin{equation*}
\overline {V}=\overline W \Leftrightarrow\, V=R_{j-i}W \text { and  \,
}\tau (V)=\tau (W).
\end{equation*}
\item\label{taubasicvi} $\tau (V)=1\Leftrightarrow V=f\cdot R_{j-c}$ where
$\deg f=c=\cod V$. Also $\tau (V)=0\Leftrightarrow V=0$.
\end{enumerate}
\end{proposition}
\begin{proof} To show \eqref{taubasici} it suffices to prove it for $s=\pm
1$ and apply an induction. For $s=1$ we have $\tau (R_1V)=\dim R_1V-\dim
R_{-1}R_1V$, but
$R_{-1}R_1V\supset V$, so $\tau (R_1V)\le \dim R_1V-\dim V=\tau(V)$ with
equality if and only if $R_{-1}R_1V=V$, which is equivalent to
$\overline{V}=\overline{R_1V}$. For $s=-1$, we have $\tau (R_{-1}V=\dim
R_1R_{-1}V-\dim R_{-1}V\le \dim V-\dim R_{-1}V=\tau(V)$ with equality if and only if $R_1R_{-1}V=V$, which is equivalent to
$\overline{R_{-1}V}=\overline{V}$. \par
Repeated use of \eqref{taubasici} shows the rest of the Proposition. For
example, we show \eqref{taubasiciv} for $s>0$. 
By definition $\tau(R_iV)=\dim R_{i+1}V-\dim V$ for $i=0,\ldots
,s$ so we have for $W=R_sV$,
\begin{equation*}
\dim R_1W=\dim V +\tau(W)+\tau(R_1V)+\cdots +\tau(R_sV).
\end{equation*}
That $\tau (V),\tau (R_1V), \ldots $ is nonincreasing shows that
$\dim R_1W=\dim V+(s+1)\tau(V)\Leftrightarrow \tau(V)=\tau(R_1V)=\cdots
=\tau (R_sV)$, as claimed. This completes the proof of \eqref{taubasiciv}.
For \eqref{taubasicvi}, evidently $\tau (V)=0\Leftrightarrow V=0$. When
$\tau (V)=1$, then Lemma
$\overline{V}=(f)$ by Lemma \ref{taugen}. Letting $c=\deg f$
 we thus have $R_{c-j}V=\langle f\rangle$ and
$R_{j-c}f=\overline{V}_j=V$, whence $c=\cod V$, as claimed. This
completes the proof of \eqref{taubasicvi}.
\end{proof}\par
\begin{example}We show here the need to use the $\dim (R_{s+1}V)$ in
Proposition
\ref{taubasic}\eqref{taubasiciv} to decide if $R_sV$ is equivalent to
$V$ when
$s>0$, and the need for
$R_{s-1}V$ when
$s\le 0$.  Let $V=\langle x^4,x^3y,y^4\rangle\subset R_4
$, then
$R_{-1}V=\langle x^3\rangle$, and $\overline{V}=(x^3,y^4),$ so
$\tau(V)=2$ while
$\overline{R_{-1}V}=(x^3)$, yet we have $\dim R_{-1}V=\dim (V)-\tau(V)$.
Thus, the dimension of $W=R_sV$ is not enough to test the equivalence of
$W$ and $V$. Here
$\dim R_{-2}V=0\not=
\dim V-2\tau(V)$, corresponding to
$\overline{V}\ne \overline{R_{-1}V}$. 	Here
$\overline{V}=\overline{R_1V}$, and $\dim R_1V=5=\dim V+\tau(V)$, $\dim
R_2V=\dim V+2\tau(V)$, but $R_2V=R_6$ so
$\overline{V}\not=\overline{R_2V}$. Here $j=4$, $\overline{V}$ is a
complete intersection, satisfying
$H(\Anc (V))=(1,2,3,3,2,1)$, $E(H)=\Delta H=(-1,-1,-1,0,e_4=1,1,1)$. As
in Proposition \ref{etau} \eqref{ehilbineq} the subsequence
$(-1,-1,-1,0,1=e_4)$ of
$E(H)$ is non-decreasing, while the subsequence ($1=e_4,1,1)$ is
non-increasing, and
$\tau (V)=2=e_4+1=e_5+1$ (see Proposition \ref{etau} \eqref{ehilbineqii}).
\end{example}
We define the greatest common divisor $\GCD (V)$ as the
principal ideal in $k[x,y]$ with a generator of highest degree,
such that $\GCD (V)$ contains $V$ (the generator divides each element of
$V$).  We will now show directly for $R=k[x,y]$ that
$\lim _{i\to
\infty}\overline{R_iV}=\GCD (V)$, a special case of $\lim _{i\to
\infty}\overline{R_iV}=\Sat(V)$ in Lemma~\ref{satV}.
\begin{proposition}\label{GCD} Assume that $H=H(R/\overline{V})$ satisfies
$\lim_{i\to
\infty}H_i=c$. Then we have 
\begin{align}\label{sumtautail}
\sum_{i\ge 0}\left( \tau (R_iV)-1\right)& =\cod
V -c=(j+1-d)-c,\\
\sum_{i\le 0} \tau (R_i\cdot V)&=\dim V = d.\label{sumtaunose}
\end{align}
 The degree $\deg \GCD (V) =c$. For $i\ge \cod V-\tau (V)+2$, we have 
\begin{equation}\tau (R_i\cdot V)=1 \text{ and
}\overline{R_i\cdot V}=\GCD (V)
\end{equation}
\end{proposition}
\begin{proof} Let $k\ge 0$ satisfy $H_{k+j}=c$; then evidently $\tau
(R_k\cdot V)=1$ and by Proposition \ref{taubasic}\eqref{taubasicii} we
have
$c=\deg \GCD (R_k\cdot V)$ and evidently since $k\ge 0$, we have $\GCD
(R_k\cdot V)=\GCD ( V)$. Now, equation
\eqref{sumtautail} is a consequence of \eqref{tauseq2}, and Equation
\eqref{sumtaunose} follows from
\eqref{tauseq1}. We now turn to the explicit bound on $i$ for
achieving $\tau (R_i\cdot V)=1$. Suppose on the contrary that for
an integer $i\ge 2$ we have
$\tau (R_i\cdot V)\ge 2$.  Proposition
\ref{taubasic}\eqref{taubasicii} shows that the sequence
$\tau (V),\tau (R_1\cdot V),\ldots $ is montone, hence we have from
\eqref{sumtautail},
\begin{equation*}
 \tau (V)-1+i\le (\tau(V)-1)+(\tau(R_1\cdot V)-1)+\cdots (\tau
(R_i\cdot V)-1)\le
\cod V,
\end{equation*}
implying $i\le \cod V-(\tau (V)-1)$. Thus we have the explicit bound $\tau
(R_iV)=1$ for
$i\ge \cod V-\tau (V)+2$, as claimed. By Lemma \ref{taugen} we have
for such $i,\, \overline{R_i\cdot V}=(f)$. As above we conclude
by Proposition
\ref{taubasic}\eqref{taubasicvi} that for such $i$, we have $f=\GCD (
R_i\cdot V)=\GCD ( V)$.

\end{proof}\par
Recall that when $H=H(R/I)$ is the Hilbert function of a graded quotient
of $R$, we denote by $E( H)$ the first difference sequence $E(H)=\Delta
H=(e_0=-1,e_1,\ldots ,e_i,\ldots )$ where $e_i=(\Delta
H)_i=H_{i-1}-H_{i}$. We set $\mu (H)=\min\{i\mid H_i<i+1\}$, which is the
order of any ideal $I\subset R$ with $H(R/I)=H$. Recall that since $H$ is
an $O$-sequence with $H_1\le 2$, $H$ must satisfy
\eqref{Hcond}, so $0\le H_i\le i+1$, and for $I_i\ne 0$, $H_{i+1}\le H_i$.
Thus, $H\ne H(R)$ (or $I\ne 0$) implies $\lim_{i \to \infty}H_i= c_H\ge 0$
with $c_H$ a non-negative constant. When $H=H(R/\overline{V})$ we have by
Proposition \ref{GCD}, $c_H=\deg \GCD (V)$.
\begin{proposition}\label{etau} 
Let $V\subset R_j$ be a vector subspace satisfying $\dim V=d$, and let
$H=H(R/\overline{V})$ as above be the Hilbert function of the ancestor
algebra of $V$, and let $c=c_H$. The first difference sequence
$E(H)$ satisfies
\begin{align}\label{ehilbineq}
&e_i\le e_{i+1} \text { for } i\le j,\text{ and } e_i\ge e_{i+1} \text{
for }
 i\ge j;\\ 
\text { also }& \sum_{i\le j} (e_i+1)=d\,  \text {  and } \sum_{i>j}e_i=
(j+1-d)-c. \label{ehilbineqi}
\end{align}
Let $V\subset R_j$ and let $H=H(R/\overline{V})$. Then
$\tau(R_{i-j}\cdot V)$ satisfies
\begin{equation}\tau (R_{i-j}\cdot V)=
\begin{cases}e_i+1=\nu(\overline{R_{i-j}\cdot V})=\#\{\text{generators of
$\overline{V}$ of degree }\le i\} & \text { if $i\le j$}\\
e_{i+1}+1 & \text { if $i\ge j$ }. 
\end{cases}\label{eqetau}
\end{equation}
We have $e_j=\tau (V)-1$ and
\begin{equation}\label{ehilbineqii}
0\le e_j=e_{j+1}\le \min\{j+1-d,d-1\},
\end{equation}
with equality $e_j=d-1$ if and only if $R_{-1}V=0$. Also, $e_{j+1}=\cod V$
if and only if $R_1V=R_{j+1}$.\par
\end{proposition}
\begin{proof} By
applying the first part of equation
\eqref{tauseq2} to $R_{i-j}\cdot V$ when $i<j$, we obtain  
\begin{equation*}
\tau
(R_{i-j}\cdot V)=\cod R_{i-j-1}\cdot V-\cod R_{i-j}\cdot V +1=e_i+1
\end{equation*}
 which is the first part of equation
\eqref{eqetau}. For any $i$ we have by Lemma 2.2 
$\tau (R_{i-j}\cdot V)=\nu (\overline{R_{i-j}\cdot V})$; when $i\le j$ we
have also the second part of equation \eqref{eqetau} since
\begin{align*}
\nu (\overline{R_{i-j}\cdot V})&=\sum_{u\le i}\left( \dim
R_{u-j}\cdot V-\dim R_1\cdot R_{u-j-1}\cdot V\right) \\
&=  \# \{ \text{generators of $\overline{V}$ having degree $\le
i$ }\} .
\end{align*}
 By applying the second part of
equation
\eqref{tauseq2} to $R_{i-j}\cdot V$ when $i\ge j$ we obtain 
\begin{equation*}
\tau
(R_{i-j}\cdot V)=\cod R_{i-j}\cdot V-\cod R_{i-j+1}\cdot V +1=e_{i+1}+1,
\end{equation*}
 which is the last part of
equation
\eqref{eqetau}. Equation \eqref{ehilbineq} now follows from Proposition
\ref{taubasic}\eqref{taubasicii}, and equation \eqref{ehilbineqi},
follows from the definition of $E(H)$ as a first
difference of $H$. The equation
\eqref{ehilbineqii} and remaining claims follow from \eqref{eqetau}.
\end{proof}\par\noindent
\begin{definition}\label{defaccept} Let $d,j$ be positive integers
satisfying $d\le j$. We say that a proper
$O$-sequence
$H$ (a sequence $H$ satisfying \eqref{Hcond}) is \emph{acceptable} for an
ancestor algebra in two variables of a $d$-dimensional subspace of $R_j$
if
$H$  satisfies
\eqref{ehilbineq},
\eqref{ehilbineqi}, and \eqref{ehilbineqii} of Proposition~\ref{etau}. 
\end{definition}\noindent
The sequence $H=0$ occurs for $V=R_j$, and $H=H(R)=(1,2,\ldots )$ occurs
for $V=0$, but we will omit these cases henceforth.
\begin{corollary}\label{accept}
 Let $j$ be a positive integer. A proper
$O$-sequence
$H$ of \eqref{Hcond} is acceptable for an ancestor ideal of a degree-$j$
vector space iff the first difference $E=\Delta (H)$ satisfies
\begin{align}\label{edfirst} e_j&=e_{j+1}\ge e_{j+2}\ge \cdots \ge
e_{\sigma (V)}=0
\\\label{edsecond} e_j&\ge e_{j-1}\ge e_{j-2}\ge \cdots \ge e_1 \ge e_0=
-1, \text { and }\\\label{edthird}
\sum_{i\le j} &(e_i+1)+\sum_{i>j}e_i+c_H=j+1.
\end{align}
\end{corollary}
\begin{proof} Immediate from Definition \ref{accept}, and
\eqref{ehilbineq},
\eqref{ehilbineqi}, and \eqref{ehilbineqii}. Here $d=\sum_{i\le
j}(e_i+1)$.
\end{proof}\par
 In the following definition we use \emph{partition} of
$n$ in the usual sense of 
$n=n_1+n_2+\cdots +n_u, n_1\ge n_2\ge \cdots \ge n_u >0$. Part of the reason for
our choice of $P,Q$ is that we later show they are the duals of the pair
of partitions
$(A,B)$ determined by the generator degrees, and the relation degrees of
ancestor algebras
$\Anc (V)$ satisfying $H(\Anc (V))=H$ (Lemma
\ref{writeA}). Recall that the order $\mu (H)$ of an $O$-sequence is the
smallest integer such that $H_i\ne i+1$. We let $s(H)=\min\{ i\mid
H_i=c(H)\}$. Also given $j,H$, with $H$ acceptable, we define
$\tau (H)=H_{j+1}-H_j+1=e_{j+1}+1=e_j+1$.
\begin{definition}\label{defP} Given positive integers $d,j$ with $d\le
j$ and an acceptable
$O$-sequence
$H$ as in Definition \ref{defaccept}, and letting $\tau =\tau
(H)=e_j(H)+1$, we define a pair of partitions $(P=P(H),Q=Q(H))$ of
$(d,j+1-d-c(H))$ as follows. Let $V$ satisfy $H(R/\overline{V})=H$. Then
$P(H), Q(H)$ satisfy
\begin{align}\label{edefP}
P(H)&=(\tau,\tau (R_{-1}\cdot
V)=e_{j-1}(H)+1,\tau(R_{-2}\cdot V)=e_{j-2}(H)+1,\ldots ,e_\mu (H)+1,\\
Q(H)&=(\tau-1=e_{j+1}(H),e_{j+2}(H),e_{j+3}(H),\ldots
,e_s (H)).\label{edefQ}
\end{align}
Recall from
Definition
\ref{defpo} that $\mathcal H(d,j,2)$ is the set of sequences
possible for the Hilbert function of $\Anc (V), V$ a
$d$-dimensional subspace of $R_j, R=k[x,y]$; understanding that $r=2$ we
will denote this set by $\mathcal H(d,j)$. We will likewise denote by
$\mathcal P(d,j)$ the partial order $\mathcal P(d,j,2)$ on $\mathcal
H(d,j,2)$ from Definition
\ref{defpo}.
We will denote by
$\mathcal H(d,j)_\tau$ the subset of $\mathcal H(d,j)$ for which
$e_j=\tau -1$.
\end{definition}
We will shortly show that the $O$-sequences that are acceptable in the
sense of Definition \ref{defaccept} are exactly those that occur as the
Hilbert function of an ancestor algebra (Theorem 
\ref{allH}).  So each pair $(P,Q)$ of partitions described
in the Lemma below actually occurs as $P=P(H), Q=Q(H)$ for some acceptable
$H$. 
\begin{lemma}\label{partitions} For
\eqref{partitionsi}-\eqref{partitionsiii} below we suppose that the
$O$-sequence
$H$ is proper and acceptable, as in Definition \ref{defaccept}, and let
$\tau=\tau (H)$. Then
\begin{enumerate}[i.]
\item\label{partitionsi} The partition $P=P(H)$ of Definition \ref{defP}
is a partition of $d$ having largest part $\tau $. The partition $Q=Q(H)$
is a partition of $j+1-d-c$ having largest part $\tau -1$.
\item\label{partitionsiii} Let $(\mu (H), s (H))=(\mu, s )$.
Then
$P(H)$ has
$j+1-\mu$ parts, and  $Q(H)$ has $s -j$ parts.\par
\item\label{partitionsii} $H$ is uniquely determined by $(j,
P(H),Q(H))$. 
\item\label{partitionsiv} Let $d,j$ be positive integers, with $d\le j$
There is a one-to-one onto correspondence
$H\to (P(H),Q(H))$ between the subset of acceptable
$O$-sequences $H$ satisfying $(\mu
(H),s (H))=(\mu,s)$ and
$c(H)=c$, and the set of pairs of partitions
$(P,Q)$ satisfying
\eqref{partitionsi}
and
\eqref{partitionsiii}. There are similar one-to-one correspondences
between the set of partitions 
$P$ and the set of sequences $N=N_H$, and also between the set of
partitions
$Q$ and the set of sequences $T=T_H$ (Definiton
\ref{defNT}).
\end{enumerate}
\end{lemma}
\begin{proof} The claim in \eqref{partitionsi} that $P$ partitions $d$
is \eqref{sumtaunose}. That 
the parts of $P$ are less than $\tau$ follows from Proposition
\ref{taubasic}\eqref{taubasicii}. That $Q$ partitions
$j+1-d-c$ follows from \eqref{ehilbineqi}; that
$e_{j+1}=\tau -1$ is \eqref{ehilbineqii}. That the parts of Q are
no greater than $\tau-1$ follows as before from Proposition
\ref{taubasic}\eqref{taubasicii}. The claim of
\eqref{partitionsiii} is immediate from the definitions, counting the
nonzero parts of $P,Q$. For \eqref{partitionsii}, we note that
the triple $(P,Q,j)$ determines $(P,Q,\tau )$ so determines $E(H)$, and
also $d,j$, hence
$c=c(H)$; then $H_i=c+\sum_{i< k} e_k$ determines
$H$. The proof of
\eqref{partitionsiv} is also immediate.
\end{proof}\par
The following Proposition and Corollary describe which ideals are ancestor
ideals, in terms of the degrees of the generators and relations. In a
related result, we
 determine the graded Betti numbers of the ancestor algebra $\Anc (V)$ in
terms of the Hilbert function
$H(\Ann(V))$ (Lemma \ref{writeA}).
\begin{proposition}{\sc Ancestor ideals.}\label{ancideals} Let $I$ be a
graded ideal of
$R=k[x,y]$. The following are equivalent:
\begin{enumerate}[i.]
\item\label{ancidealsi} $I$ is the ancestor ideal of $I_j$.\par
\item\label{ancidealsii} $I$ is homogeneously generated by elements of
degree no greater than $j$, and for each $i$ satisfying $0\le i\le j$ we
have $ \tau (I_i)=
\#
\{$generators of  $I$  having degree less or equal $i\}$.
\item\label{ancidealsiii} $I$ is generated by forms of degree at most
$j$, and with relations of degrees at least $j+1$.
\item\label{ancidealsiv} $I$ has a generating set $f_1,\ldots ,f_\nu$ of
degrees at most $j$ and
\begin{equation}
  I_{j+1}=
\bigoplus_{1\le i\le \nu} R_{j+1-\deg f_i}f_i.
\end{equation}
\item\label{ancidealsv} $H(R/I)$ satisfies equation \eqref{ehilbineq},
and $I$ has the minimum possible number of generators for a graded ideal
defining a quotient
$R/I$ of Hilbert function $H$, namely 
\begin{equation}
\nu (I)=e_j+1 =
H_{j-1}-H_j+1=H_j-H_{j+1}+1=e_{j+1}+1.
\end{equation}
\end{enumerate}
\end{proposition}
\begin{proof} We show first that (i)-(iv) are equivalent, and then (i,ii)
$\Leftrightarrow$ (v). That (i) $\Rightarrow $ (ii) is from equation
\eqref{eqetau}. Assume (ii). Then we have for $i\le j$, 
\begin{align*}
\cod R_{-1}I_i-\cod I_i&=\tau (I_i)-1\\
&=\tau (I_{i-1})-1+\#\{\text {generators of degree $ i$ }\}\\
&=\cod (I_{i-1})-\cod (R_1\cdot I_{i-1})+\dim I_i-\dim (R_1\cdot
I_{i-1})\\  &=\cod I_{i-1}-\cod I_i,
\end{align*}
hence $\cod R_1\cdot I_i=\cod I_{i-1}$. Since always
$R_{-1}\cdot I_i\supset I_{i-1}$ the equality of dimensions shows
$R_{-1}\cdot I_i=I_{i-1}$ for $i\le j$: this and $I$ generated by degree
$j$ shows that $I$ is the ancestor ideal of $I_j$, so (ii) implies (i).
Suppose $i\le j$. We have
\begin{align*} \dim I_{i+1}&=
\dim I_i+ \nu( I_{\le i+1}) - \#\{ \text{relations
of $I$ in degrees } \le i+1.\}\\
\tau (I_i)&=\nu (I_{\le i})-\# \{\text{relations
of $I$ in degrees } \le i+1\},
\end{align*}
hence we have (ii) $\Leftrightarrow $ (iii). The condition  (iii) is
evidently equivalent to (iv). We have shown (i)-(iv) equivalent.\par
Assuming (i), (v) is a consequence of Proposition \ref{etau} equation
\eqref{ehilbineq} and Theorem \ref{basic}\eqref{basicii}. Assuming (v)
we have that $I$ has a generating set of degrees no greater than $j$,
and for
$i\le j+1$,
\begin{equation*}
 \dim R_i-\dim R_1\cdot I_{i-1}=\#\{\text { 
generators of degree $i$ }\},
\end{equation*}
implying (ii). This completes the proof. 
\end{proof}
\begin{corollary}\label{ancidealchar} The ideal $I\subset k[x,y]$ is an
ancestor ideal if and only if the highest degree $\beta_1$ of any generator and the
lowest degree
$\beta_2$ of any relation satisfy $\beta_1+2\le \beta_2$. Then $I$ is the
ancestor ideal of $I_j$ for each $j$ satisfying $\beta_1\le j\le
\beta_2-2$.
\end{corollary}
\begin{proof} The Corollary is immediate from (i) $\Leftrightarrow$ (iii)
in Proposition \ref{ancideals}.
\end{proof}
\begin{example} Let $H=(1,2,3,3,2,1)$ and let $I=(x^3,y^4)\subset k[x,y]$
Then $I$ is a complete intersection, with a single relation in degree 7.
It follows from Corollary \ref{ancidealchar} that $I$ is an ancestor ideal
both for
$I_4=\langle x^4,x^3y,y^4\rangle$ and for $I_5$.
\end{example}
We will need the following well-known result \cite{Mac1,I2}
\begin{corollary}\label{commonfactor} Let $I\subset R=k[x,y]$ be an ideal
satisfying
$H(R/I)=T, \lim_{i\to \infty}T_i=c $ where $ c=c_T>0$. Then $I=f\cdot I'$
where the common factor $f$ satisfies $\deg f=c$, and where $R/I'$
is an Artinian quotient of Hilbert function $T:c$, where
\begin{equation}\label{T:c}
(T:c)_i=T_{i+c}-c.
\end{equation}
\end{corollary}
\begin{proof}
Let $T_s =c, T_{s -1}>c$, and suppose $\mu =\mu (T)=\min\{ i\mid
T_i\ne i+1\}$ be the order of any ideal
$I$ of 
$R$ having Hilbert function $H(R/I)=T$
 (so
$I_\mu\ne 0, I_{\mu -1}=0$). Then we have
\begin{equation}\label{anyidealancinclude}
\overline{I_1}\subset
\overline{I_2}\subset
\ldots \subset
\overline{I_i}\subset \cdots \subset\overline{ I_s}=(f),f=\GCD
(I_s).
\end{equation}
Here $\overline{I_s}=(f)$ since evidently $\tau (I_s) = \cod
I_s-\cod I_{s +1}+1=1$, and we have $f\mid I$. The Corollary
follows.
\end{proof}\par
We turn now to characterizing the Hilbert functions of
level algebras and the algebras $R/(V)$.
\begin{lemma}\label{NT} The Hilbert function $N$ of a level algebra
$\LA (V)$ determined by the vector subspace $ V\subset R_j, \dim V =d $
satisfies 
\begin{align}\label{eN}
\tau (V)&\le \min\{
d,j+2-d\}, \, N_j=j+1-d, N_i=0 \text { for } i>j, \text { and }\notag\\
e_{j+1}(N)&=j+1-d\ge e_j(N)=\tau (V)-1\ge e_{j-1}(N)\ge \cdots .
\end{align}
The Hilbert function $T=H(R/(V))$ for the algebra $R/(V)$ determined by
the vector subspace
$V\subset R_j,
\dim V=d$ satisfies
\begin{align}\label{eT}
\tau (V)&\le \min\{
d,j+2-d\}, T_j=j+1-d, T_i={i+1} \text{ for }i<j,\,\text { and }\notag\\
 e_j(T)&=d-1\ge
e_{j+1}(T)=\tau (V)-1\ge e_{j+2}(T)\ge \cdots .
\end{align}
\end{lemma}
\begin{proof}  Immediate from the definitions of $\LA (V), \GA(V)$ and
Proposition
\ref{etau}, equation
\eqref{ehilbineq}.
\end{proof}
\begin{definition}\label{defNT} Let $d,j$ be positive integers satisfying
$d\le j$. Let $H$ be an acceptable $O$-sequence as in Definition
\ref{defaccept}. The
\emph{nose}
$N_H$ is the sequence
\begin{equation}\label{eNH}
N_H=(H_0,\ldots ,H_{j-1},H_j=j+1-d,0),
\end{equation}
and the \emph{tail} $T_H$  (the Hilbert function is looking to the left!)
is the sequence
\begin{equation}\label{eTH}
T_H=(1,2,\ldots
,j,H_j=j+1-d,H_{j+1},\ldots ,H_i, \ldots ).
\end{equation}
 A pair of sequences
$(N,T), N=(1,\ldots ,N_j,0),\, T=(1,2,\ldots, j,T_j,T_{j+1},\ldots )$ is
\emph{compatible} for
$(d,j)$, if
$N_{j-1}-N_j=\tau-1=T_j-T_{j+1}$, and each of $N,T$ can arise as above
from acceptable $O$ sequences $H,H'$: $N=N_H, T=T_{H'}$. For $(N,T)$
compatible, we  define
$H(N,T)$ by 
\begin{equation}\label{joinhilb}
H(N,T)=\begin{cases}N_i& \text { for } i\le j\\
T_i &\text {for }i\ge j.
\end{cases}
\end{equation}
We let $\LA_N(d,j)$ parametrize all level algebras $\LA(V), V\subset R_j,
\dim V=d$), as a subscheme of
$\Grass(d,R_j)$. We define
$\GA_T(d,j)\subset
\Grass(d,R_j)$ similarly as the parameter variety for all graded algebras
$\GA(V)=R/(V),
 V\subset R_j, \dim V=d$, having Hilbert function $H(\GA(V))=T$. As we
shall see, the maps
$V\to
\LA(V)$ and
$V\to \GA(V)$ give open dense immersions from $\LA_N(d,j)$ to $G(N)$, the
the projective variety paremetrizing graded
ideals $I$ of Hilbert functions $H(R/I)=N$, and from
$\GA_T(d,j)$ to $G(T)$ (Theorem \ref{dimstrata} \eqref{dimstrataii}).
\end{definition}\noindent
 \begin{urem} Suppose that $H$ satisfies $H=H(\Anc (V))$; then $\LA (V)$,
$\GA (V)$, respectively, have Hilbert functions $N_H,T_H$, respectively.
Also, we have
$H(N_H,T_H)=H$ in the sense of equation
\eqref{joinhilb}.
\end{urem} 
Recall that $\Grass_\tau (d,j)$ denotes the subfamily of $\Grass (d,R_j)$
parametrizing $d$-dimensional vector subspaces $V\subset R_j$ with $\tau
(V)=\tau$. We will later show that $\Grass_\tau (d,j)$ is irreducible.
We let $\rem (a,b)=b-\lfloor b/a\rfloor\cdot a$. For an
integer
$\tau$ satisfying $1\le \tau \le \min (d,j+2-d)$, we define $H_\tau
(d,j)$ as the Hilbert function corresponding to the pair of partitions
$(P_\tau (d,j),Q_\tau (d,j))$ of $(d,j+1-d)$ for which $P$ has at most one
of its parts different from $\tau$, 
$Q$ has at most one part different from $\tau-1$. Thus,
\begin{equation}\label{eHtau}
 P_\tau (d,j)=(\tau ,\ldots 
\tau ,\rem (\tau ,d)), \,\, Q_\tau (d,j)=(\tau-1,\ldots ,\tau -1, \rem
(\tau -1, j+1-d)). 
\end{equation}
Here $P_\tau (d,j)$ has $\lfloor d/\tau \rfloor$ parts of size $\tau $,
and if
$\rem (\tau,d)\ne 0$ one further part; likewise the partition
$Q_\tau (d,j)$ has $\lfloor (j+1-d)/(\tau -1) \rfloor$ parts of size 
$\tau -1$ and at most one further part.  We have, letting $a=j+1-d$,
\begin{equation}\label{eHtaub}
H_\tau (d,j)_i=
\begin{cases}&\min\{ i+1,a+(\tau -1)(j-i)\} \text { for } i\le j\\
&\max\{0,a-(\tau -1)(i-j)\}  \text { for } i>j.
\end{cases}
\end{equation}
\par We now show our main result characterizing the Hilbert function
strata of the three algebras attached to $V$. In each of equations
\eqref{edimN},\eqref{edimT},\eqref{ecodN},\eqref{ecodT}, below
the term on the far right has the same form as the terms in the sum
enclosed in parentheses; we have broken out the single term for clarity,
since for example
$e_{j+1}(N)=j+2-d-\tau
\ne e_{j+1}(H)=\tau-1$. In the equations below $e_i=E(H)_i=H_{i-1}-H_i$
throughout. We will show analogous equations for the
codimensions of the strata in terms of the graded Betti numbers in
Section~2.2, Theorem
\ref{codpartition}.  Note that the
dimension equations
\eqref{edimH},\eqref{edimN},\eqref{edimT} are written essentially in
terms of the partitions
$P,Q$ which are determined by $E(H)$. 
\begin{theorem}\label{dimstrata} Let $r=2$, let $k$ be an infinite field, 
and fix positive integers
$d\le j$. Let
$H$ be a proper acceptable
$O$-sequence in the sense of Definition \ref{defaccept}. Then 
\begin{enumerate}[A.]
\item\label{dimstrataii} Assume $k$ is algebraically closed. Each of the
schemes
$\Grass_H(d,j),
\LA_N(d,j)$, $\GA_T(d,j)$ has an open cover by opens in affine spaces of
the given dimension. Each such scheme is irreducible, rational and smooth.
Each is an open dense subscheme of the corresponding scheme $G(H),G(N),$
or $G(T)$ parametrizing all graded ideals of the given Hilbert function.
\item\label{dimstrataiii}
Let
$\lim_{i\to\infty} H_i=c_H$. The dimensions of
$\Grass_H(d,j)$, and of the related varieties satisfy
\begin{align}\label{edimH}
\dim \Grass_H(d,j)&=c_H+\sum_{i\ge \mu(H)}
(e_i+1)(e_{i+1}),\\
\dim \LA_N(d,j)&=\left( \sum_{\mu (N)\le i< j}
(e_i+1)(e_{i+1})\right) +(e_j+1)(j+1-d)\label{edimN}\\
\dim \GA_T(d,j)& =c_T+\left( \sum_{i\ge j+1}(e_i+1)(e_{i+1})\right)
+d(e_{j+1}).\label{edimT}
\end{align}\noindent
\item\label{dimstrataiv} The codimension of $\Grass_H(d,j)$ and of related
varieties in
$\Grass (d,R_j)$ satisfy
\begin{align}\label{ecodH}
\cod \Grass_H(d,j)&=\cod \LA_N(d,j)+\cod \GA_T(d,j)-\cod \Grass_\tau
(d,j),\\ 
\cod \LA_N(d,j) &=\left( \sum_{\mu (N)\le i< j}
(e_{i+1}-e_i)(i-N_{i-1})\right) +(d-\tau )(j+2-d-\tau )\label{ecodN}\\
\cod \GA_T(d,j)&= (2d-2-j)c_T+\left( \sum_{i\ge j+1}
(e_i-e_{i+1})(T_{i+1})\right) +(d-\tau )(j+2-d-\tau),\label{ecodT}\\
\cod \Grass_\tau (d,j)&=(\dim V-\tau)(\cod V-(\tau-1))=(d-\tau
)(j+2-d-\tau
).\label{ecodtau}
\end{align}
\end{enumerate}
\end{theorem}
\begin{proof} That each such $H$ occurs as $H(R/\overline{V})$ for some
such $V$ is a consequence of Proposition \ref{ancideals}
\eqref{ancidealsi} equivalent to \eqref{ancidealsv}, and Theorem
\ref{basic} \eqref{basiciii}. That each scheme has a cover by opens in
affine spaces of the given dimension, and the dimension formulas
themselves also follow from Theorem~\ref{basic}, applied to the relevant
Hilbert functions $H,N,$ or $T$, respectively. In each case
the schemes parametrize those ideals of the given Hilbert function
having the minimum possible number of generators, hence when $k$ is
algebraically closed, they are by Theorem
\ref{basic} open dense subschemes of the schemes
$G(H),G(N)$, or $G(T)$, respectively, that parametrize all graded ideals
of the Hilbert function (not just those that are $\overline{V},L(V)$, or
$(V)$, respectively with $V=I_j$). The codimension formulas are
consequences of the dimension formulas, as we will now show. We
begin by verifying
\eqref{ecodN}, whose right side we denote by $L(N)$. Since for
$I=\overline{V}\mid H(R/I)=H$ we have by Proposition
\ref{ancideals}\eqref{ancidealsii},\eqref{ancidealsiii} there are no
relations among the generators in degrees less or equal j+1, we have
\begin{equation*} 
i-N_{i-1}=\dim I_{i-1}=\tau (I_{i-1})+\tau
(I_{i-2})+\cdots =(e_{i-1}+1)+(e_{i-2}+1)+\cdots .
\end{equation*}
We have, noting that $\sum_{i<j}(e_i+1)=\dim I_{j-1}=d-\tau$,
\begin{align*}
\begin{split}
\dim \LA_N+ L(N)&=\sum_{i<j}(e_{i+1}-e_i)\left(
(e_{i-1}+1)+(e_{i-2}+1)+\cdots
\right)\\&\qquad
+\sum_{i<j}(e_i+1)e_{i+1}+(e_j+1)(j+1-d)+(d-\tau)(j+2-d-\tau)
\end{split}\\
&=\sum_{i<j}e_j(e_i+1)+(e_j+1)(j+1-d)+(d-\tau)(j+2-d-\tau)\\
 &=(\tau -1)(d-\tau)+\tau (j+1-d)+(d-\tau)(j+2-d-\tau)\\
&=d(j+1-d)=\dim \Grass (d,R_j).
\end{align*}
	It follows that $L(N)=\cod \LA (N)$, which is \eqref{ecodN}.\par
We now show \eqref{ecodT}, first when $c_T=\lim_{i\to\infty}T_i=0$.
Letting $L(T)$ denote the right side of \eqref{ecodT}, with the last term
on the right included in the sum (here $e_j(T)=j-(j+1-d)=d-1$), and noting
that since
$c_T=0, T_{i+1}=e_{i+2}+e_{i+3}+\cdots$, we have in this case
\begin{align*} \dim \GA_T(d,j)+L(T)&=\sum_{i\ge j+2}(e_j(T)+1)\cdot
e_i+d(e_{j+1})\\
&=d (T_{j+1})+d(\tau-1)=d(j+1-d)=\dim \Grass (d,R_j),
\end{align*}
thus we have $L(T) =\cod \GA_T(d,j)$ when $c_T=0$. When $c_T>0$, the
formula results from a comparison with the same sums for $T'=T:c$ (see
Corollary \ref{commonfactor}).\par We now show the formula
\eqref{ecodtau} for
$\cod
\Grass_\tau (d,j)$. Since $\Grass_\tau (d,j)=\bigcup_{\tau
(H)=\tau}\Grass_H(d,j)$, we will need to use that its largest-dimensional
stratum is
$\Grass_{H_\tau}(d,j)$, where $H_\tau =H_\tau (d,j)$ is defined above in
equation
\eqref{eHtaub}. Although this fact can be seen from equation
\eqref{edimH}, it is more readily apparent from \eqref{eHtau} and the
codimension formula
\eqref{coda}
in  terms of the partitions $(A,B)=(P^\ast,Q^\ast)$ of Theorem
\ref{writeA}; it is also, of
course, a consequence of the irreducibility of
$\Grass_\tau(d,j)$, with $\Grass_{H_\tau}(d,j)$ being a dense open
subscheme, shown below for $k$ algebraically closed in Corollary
\ref{grasstauirred}. We have by \eqref{edimH} and \eqref{eHtau},
\begin{align}
\dim \Grass_{H_\tau}(d,j)&=\sum_{i<
j}(e_i+1)(e_{i+1})+\sum_{i\ge j}(e_i+1)e_{i+1}\notag\\
&=\sum_{i< j}(e_i+1)\cdot(\tau -1)+\tau \cdot\sum_{i\ge j}e_{i+1}\notag\\
&=(d-\tau)(\tau -1)+\tau (j+1-d)=\tau (j+2-\tau)-d,\label{edimtau}
\end{align}
whence we have $\cod \Grass_{H_\tau}(d,j)=(d-\tau)(j+1-d-(\tau -1))$,
which is \eqref{ecodtau}, with, as mentioned, the dense open subscheme
$\Grass_{H_\tau}(d,j)$ in place of $\Grass_\tau (d,j)$.\par
We now show \eqref{ecodH}, which is equivalent to the analogous equation
with dimension replacing codimension. We have evidently from
\eqref{edimH},\eqref{edimN}, and \eqref{edimT}, since
$e_j(H)=e_{j+1}(H)=\tau -1$,
\begin{align*}
 \dim
\LA_N(d,j)+\dim \GA_T(d,j)&=\dim
\Grass_H(d,j)+(e_j+1 )(j+1-d)+d(e_{j+1})-(e_j+1)(e_{j+1})\\
&=\dim\Grass_H(d,j)+\tau (j+1-d-(\tau-1))+d(\tau -1)\\
&=\dim\Grass_H(d,j)+\dim\Grass_\tau (d,j),
\end{align*}
using \eqref{edimtau}. This completes the proof of Theorem
\ref{dimstrata}.
\end{proof} 
\begin{corollary}\label{NTintprop}
Let $d,j,\tau$ be positive integers with $d\le j$, and let $H$ be an
acceptable $O$-sequence in $\mathcal H(d,j)_\tau$. Let
$N=N_H, T=T_H$ be the sequences of equations
\eqref{eNH},
\eqref{eTH} or Definition~\ref{defNT}. Then
 $\LA_N(d,j)$ and
$\GA_T(d,j)$ intersect properly in $\Grass_\tau (d,j), \tau =e_j+1$,
and
$\LA_N(d,j)\cap \GA_T(d,j)=\Grass_H(d,j)$.
\end{corollary}
\begin{theorem}\label{allH} Let $d,j$ be positive integers with $d\le j$.  Let $(P,Q)$
be a pair of partitions satisfying 
\eqref{partitionsi} and \eqref{partitionsiii} of Lemma
\ref{partitions}. 
\begin{enumerate}[i.]
\item\label{allHi}
The set of proper $O$-sequences $H$ as in equation \eqref{accept} that are
\emph{acceptable} for $(d,j)$ as in Definition \ref{defaccept}, is
identical with $\mathcal H(d,j)=\mathcal H(d,j,2)$, the set that occur as
the Hilbert functions
$H(\Anc (V))$ for some
$d$-dimensional vector space $V\subset R_j$. 
\item\label{allHii} All proper $O$-sequences $H$ satisfying the conditions
of Corollary
\ref{accept} occur as the Hilbert function of an ancestor algebra of a
proper vector subspace $V\subset R_j$.
\item\label{allHiii} Fix $\tau=\tau (H)$. The pairs of partitions
$(P,Q)$ of $(d,j+1-d -c)$ where $c\le j+1-d-\tau$,
satisfying the condition of Lemma
\ref{partitions}\eqref{partitionsi} (that $P$ has at least one part $\tau$ and no
larger parts, and $Q$ has at least one part $\tau -1$ and no larger parts)
 are exactly the pairs that occur as the partitions
$P(H),Q(H)$ for those Hilbert functions $H\in \mathcal H(d,j)$ satisfying
$\tau=e_j+1$ fixed and $c_H=c$.
\end{enumerate}
\end{theorem}
\begin{proof}
The Corollary \ref{NTintprop} is immediate from Theorem \ref{dimstrata}.
Theorem \ref{allH}
\eqref{allHi} follows
from Proposition \ref{etau} and
\eqref{edimH}: the lowest value for  $\dim
\Grass_H(d,j), H$ acceptable is one, which occurs only for $d=j,
H=(1,1,\ldots )$. Theorem
\ref{allH}\eqref{allHii},
\eqref{allHiii} follow from Theorem \ref{allH}\eqref{allHi} and Lemma
\ref{partitions}.
\end{proof}\par
We now use our results to count the number of level algebra and related
Hilbert functions, given $(d,j)$. We first define the $q$-binomial
series, a power series in $q$ 
\begin{equation}\label{eqbinom}
\binom{{\bf a+b}}{{\bf
a}}=\frac{(q^{a+b}-1)(q^{a+b}-q)\cdots
(q^{a+b}-q^{b-1})}{(q^b-1)(q^b-q)\cdots (q^b -q^{b-1})}.
\end{equation} 
Recall that the number $p(a,b,n)$ of
partitions of $n$ into at most $b$ parts, each less or equal to $a$
is given by the coefficent of $q^n$ in the $q$-binomial series
${\bf \binom{a+b}{b}}$ \cite[Proposition 1.3.19]{St2}. We denote by
$p(n)$ the number of partitions of $n$, and by $p_k(n)$ the number of
partitions of $n$ into exactly $k$ parts (or, equivalently, partitions
of  $n$ with a largest part equal to $k$). Evidently, there are
$p(a-1,b-1,n-a-(b-1))$ partitions of $n$ into exactly $b$ parts, with
largest part $a$.
\begin{corollary}\label{count} Let $d,j$ be positive integers with $d\le
j$.  We assume $V\subset R_j, \dim V=d$.
\begin{enumerate}[A.] 
\item The level algebra Hilbert functions $N$ of socle degree
$j$ with $N_j=j+1-d, \tau (I_j)=\tau$ correspond one to one as in
\eqref{edefP} with the
$p_\tau (d )$ partitions
$P$ of
$d$ with largest part $\tau$. Here $\tau $ runs through all
integers less or equal $\min\{ d,j+2-d\}$.
\item The level algebra Hilbert
functions $N$ of socle degree
$j$ with $N_j=j+1-d, \tau (I_j)=\tau$ having order $\mu (N)=\mu$ correspond one to one as in
\eqref{edefP} with the
$p(\tau -1,j-\mu,d-\tau-(j-\mu))$ partitions
 of
$d$ into exactly
$j+1-\mu$ nonzero parts with largest part $\tau$. There are $p(\tau, j+1-\mu,d)$ level algebra Hilbert functions $N$
with $(\tau (N)\le 
\tau,\mu (N)\ge \mu )$, and fixed $(d,j)$.
\item The Hilbert functions
$T$ for Artinian algebras
$A=R/(V), \tau (V)=\tau $ correspond one to
one as in \eqref{edefQ} to the $p_{\tau-1}(j+1-d)$ partitions $Q$ of
$j+1-d$ having largest part
 $\tau -1$.
\item The Hilbert functions
$T$ for Artinian algebras
$A=R/(V), \tau (V)=\tau $, where $T_{s-1}\ne
0$  but
$T_s=0$ correspond one to one as in \eqref{edefQ} to
the  $p(\tau-1,s-j-1,j+1-d-\tau -(s-j-1))$ partitions of
$j+1-d$ into
$s-j$ parts, with largest part $\tau -1$. There are $p(\tau -1,s-j,j+1-d)$
such Hilbert functions $T$ with $(\tau (T)\le \tau, s(T)\le s)$ and
fixed $(d,j)$. 
\item There are $p_\tau (d)\cdot p_{\tau -1}(j+1-d-c)$ acceptable Hilbert
functions $H$ as in Defintion \ref{defaccept}, having $\tau (H)=\tau,
c_H=c$. This is the subset of $\mathcal H(d,j)$ delimited in
Theorem
\ref{allH}\eqref{allHiii}.
\end{enumerate}
\end{corollary}
\begin{proof} The Corollary follows immediately from Theorem
\ref{allH}, and Lemma
\ref{partitions}.
\end{proof}\par
\subsection{Minimal resolutions of the three algebras of
$V$, and partitions}\label{minressub}
 In this section we relate the sets of graded Betti
numbers of the ancestor algebra $\Anc (V)$, the level algebra $\LA(V)$,
and the usual graded algebra
$\GA (V)$ determined by a vector space of degree-$j$ homogeneous
elements of $R$. These depend on several partitions
$A,B$ derived from the Hilbert function $H(\Anc (V))$ -- from the
generator and relation degrees of the ancestor ideal $\overline{V}$. We
also give further codimension formulas for the Hilbert function strata, in
terms of the graded Betti numbers, or natural invariants of the
partitions. The following results were not in the original preprint
\cite{I1}. They are inspired by the special case \eqref{codf}
below, a formula for
$\cod
\GA_T(d,j)$ in
\cite{GhISa}, which arose from a
geometric tradition in studying the restricted tangent bundle from
projective space to an embedded rational curve (see also \cite{Ra,Ve}). We
will suppose that $V\subset R_j$ satisfies
$H(R/\overline{V})=H$; unless otherwise stated we will suppose also  that
$\lim_{i\to\infty}H_i=0$,  Then, as we shall see in Lemma~\ref{writeA},
the ancestor algebra
$\Anc(V)= R/\overline{V}$, the algebra $\GA (V)=R/(V)$ and the level
algebra $\LA (V)$ determined by $V$ have graded Betti numbers given by
certain sequences/partitions $A,B$ as follows, 
\begin{align}\label{eminresa}
0\to &\sum_{i=1}^{\tau-1} R(-j-1-b_i)\to\sum_{i=1}^\tau R(-j-1+a_i)\to R
\to R/\overline{V}\to 0,\\
0\to &R(-j-2)^{j+1-d}\to\sum_{i=1}^\tau R(-j-1+a_i)\oplus
R(-j-1)^{j+2-d-\tau}\to R
\to
\LA (V)\to 0, \label{eminresc}
\\
0\to &\sum_{i=1}^{\tau -1} R(-j-1-b_i)\oplus R(-j-1)^{d-\tau}\to
R(-j)^d\to R
\to R/(V)\to 0, \label{eminresb},
\end{align}
where we assume that the sequences $A=(a_1,\ldots ,
a_\tau)$ and $B=(b_1,\ldots ,
b_{\tau -1})$ defined by \eqref{eminresa} are listed in decreasing order
$a_1\ge
\cdots
\ge a_\tau$ and $b_1\ge \cdots \ge b_{\tau -1}$.
\begin{definition}\label{defA-D}
When $\lim_{i\to\infty}H_i=0$, we define partitions $A,B$ given $V$ 
by
\eqref{eminresa}; we will show that they depend only on $H$, and
evidently they are the same that occur in
\eqref{eminresc} and
\eqref{eminresb} (See Lemma \ref{writeA}). By
$A+\underline{1}$ we mean the partition whose parts are
$A+\underline{1}=(a_1+1,a_2+1,\ldots )$. We denote by $C$ the partition of
$j+2$ having
$j+2-d$ parts given by
$(A+\underline{1})\cup [1^{j+2-d-\tau}]$, namely
$A+\underline{1}$ with
$j+2-d-\tau$ parts of size one adjoined; and
we denote by $D$ the partition of $j$ having $d-1$ parts given by
$(B+\underline{1})\cup [1^{d-\tau}]$, namely $B+\underline{1}$ with
$d-\tau$ ones adjoined. \par
When $\lim_{i \to \infty}H_i=c_H\ge 0$ we define $A,B$ from the minimal
resolution of $\overline{V:f}$, where $f=\GCD (V)$; then $A,B$ depend
only on $H:c_H$ (see \eqref{T:c}). We define $C,D$ in this case as above
from $A,B$; here $C$ again partitions $j+2$, but $D$ partitions
$j+2-d-\tau-c$.
\end{definition}
 Evidently, the
generator degrees of the ideal $L(V)$ defining $\LA(V)$ in
\eqref{eminresc} are $\underline{j+2}-C$
and the relation degrees of $(V)$ in \eqref{eminresb} are
$\underline{j}+D$. We have chosen $A$ and $B$, then $C$ and $D$ in a
symmetric fashion so that they partition integers depending only on $d$
and $j$; this allows application of Lemma \ref{dualpartition} later.
 As we
shall see, the partitions
$A,C$ depend only on
$N=N_H$, determined by
$H_{\le j}$; and
$B,D$ depend only on $T=T_H$, determined by $H_{\ge j}$ (see Definition
\ref{defNT}). To describe this dependence simply, we use the dual
partition.
\begin{definition}\label{dualpart}
Let $A=(a_1,\ldots ,a_k),\  a_1\ge a_2\ge \ldots $ be a partition of
$a=\sum a_i$ into
$k$ non-negative parts (some may be zero). Recall that the
{\emph{Ferrers graph}}
${\mathfrak F} (A)$ of
$A$ has
$k$ rows, the $i$-th row of length $a_i$. We denote by
$A^\ast = (a_1^\ast ,a_2^\ast ,\ldots )$ the {\emph{dual partition}} of $a$,
whose Ferrers graph is obtained by switching rows and columns in the
Ferrers graph $\mathfrak F (A)$. Here also, $a_i^\ast$ is the number of
parts of
$A$ of length greater or equal $i$.\par 
\end{definition}
\begin{lemma}\label{writeA} Let $d,j$ be positive integers satisfying
$d\le j$, and let $H$ be an acceptable $O$-sequence as in Definition
\ref{defaccept}, and suppose that  $c_H=\lim_{i\to\infty} H_i=0$. Then the
algebras
$\Anc (V), \LA (V)$, and $ R/(V)$ have minimal resolutions whose graded
Betti numbers are given by 
\eqref{eminresa},
\eqref{eminresc}, \eqref{eminresb}. We have
\begin{equation}\label{sumgendegs}
\sum_{1=1}^\tau a_i=d ;
\end{equation}
 $A$ satisfies $a_i\ge 1$, and $A$ has dual partition
$A^\ast =P=(\tau, \tau (R_{-1}\cdot V),
\tau (R_{-2}V),\ldots )$ of $d$, and 
\begin{equation}\label{AdualE}
a^\ast_i=\tau (R_{-i+1}\cdot
V)=e_{j+1-i}(H)+1.
\end{equation}
 Also
\begin{equation}\label{sumreldegs}
\sum_{i=1}^{\tau-1}b_i=j+1-d;
\end{equation}
$B$ satisfies $b_i\ge 1$, and $B$ has dual
$B^\ast=Q=(e_{j+1}(H),\ldots )$ of
$j+1-d$, and $b^\ast_i=e_{j+i}$. We have for $i\ge 0$
\begin{align}\label{edimgens}  \dim I_{j-i}&=\sum_u \mid a_u-i\mid^+,\\
 H_{j+i}&=\sum_u
\mid b_u-i\mid^+.\label{ecodrels} 
\end{align}
Likewise, the partition $C$ has dual the partition
$(E(N)_{j+1}+1,E(N)_j+1,\dots )$ of $j+2$
\begin{equation}\label{eCdual}
C^\ast=(j+2-d,\tau (V),\tau (R_{-1}V),
\tau (R_{-2}V),\ldots ),
\end{equation}
 and $D$ has dual the partition $E(T)_{\ge j}$ of $j$
\begin{equation}\label{eDdual}
D^\ast=(d-1,e_{j+1},e_{j+2}\ldots
).
\end{equation}
When $\lim_{i\to \infty}H_i=c_H>0$, then $A$ from Definition
\ref{defA-D} satisfies all the statements above, including
\eqref{sumgendegs},\eqref{AdualE},\eqref{edimgens}; and $B$ is a
partition of
$j+1-d-c_H$ into $\tau -1$ parts. Also, $B^\ast$ satisfies the same
condition above, and $H_{j+i}=c_H+\sum_u\mid b_u-i\mid^+$ in place of
\eqref{ecodrels}. Also, $C^\ast$ satisfies \eqref{eCdual}, and $D^\ast$
satisfies
\eqref{eDdual}.
\end{lemma}
\begin{proof} We first assume $\lim_{i\to\infty}H_i=0$, The definition of
$\overline{V}$ shows that it is generated in degrees less or equal $j$,
and Proposition
\ref{ancideals} shows that
$\overline{V}$ has no relations in degrees less or equal $j+1$. Thus,
equation \eqref{eminresa} defines ordinary partitions A and B, with
nonzero parts. Given the definition of $A,B$ in
\eqref{eminresa}, the graded Betti numbers shown in \eqref{eminresc},
\eqref{eminresb} for the level algebra $\LA(V)$ and the algebra
$\GA(V)=R/(V)$ follow immediately from the definitions of these algebras
from $\overline{V}$ in Definition \ref{basics}, and the relations among
them given in Remark
\ref{ancid}. For example, since the ideal $L(V)$ defining the level
algebra $\LA (V)$ satisfies $L(V)=\overline{V}+M^{j+1}$ one obtains $L(V)$
it by adding
$H_{j+1}=(j+1-d-(\tau -1)=j+2-d-\tau$ generators of degree $j+1$, and
evidently all the relations are in degree $j+2$, since the socle of
$R/L(V)$ lies solely in degree $j$; this shows
\eqref{eminresc}.\par Proposition
\ref{etau} shows that for
$i\ge 0,
\,
\tau (R_{-i}\cdot V)=e_{j-i}(H)+1$, so $\tau (R_{-i}\cdot V)$ depends
only on initial portion
$N_H$ of $H$. We have from Proposition
\ref{ancideals} \eqref{ancidealsiii}, and the definition of $A^\ast$ that
for
$i\ge 1$, 
\begin{equation*}
\tau(R_{-i+1}V)= \# \{u\mid a_u\ge i\} =a^\ast_{i}.
\end{equation*}
It follows from \eqref{sumtaunose} that $\sum a_i=\sum_{i=1}
a^\ast_i=d$, which is \eqref{sumgendegs}. \par
Concerning $B$, we have from \eqref{eminresa}, that for $i\ge 0$
\begin{align*} 
H_{j+i}&=H_j-(\tau -1)i+\sum_{u\mid b_u\le i-1}(i+1-b_u); \text { thus }\\
e_{j+i}&=\tau -1-\sum_{u\mid b_u\le i-1}(-1)=\tau -1-(\#\{\text {
relations }\}  -b^\ast_{i})\\
&=b^\ast_{i}.
\end{align*}
Thus we have 
\begin{equation*}
\sum b_i=\sum b^\ast_i=\sum_{u\ge 1}e_{j+u}=H_{j}=j+1-d,
\end{equation*}
which is \eqref{sumreldegs}. It remains to show \eqref{edimgens} and
\eqref{ecodrels}. We have for $i\ge 0$,
\begin{align}
H_{j+i}&=H_{j}-(e_{j+1}+\cdots + e_{j+i})\notag\\
&=j+1-d-(b^\ast_1+\cdots +
b^\ast_{i})=b^\ast_{i+1}+b^\ast_{i+2}+\cdots\label{BandH}\\
 &=\sum\mid b_u-i\mid^+,\notag
\end{align}
which is \eqref{ecodrels}. Since $\overline{V}$ has no relations in
degrees less or equal $j+1$, we have for $i\ge 0$,
\begin{align*}
\dim I_{j-i}=\sum_{a_u\ge i+1}(a_u-i)=\sum_{u=1}^\tau\mid a_u-i\mid^+,
\end{align*}
which is \eqref{edimgens}. This completes the proof in the case 
$\lim_{i\to\infty}H_i=0$. \par
When $\lim_{i\to\infty}H_i=c_H>0$, the assertions
at the end of the Lemma follow from Definition
\ref{defA-D} of
$A,B$ in this case that uses $V:\GCD (V)$, Corollary
\ref{commonfactor} and the Lemma for $V: \GCD (V)$.
\end{proof}\par
We denote by $\mid n\mid^+$ the integer $n$ if $n\ge 0$, or $0$
otherwise. We will denote by $\underline{n}$ the sequence $(n,n,\ldots )$
of appropriate length. For a partition $A=(a_1,\ldots ), a_1\ge a_2\ge
\cdots $ we denote by
$\ell (A)$ the sum 
\begin{equation}\label{ellA}
\ell (A)=\sum_{u\le v}\mid
a_u-a_v-1\mid^+.
\end{equation}
Recall from \eqref{ecodtau} that $\cod \Grass_\tau (d,j)$ in
$\Grass (d,R_j)$ satisfies
\begin{equation*}
\cod (\Grass_\tau (d,j))=(d-\tau)(j+2-d-\tau)=(\dim V-\tau)(\cod
V-(\tau-1)),
\end{equation*}
for any $V$ satisfying $\tau (V)=\tau$. This is a term in equation
\eqref{codd}.
\begin{theorem}\label{codpartition} Let $d,j$ be positive integers
with $d\le j$. Let $H$ be an acceptable $O$-sequence, and let
$\lim_{i\to\infty}H_i=c_H$, and let $N=N_H, T=T_H$ be the sequences
of Definition \ref{defNT}, where $c_T=c_H$. The codimensions of the
families
$\LA_N(d,j)$, $\GA_T(d,j)$, and $\Grass_H(d,j)$ in $\Grass_\tau
(d,j)$ satisfy
\begin{align}\label{codb}
\cod_\tau\LA_N&=\ell (A),\\
\cod_\tau\GA_T&=\ell (B) +(d-1)c_T,\label{codc}\\
\cod_\tau \Grass_H(d,j)&=\ell (A) + \ell (B)+(d-1)c_T.\label{coda}
\end{align}
The codimensions of these families in $\Grass (d,R_j)$ satisfy
\begin{align}\label{code}
\cod\LA_N&=\ell (C),\\\label{codf}
\cod\GA_T&=\ell(D) +(d-1)c_T,\\
\cod \Grass_H(d,j)&=\ell(C) + \ell(D)+(d-1)c_H -(d-\tau )(j+2-d-\tau
),\label{codd}\\
&=\ell(C)+\ell(B)+(d-1)c_H.\label{code}
\end{align}
\end{theorem}
\begin{proof} We first note that \eqref{codb} $\Leftrightarrow$
\eqref{code}, and \eqref{codc} $\Leftrightarrow$ \eqref{codf}; evidently
\eqref{coda} is a consequence of \eqref{codb} and \eqref{codc}, and
similarly for \eqref{codd}. Assume first that $c_H=0$. We have 
\begin{align*}
\ell (C)-\ell (A)&=\left( \sum
(a_i)\right)(j+2-d-\tau)\\ 
&=(d-\tau)(j+2-d-\tau )=\cod \Grass_\tau (d,j). 
\end{align*}
Likewise, 
\begin{align*}
\ell (D)-\ell (B)&=(d-\tau)\left( \sum_{i=1}^{\tau-1}
(b_i-2)\right)\\
&=(d-\tau )(\left( j-(d-\tau)-2(\tau-1)\right)\\
&  =\cod
\Grass_\tau (d,j).
\end{align*}
We now show \eqref{codc} when $c_H=0$. Since $\lim_{i\to \infty}T_i=0$, by
Theorem
\ref{dimstrata} equation
\eqref{ecodT} we have
\begin{equation*}
\cod \GA_T= \sum_{i\ge j+1}
(e_i-e_{i+1})(T_{i+1})+(d-1-e_{j+1})(T_{j+1}),
\end{equation*}
whence, subtracting $\cod \Grass_\tau (d,j)=(d-\tau)T_{j+1}$ and
noting
that we specify $E(H)$ below, as $e_j(H)$ is different from
$e_j(T)$, we find,
\begin{align*}
\cod_\tau \GA_T&= \sum_{i\ge j+1}
(e_i-e_{i+1})(T_{i+1})+(d-1-e_{j+1})(T_{j+1})-(d-\tau)T_{j+1}\\
&= \sum_{i\ge j}
(e_i(H)-e_{i+1}(H))(H_{i+1})=\sum_{u\ge 0}(e_{j+u}-e_{j+u+1})H_{j+u+1}\\
&=\sum_{u\ge 0}(b^\ast_{u}-b^\ast_{u+1})H_{j+u+1}\text { by Lemma
\ref{writeA} } ,\\ &=\sum^{\tau-1}_{u=1} H_{j+b_u+1}\\
&=\ell (B) \text { by \eqref{ecodrels} }.
\end{align*}
We now show \eqref{codb}.  By Theorem
\ref{dimstrata} equation \eqref{ecodN}, taking into account that the last
term on the right is $\cod \Grass_\tau (d,j)$, and by \eqref{AdualE} we
have
\begin{align*}
\cod_\tau \LA (N)&=\sum_{\mu(N)\le u<j}(e_{u+1}-e_u)(\dim
I_{u-1})=\sum_{1\le i}(e_{j-(i-1)}-e_{j-i})(\dim I_{j-(i+1)})\\
&=\sum_{1\le i}(a^\ast_{i}-a^\ast_{i+1})(\sum_u\mid
a_u-(i+1)\mid^+)\text { by Lemma \ref{writeA} and \eqref{edimgens} }
\\
 &=\sum (\# \{ a_v=i\} )(\sum_u \mid a_u-i-1\mid^+\\
&=\ell (A).
\end{align*}
The adjustment of adding $(d-1)c_H$ for the case
$\lim_{i\to\infty}H_i=c_H$ comes from a comparison with the Hilbert
function
$T': T'_i=T_{i+c}-c, c=c_H$. The partitions $B,D$ are the same for
$T$ and for $T'$, and $\dim\GA (T)=c+\dim \GA (T')$, so the codimension
of $\GA(T)$ in $\Grass(d,R_j)$ satisfies
\begin{align*} 
\cod \GA
(T)&=\cod\GA(T')+\dim\Grass(d,R_j)-\dim\Grass(d,R_{j-c})-c\\
&=\ell(D)+(d-1)c_H.
\end{align*} 
This
completes the proof.
\end{proof}
\begin{example} We take $(d,j)=(9,14)$ and $\tau =4$, then $\dim
\Grass(9,R_{14})=\dim \Grass(9,15)=9\cdot 6=54$, and $\cod \Grass_4
(9,14)=(9-4)(6-(4-1))=15$, so $\dim \Grass_4(9,14)=39.$
 Consider 
\begin{equation*}
H=(1,\dots ,12,11,9,6,3,0)\text{ with }H_{14}=6.
\end{equation*}
 Here the
sequence 
\begin{equation*}
A^\ast=(\tau,\tau
(R_{-1}\cdot V),\tau (R_{-2}\cdot V),\ldots
)=(\tau,e_{13}+1,e_{12}+1)=(4,3,2),
\end{equation*}
whose dual partition is $A=(3,3,2,1)$, with $\ell (A)=2$ while
$B^\ast =(2,2,2)$, $B=(3,3)$, for which
$\ell (B)=0$. By \eqref{eminresa} the generator degrees of
$\overline{V}$ are $(j+1-a_1,j+1-a_2,\ldots )=(\underline{j+1}-A)$. Here
the generator degrees are
$(\underline{15}-A)=(15-3,15-3,15-2,15-1)=(12,12,13,14)$. The codimension
of $\Grass_H(9,14))$ in $\Grass_4(9,14)$ is by equation
\eqref{coda} $\ell (A)+\ell (B)=2+0=2$, so
$\dim\Grass_H(9,14)=39-2=37$. The formula  \eqref{edimH} that
$\dim\Grass_H(9,14)=\sum (e_i+1)(e_{i+1})$ when applied to $E(H)_{\ge
13}=(1,2,3,3,3)$ also gives 37. Here the partition $C=(4,4,3,2,1,1,1)$
and $\ell (C)=17$, and $\cod (\Grass_H(9,14))=\ell(C)+\ell(B)=17$ in
$\Grass(9,R_{14})$ by \eqref{code}.\par Consider now
$H'=(1,\ldots ,12,11,9,6,3,2,1)$. Here
$A'=A$, but
$B'=(4,1,1)$, the dual partition to $(e_{15},\ldots )=(3,1,1,1)$, $\ell
(B')=4$, and we have
$\cod_4 \Grass_{H'}(9,14)=\ell (A')+\ell (B')=6$ in $\Grass_4(9,14)$,
giving
$\dim \Grass_{H'}(9,14)=33$.
\end{example}
\subsection{Closure of the Hilbert function strata}\label{closuresub}
We now determine the Zariski closure of $\Grass_H
(d,j)$ when $r=2$, and we show that the family $G(H)$ of graded algebra
quotients of $A$ having Hilbert function $H$ is a natural
desingularizatiion of $\overline {\Grass_H(d,j)}$
(Theorem
\ref{closureofstrata}). This is one of our main results, and certainly
the deepest.\par 
We show that the closure of a stratum $\Grass_H(d,j)$ is the union of the
more special strata $\Grass_{H'}(d,j)$, for $H'\le _\mathcal P H$, where
$\mathcal P$ is the partial order on acceptable sequences given in
Definition~\ref{defpo}. Evidently the partial order $\mathcal P$
determines related partial orders on
the sequences $N$ possible for level algebras, and to the sequences $T$
possible for graded ideals
$(V)$. For the case $r=2$ we intepret these latter
partial orders as majorization partial orders on sets of partitions (Lemma
\ref{comparepart}).
This result was suggested by an application to the
restricted tangent bundle in
\cite{GhISa}. We show that the partially ordered set $\mathcal H(d,j)$ of
acceptable Hilbert functions under the partial order $\mathcal P$ --- the
same order as that determined by Zariski closure of the varieties
$\Grass_H(d,j)$ --- is equivalent to a partially ordered set
$\mathcal PA (d,j)$ of certain pairs of partitions, under the product of
majorization partial orders (Theorem
\ref{equivpos}). \par The proof of our main
result depends on a key construction. Suppose that we are given two
acceptable Hilbert functions
$H,H'\in \mathcal H(d,j)$, with $H'\ge H$ (more special) in the partial
order $\mathcal P(d,j)$, and let $V'$ be a point of $\Grass_{H'}(d,j)$.
We build a graded ideal
$I$ of Hilbert function H, that is related as in
\eqref{eancineq} to the ancestor ideal $I'=\overline{V'}$ (Lemma
\ref{build}). This ideal $I$ determines a point of $G(H)$ lying
over the given point $V'$ of $\Grass_{H'}(d,j)$ (Theorem
\ref{closureofstrata} \ref{closureB}).
\begin{definition}\label{Porder} The \emph{length} $\mid D\mid$ of
a partition
$D$ is the sum of its parts. We recall the \emph{majorization} partial
order on partitions (see \cite{GrK}). Let $D,D'$ be two partitions
$D=(d_1,d_2,\ldots ,d_s) \mid d_1\ge d_2\ge \cdots $ and
$D'=(d'_1,d'_2,\cdots ,d'_{s'}) \mid d'_1\ge d'_2
\cdots$. We say
$D'\ge D$ if $|D'|\ge |D|$ and
\begin{equation}
\sum_{u\le i}d'_u\ge\sum_{u\le
i}d_u \text{ for all } i\mid 1\le i\le \min\{ s,s'\} .
\end{equation}
Let $D$ have $r_i$ parts of size $v_i, v_1>v_2> \cdots
>v_k$. We define for each
$s, 1\le s\le k$ the partition $D_s$ with $r_i$ parts of size $v_i,1\le
i\le s$,  and no other parts. The 
\emph{polygon} of 
$D$ is the convex graph with vertices $(0,0)$ and
\begin{equation}
 (\sum_{i=1}^s r_i, \sum_{i=1}^s r_iv_i), 1\le s\le k,
\end{equation}
the height of the $s$-th vertex being the length $| D_s|$ of $D_s$.
We define the Harder-Narasimham partial order \cite{HN} on partitions
having the same number of parts, by
$D'\ge_{HN} D$ if and only if the polygon of $D'$ is never below the polygon of $D$.
\end{definition}
The Harder-Narasimham order as stated above is a special case for bundles
of the form
$\oplus \mathcal O_{{\mathbb P}^1}(v_i)^{r_i}=\oplus \mathcal O_{{\mathbb P}^1}(d_i)$ over
${\mathbb P}^1$ of an order defined more generally by Harder-Narasimham
(see \cite{HN}). This is relevant since the 
partition $C$ corresponds to the generator degrees of
the ideal $L(V)$ defining the level algebra $\LA(V)$, and $D$ corresponds
to the relation degrees of the ideal $(V)$ determining $\GA(V)$. The
latter corresponds to the decomposition into a direct sum of line bundles
of the ``restricted tangent bundle'' to the rational curve $X$ in
$\mathbb P^{r-1}$ determined by
$V$, studied in
\cite{GhISa,Ra,Ve}; the former corresponds to the decomposition of another
natural bundle over $X$, of rank $j+2-d$. It is a general result that
specialization in a family $\mathcal V(t), t\not= t_0$ of vector bundles
having fixed Harder-Narasimham polygon over $X$ yields a bundle
$V(t_0)$ of equal or higher Harder-Narasimham polygon
\cite{BPV}. Both L. Ramella and F. Ghione et al show a converse for
the restricted tangent bundle, related to Theorem
\ref{closureofstrata} \ref{closureA} for the closure of $\GA_T(d,j)$.
\par We need a preparatory result, before giving some equivalent versions
of the partial order $\mathcal P (d,j)$.
\begin{lemma}\label{dualpartition} If $D,D'$ are two partitions of the
same integer n, then 
\begin{equation}\label{dualorder}
D'\ge
D\Leftrightarrow {D'}^\ast\le D^\ast.
\end{equation}
\end{lemma}
\begin{proof}
 It suffices to consider adjacent
partitions $D'> D$ in the partial order: then $D'$ is obtained from $D$
by increasing a part of $D$ by one and decreasing the next
smaller-or-equal block by one. A basic case is $D=(d_1,\ldots
,d_{s+1})=(a,1,\ldots ,1)$ and $D'=(d'_1,]dots ,d'_s)=(a+1,1,\ldots ,1)$.
Then $D\ast=(s+1,1,\ldots ,1)$ with $a-1$ ones, and
${D'}^\ast=(s,1,\ldots ,1)$ with $a$ ones, whence we have ${D'}^\ast <
D$. The general case has $s+1$ relevant parts for $D$, $(d_i,\ldots
,d_{i+s}=(k+a,k+1,\ldots ,k+1)$ with $d_{i-1}>d_i$, and $s+1$ relevant
parts for
$D'$,
$d'_{i+1},\ldots ,d'_{i+s}=(k+a+1,k+1,\ldots ,k+1,k)$; then $D^\ast$ has
relevant parts $(d^\ast_{k+1},\ldots ,d^\ast_{k+a+1})=(i+s,i+1,\ldots
,i+1,i)$ and ${D'}^\ast$ has corresponding parts $(i+s-1,i+1,\ldots,
i+1,i+1)$, whence ${D'}^\ast < D^\ast$. 
\end{proof}\par
  We
say a Hilbert function sequence
$T'\ge T$ if for each i,
$T'_i\ge T_i$. Recall from Definition \ref{defpo} the
partial order
$\mathcal P=\mathcal P(d,j)$ on
$\mathcal H(d,j)$:
\begin{equation}\label{eHH'}
H'\ge_{\mathcal P}H\Leftrightarrow H'_i\le H_i \text { for } i\le j \text
{ and } H'_i\ge H_i \text { for } i\ge j.
\end{equation}
 Recall from Definiton \ref{defNT} that
$(N_H)_i=H_i$ for
$i\le j$ and $0$ otherwise, and $(T_H)_i=H_i$ for $i\ge j$ and
$(T_H)_i=i+1$ for $i<j$.
In terms of the pair $N_H,T_H$  we thus have
\begin{equation*}
H'\ge_{\mathcal P}H\Leftrightarrow N'\le N \text { and } T'\ge T,
\end{equation*}
where $N'\le N$ and $T'\ge T$ in the termwise partial order on sequences.
\par
 We now determine the analogues of the partial order
$\mathcal P(d,j)$, for the pairs of partitions $(P,Q)$ from Definition
\ref{defP}, and the pairs $(A,B)$ or
$(C,D)$ from Definition
\ref{defA-D}. In the
Lemma below
$H',N', A', B',\ldots $ are more special than $H,N,A,B,\ldots $, as we
shall show in  Theorem \ref{closureofstrata}. The implications $T'\ge
T\Leftrightarrow D'\ge D\Leftrightarrow D(T')\ge D(T)$ from Lemma
\ref{comparepart}\ref{comparepartB} are shown for
 $c(T)=c(T')=0$ in \cite{GhISa}. Recall
that we showed 
$P=A^\ast$ and $ Q=B^\ast$ in Lemma \ref{writeA}.
\begin{lemma}\label{comparepart} We fix positive integers $d,j$ with $d\le j$. We
treat separately the Hilbert functions for the level algebra $\LA (V)$, graded
algebra $\GA (V)=R/(V)$ and ancestor algebra $\Anc(V)$. 
\begin{enumerate}[A.]
\item\label{comparepartA}
The following are equivalent:
\begin{enumerate}[i.]
\item\label{comparepartNi} $N'\le N$ (note: $N'$ is more special!),
\item\label{comparepartNii} $ A(N')\ge A(N)$, or
equivalently $C(N')\ge C(N)$,
\item\label{comparepartNiii}
$P(N')\le P(N)$ (i.e. ${A'}^\ast\le A^\ast$) or, equivalently
$C(N')\ge_{HN}C(N)$;
\end{enumerate}
\item\label{comparepartB}
The following are equivalent:
\begin{enumerate}[i.]
\item\label{comparepartTi} $T'\ge T$ (note: $T'$ is more special!),
\item\label{comparepartTii} (Only when $c(T)=c(T')$) $ B(T')\ge B(T)$, or,
equivalently
$D(T')\ge D(T)$,
\item\label{comparepartTiii}$Q(T')\le Q(T)$,(i.e. ${B'}^\ast\le B^\ast$)
or, equivalently
$D(T')\ge_{HN} D(T)$;
\end{enumerate}
\item\label{comparepartC} The following are equivalent:
\begin{enumerate}[i.]
\item\label{compareparti}$H'\ge_{\mathcal P} H$. meaning both $N'_H\le
N_H$ and $T'_H\ge T_H$,
\item\label{comparepartiv}  $P(H')\le P(H)$ and $Q(H')\le Q(H)$,
(i.e. both ${A'}^\ast\le A^\ast$ and ${B'}^\ast\le B^\ast$).
\item\label{comparepartii} (Only when $c_H=c_{H'}$) $A(H')\ge A(H)$ and $B(H')\ge
B(H)$
\item\label{comparepartv} (Only when $c_H=c_{H'}$)
$C(H')\ge_{HN} C(H)$ and
$D(H')\ge_{HN} D(H)$.
\end{enumerate}
\end{enumerate}
\end{lemma}
\begin{proof}
We first show \eqref{comparepartNi} $\Leftrightarrow$
\eqref{comparepartNii}   $\Leftrightarrow$
\eqref{comparepartNiii} and \eqref{comparepartTi} $\Leftrightarrow$
\eqref{comparepartTii}  $\Leftrightarrow$
\eqref{comparepartTiii}. From equation \eqref{AdualE} that
$a^\ast_i=e_{j+1-i}(H)+1$ we have for
$i\ge 1$
\begin{equation}\label{HAdual}
H_{j-i}=j+1-d+(a^\ast_1-1)+\cdots + (a^\ast_i -1)=j+1-d-i+\sum_{u=1}^i
a^\ast_u,
\end{equation}
whence we have $N_H$ satisfies, using
\eqref{dualorder}
\begin{equation}\label{AequivN}
N_{H'}\le N_H\Leftrightarrow A^\ast (N')\le A^\ast (N)\Leftrightarrow
A(N')\ge A(N).
\end{equation} 
Since $A'\ge A\Rightarrow \tau '={a'}^\ast_1\le a^\ast_1=\tau$, we have
$C'=\underline{1}+A'\cup 1^{(j+2-d-\tau ) '}\ge C=\underline{1}+A\cup
1^{(j+2-d-\tau )}$.  From Lemma \ref{writeA} we have that
$b^\ast_i=e_{j+i}$, and as in
\eqref{BandH}
\begin{equation*}
H_{j+i}=j+1-d-\sum_{u=1}^i b^\ast_{i-1},
\end{equation*}
whence we have using \eqref{dualorder}
\begin{equation}\label{BandT}
T'\ge T\Leftrightarrow {B'}^\ast\le B^\ast
\Leftrightarrow B'\ge B. 
\end{equation}
This completes the proof of the Lemma except for the  equivalences
involving $\ge_{HN}$, which we now show. Note that for the partitons $C$
or
$D$ both the number of parts and sum are fixed by the triple $(d,j,\tau
)$. That
\eqref{comparepartii}
$\Rightarrow$ \eqref{comparepartv} follows, since, considering $D$, the
vertices of the  polygon of $D$ are a subset of the vertices of the graph
of the sum function $\sum D$ of $D$, used in comparing $D$ and $D'$: thus
$D'\ge D\Rightarrow D'\ge_{HN} D$. The converse follows from the
extremality of the vertices of the graph of $\sum D$ chosen as vertices
of the Harder-Narasimham polygon. 
\end{proof}
\begin{example}\label{nsimple}{\sc $\mathcal P(d,j)$ is not a
simple order on $\mathcal H(d,j)$}
\begin{enumerate}[A.]
\item Let
$d=3,j=5$, so $H_5=j+1-d=3$. Let $H=(1,2,3,4,4,3,2,1,0)$, where $\tau =
1$,      and $\mu (H)=4$, $H'=(1,2,3,4,5,3,1,1,1,\ldots )$
 where $\tau=2$ and $\mu (H')=5$. Then $H$ and $H'$ are incomparable
in the order $\mathcal P(3,5)$ since $H_6>H'_6$ but $H_8<H'_8$.
Neither stratum is in the Zariski closure of the other. The two strata are
\emph{geometrically incomparable} in the sense that no
element of either stratum can be in the closure of a subfamily of the
other stratum, by Corollary
\ref{closure}. This example essentially involves just the \emph{tail} of
$H$, namely $T(V)=H(R/(V))$, with
$(V)$ the ideal generated by $V$
(see Definition
\ref{defNT}). 
\item We give an example of similar behavior for the level algebra strata
$\LA_N(d,j)$ --- the family of level algebras of socle degree j and type
$d$ having Hilbert function $N$. Here
$N$ is the \emph{nose} of $H$ as in Definition \ref{defNT}. To create the
example, we begin with two partitions $P: 10=4+2+2+2$ and
$P': 10=3+3+3+1$, that are incomparable in the \emph{majorization}
partial order of Definition \ref{Porder}. Thus, 
their associated sum sequences
$\sum P=(4,6,8,10),\,
\sum P'=(3,6,9,10)$ are incomparable in the termwise order on sequences.
 By
Definition \ref{defP} the corresponding sequences $E=\Delta N, E'=\Delta
(N')$ are $(3,1,1,1)$ and $2,2,2,0$, respectively, and by
Lemma \ref{partitions}\eqref{partitionsi} the dimension $d$ satisfies
$d=|P|=10$. By \eqref{ehilbineqii} the simplest such case satisfies
$j+1-d=p_1-1=4-1=3$, where 
$p_1$ is the largest part of $P$, so we have $(d,j)=(10,12),
\,
\mu (N)=\mu(N')=9, N=(1,2,\ldots ,8, 9,8,7,6,3,0)$ and $ N'=(1,2,\ldots
,8,9,9,7,5,3,0)$. Thus, $N$ and $N'$ are incomparable in the
partial order $\mathcal P_N(10,12)$ on the set of nose sequences $\{
N_H\mid H\in
\mathcal H(10,12)\}$ induced from the partial order
$\mathcal P(10,12)$ on acceptable $O$-sequences $H$. Again Corollary
\ref{closure} implies that
$\LA_N(10,12)$ and
$\LA_{N'}(10,12)$ are
\emph{geometrically incomparable} in the sense that no subfamily of either
stratum can have as limit a space
$V$ in the other stratum.  This example illustrates (Lemma
\ref{comparepart} (A)). 
\end{enumerate}
\end{example}
The following lemma is the crux of the proof that the morphism $\pi :
G(H)\to
\overline{\Grass_H(d,j)}$ is surjective (Theorem \ref{closureofstrata}).
The proof we give is basically that of the original preprint,
but we have supplied further details and made an improvement. Note that
although the Hilbert functions $H,H'$ that occur are acceptable, the
ideals $I,I'$ are \emph{not} assumed to be ancestor ideals, Thus in the
proof we are rather careful about how we use previous results. In
particular, a key step, the last in the section concerning $N$ is to show
in equation
\eqref{N1key} that
$\cod R_1\cdot I(1)_{u-1}$ satisfies a certain inequality (a similar step
for
$T$ occurs in \eqref{T1key}); the apparent clumsiness --- or perhaps we
should say, subtlety --- of the argument here is in part due to $I'$ not
being an ancestor ideal!
\begin{lemma}\label{build} Let $d,j$ be positive integers satisfying
$d\le j$, Assume that $H$ and $H'$ are
acceptable
$O$-sequences for the pair $(d,j)$ (Definition \ref{defaccept}) satisfying
$H'\ge _{\mathcal P(d,j)}H$. When $c_H=c_{H'}$ let $k$ be an arbitrary
field; otherwise assume $k$ is algebraically closed.  Let
$I'$ be a graded ideal of Hilbert function $H(R/I')=H'$. Then there is a
graded ideal $I$ of Hilbert function $H(R/I)=H$, satisfying $I_j=V'$, or,
equivalently, by Lemma \ref{AncandI}, satisfying
\begin{equation}\label{inclusion}
I+M^{j+1}\subset I'+M^{j+1} \text { and } I\cap M^j\supset I'\cap M^j.
\end{equation}
Let $N$ and $N'$ satsify the condition \eqref{eN} of Lemma \ref{NT} for a
fixed pair
$(d,j)$ and let
$I'$ be an ideal of Hilbert function $H(R/I')=N'$; then there is an ideal 
$I$ of Hilbert function
$H(R/I)=N$ satisfying $I\subset I'$. Likewise, let $T,T'$ satisfy the
condition \eqref{eT} of Lemma \ref{NT} and let $I'$ be an ideal of Hilbert
function
$H(R/I')= T$, then there is an ideal $I$ satisfying $H(R/I)=T$, and such
that
$I\supset I'$. 
\end{lemma}
\begin{proof}
Since $\dim I_j=\dim I'_j$ we have $I_j=I'_j$; thus we may prove the
result for $H$ by proving that for $N$ and $T$ separately. Our 
overall method is to construct a sequence of ideals $I'=I(0),I(1),\ldots
,I(s)=I$ of different Hilbert functions $H(R/I(u))=H(u)\in \mathcal
H(d,j)$ between
$H'=H(0)$ and
$H=H(s)$, using the properties of the $\tau$ invariant. \par We begin by 
considering a pair of Hilbert functions
$N\le N'$, each satisfying the condition relevant to $N$ in Lemma
\ref{NT}, and a given graded ideal $I'$ satisfying $H(R/I')=N'$. We will
construct an element of
$G(N)$, a graded ideal of Hilbert function $N$ satisfying $I\subset I'$.
We may assume that all the ideals considered contain $M^{j+1}$.
We first prepare to choose a Hilbert function
$N(1)$ of
$R/(I(1)$ differing from $N'$ in the highest possible degree. Then we will
determine the ideal
$I(1)\subset I'$. Let
$t<j$ be the largest integer, such that there is a permissible sequence
$N(1)$ for a level algebra in the sense of Lemma \ref{NT},
such that $N(1)_t\ne
N'_t$ and satisfying both
\begin{align}\label{interpolate}
  N'&\le N(1)\le N:  \text { that is } \forall i\le j \,\,  N'_i\le
N(1)_i\le N_i,\text { and }\notag\\
N(1)_i&=N'_i \quad \forall i\mid t<i\le j. 
\end{align}  
Let $E'=\Delta (N')$ be the difference sequence, and let $a$ be the
largest non-negative integer such that
\begin{equation*}
e'_t=e'_{t-1}=\cdots =e'_{t-a}.
\end{equation*}
{\it Claim a.} The sequence $N(1)$, defined by
\begin{equation}\label{minimal}
N(1)_i=
\begin{cases}
&N'_i\, \text{ unless }t-a\le i\le t\\
&N'_i+1 \, \text { for } t-a\le i \le t,
\end{cases}
\end{equation}
is a permissible sequence, in the sense that $N(1)$ satisfies \eqref{eN}
of Lemma
\ref{NT}.
Also, let $N''\ge N'$ termwise (so
$N''\le N'$ is a permissible
sequence for which
$\exists k,t-a\le k\le t$ with $N''_k\ne N'_k$). Then
$N''_i\ge  N(1)$. \smallskip\par\noindent
{\emph{Proof of Claim a.}} Because
$e'_i$ is non-increasing as
$i\le j$ decreases, the integer $t$ identifies the largest part
$e'_{t+1}\ne e_{t+1}$, and we have $e'_{t+1}<e_{t+1}$. By the definition
of $N(1)$ we have
$e(N(1))_i=e'_i$ unless $i=t+1$ or
$i=t-a$. We have 
\begin{equation*}
e(N(1))_{t+1}=e'_{t+1}+1\le e_{t+1}
\le
e_{t+2}=e'_{t+2}=e(N(1))_{t+2}.
\end{equation*}
and
\begin{equation*}
e(N(1))_{t-a}=e'_{t-a}-1
\ge e'_{t-a-1}=e_{t-a-1}.
\end{equation*}
Since both $N$ and $N'$ are permissible, the above inequalities shows that
$N(1)$ also is a permissible Hilbert function satisfying the condition
\eqref{eN} of Lemma
\ref{NT}.\par Suppose by way of contradiction that $N''$ is a permissible
sequence (for
$\LA (d,j)$ satisfying 
$N''\ge N'$ termwise, but not satisfying $N''\ge N(1)$, and let $u$ be the
smallest integer, $t-a\le u\le t$ such that $N''_u=N'_u$. If $t-a<u<t$
the difference $e''_u> e'_u=e'_{u+1}\ge e''_{u+1}$, contradicting the
assumption that $N''$ is permissible for $\LA (d,j)$. This completes the
proof of the Claim a.
\par
We now choose an ideal $I(1)\subset I'$ with $H(R/I(1))=N(1)$. Clearly
$I(1)_i=I'_i$ unless $t-a\le i\le t$, so we need only choose
$I(1)_{t-a},\ldots ,I(1)_t$. We construct $I(1)$ beginning with
\emph{lower} degrees. Suppose that
$u$ satisfies
$t-a\le u\le t$ and $I(1)_0,\ldots ,I(1)_{u-1}$ have been chosen so that
(here we regard $I(1)_u\subset R_u$)
\begin{equation*}
R_1\cdot I(1)_{v-1}\subset I'_v,\,\,  I(1)_v\subset I'_v, \text
{ and }\cod(I(1)_v)=N(1)_v\text { for } v<u.
\end{equation*}
Now  $R_1\cdot I_{u-1}\subset R_1\cdot I'_{u-1}\subset I'_u$, the first
inclusion by assumption, and the second since $I'$ is an ideal. We need to
choose a vector space
$I(1)_u$ between
$R_1\cdot I(1)_{u-1}$ and
$I'_u$, having codimension $N(1)_u$ in $R_u$. This is possible if and only if $\cod
(R_1\cdot I(1)_{u-1})
\ge N(1)_u$. We have
\begin{align*}
\dim R_1\cdot I(1)_{u-1}-\dim I(1)_{u-1}&=\tau (I(1)_{u-1})=\dim
I(1)_{u-1}-\dim R_{-1}\cdot I(1)_{u-1}\\
&\le \dim I(1)_{u-1}-\dim I(1)_{u-2}\,\,\text { by
\eqref{einclusion} }\\&=1+e_{u-1}(N(1))\\ &\le 1+e_u(N(1)), \text{ since
$N(1)$ is permissible }.
\end{align*}
Thus
\begin{align}
u+1-
\dim R_1\cdot I(1)_{u-1}&\ge u-\dim I(1)_{u-1}-e_u(N(1)),\notag\\
\cod R_1\cdot I(1)_{u-1}&\ge N(1)_{u-1}-e_u(N(1))= N(1)_u\label{N1key}
\end{align}
by our choice of $N(1)$. Therefore, we may choose $I(1)_u$ such that
$I'_u\supset I(1)_u\subset R_1\cdot I(1)_{u-1}$, satisfying $\cod
I(1)_u=\cod I'_u+1$. Continuing this process, we may choose
an ideal $I(1)\subset I(0)=I'$ of  Hilbert function $H(R/I(1))=N(1)$, as
claimed. Continuing in this manner, we eventually construct $I(s)$ of
Hilbert function
$H(R/I(s))=N(s)=N$, and satisfying
$I(s)\subset I'$, as claimed. This completes the proof of the Lemma for
the pair
$(N,N')$.\smallskip
\par
We now turn to choosing an ideal $I$ of Hilbert function $H(R/I)= T$
given $I'$ satisfying $H(R/I')=T'$. Although proof of this portion of the
Lemma involving
$\GA_T(d,j)$ for
$T,T'$ eventually zero  appears already in
\cite[Section 4B]{I2}, we include the argument with further details here
for completeness.  For now we assume that $T,T'$ are eventually zero:
that $c_T=c_{T'}=0$. We will also now assume that our ideals $I\subset
M^j$, by intersecting with
$M^j$ if necessary. We first choose the Hilbert function
$T(1)$ of
$R/(I(1))$, differing from $T'$ in the lowest degree possible, and then
the corresponding ideal
$I(1)$.
 \par Let
$t>j$ be the smallest integer, such that there is a permissible sequence
$T(1)$ satisfying the condition \eqref{eT} of Lemma \ref{NT} for $T$, and
such that $T(1)_t\ne
T'_t$ and satisfying both
\begin{align}\label{interpolateb}
  T'&\ge T(1)\ge T:  \text { that is } \forall i\ge j \,\,  T'_i\ge
T(1)_i\ge T_i,\text { and }\notag\\
T(1)_i&=T'_i \quad \forall i\mid j \le i<t. 
\end{align}  
Let $E'=\Delta T'$ be the difference sequence, and let $a$ be the largest
non-negative integer such that
\begin{equation}\label{ee'T}
e'_{t+1}=e'_{t+2}=\cdots =e'_{t+a}.
\end{equation}
{\emph {Claim b}}. The sequence $T(1)$, defined by
\begin{equation}\label{minimalb}
T(1)_i=
\begin{cases}
&T'_i\, \text{ unless }t\le i\le t+a-1\\
&T'_i-1 \, \text { for } t\le i \le t+a-1,
\end{cases}
\end{equation}
is a permissible sequence satisfying the condition \eqref{eT} of
Lemma
\ref{NT}. Furthermore, let $T''\le T'$ (termwise) be a permissible
sequence for which
$\exists k,t< k\le t+a$ with $T''_k\ne T'_k$. Then
$T''\le T(1) $. \smallskip\par\noindent
{\emph{Proof of Claim b}}. Because
$e'_i$ is non-increasing as
$i\ge j$ increases, the integer $t$ identifies the largest difference
$e'_{t}\ne e_{t}$, and we have $e'_i=e_i$ for $i$ satisfying $i\le
t-1$. Since $T'_{t}>T_{t}$,
we have $e'_{t}=T'_{t}-T'_{t-1}> T_{t}-T_{t-1}=e_{t}$ so we have
$e'_{t}>e_{t}$.
Evidently $e(T(1))_i=e'_i$ unless $i=t$ or
$t+a$. We have 
\begin{align*}
e(T(1))_{t}&=e'_{t}+1\le e_{t}\\
&\le
e_{t-1}=e'_{t-1}=e(T(1))_{t-1}.
\end{align*}
and
\begin{align*}
e(T(1))_{t+a}&=e'_{t+a}-1\\
&\ge e'_{t+a+1}=e(T(1))_{t+a+1}.
\end{align*}
Since both $T$ and $T'$ are permissible, the above inequalities show that
$T(1)$ also is a permissible sequence --- one satisfying the condition
\eqref{eT} of Lemma
\ref{NT} for $T$.\par
 Suppose by way of contradiction that $T''$ is likewise a permissible
sequence satisfying
$T''\le T'$ termwise, but $T''$ does not satisfy $T''\le T(1)$, and let
$u$ be the smallest integer, $t\le u\le t+a$ such that $T''_u=T'_u$. If
$t<u<t+a$ the difference $e''_u< e'_u=e'_{u+1}\le e''_{u+1}$,
contradicting the assumption that $T''$ is permissible for $\GA (d,j)$.
This completes the proof of the Claim b.\smallskip
\par
We now choose an ideal $I(1)\supset I'$ with $H(R/I(1))=T(1)$, beginning
with the \emph{higher} degrees. Clearly
$I(1)_i=I'_i$ unless $t\le i\le t+a-1$, so we need only choose
$I(1)_{t},\ldots ,I(1)_{t+a-1}$. Suppose that $u$ satisfies
$t+1\le u\le t+a$ and $I(1)_{u+1},\ldots ,I(1)_{t+a}$ have been chosen so
that
\begin{equation*}
 R_{-1}\cdot I(1)_{v+1}\supset I'_v,\,\,  I(1)_v\supset I'_v, \text
{ and }\ \cod I(1)_v=T(1)_v\text { for } v>u.
\end{equation*}
Now  $R_{-1}\cdot I(1)_{u+1}\supset R_{-1}\cdot I'_{u+1}\supset I'_u$, the
first inclusion is by assumption, and the second since $I'$ is an ideal.
We need to choose a vector space
$I(1)_u$ between
$R_{-1}\cdot I(1)_{u+1}$ and
$I'_u$, having codimension $T(1)_u$ in $R_u$. This is possible if and only if $\cod
(R_{-1}\cdot I(1)_{u+1})
\le T(1)_u=T'_u-1.$ We have
\begin{align*}
\dim I(1)_{u+1}-\dim R_{-1}\cdot I(1)_{u+1}&=\tau (I(1)_{u+1})=\dim
R_1\cdot I_{u+1}-\dim I_{u+1}\\ &\le \dim I(1)_{u+2}-\dim
I(1)_{u+1}\,\,\text { by
\eqref{einclusion} }\\&\le 1+e_{u+2}(T(1))\\ &\le 1+e_{u+1}(T(1)), \text{
since
$T(1)$ is permissible }.
\end{align*}
Thus
\begin{align}
u+1-
\dim R_{-1}\cdot I(1)_{u+1}&\le u+2-\dim I(1)_{u+1}+e_{u+1}(T(1)),\notag\\
\cod R_{-1}\cdot I(1)_{u+1}&\le T(1)_{u+1}+e_{u+1}(T(1))=
T(1)_u\label{T1key}
\end{align}
by our choice of $T(1)$. Therefore, we may choose $I(1)_u$ such that
$I'_u\supset I(1)_u\subset R_{-1}\cdot I(1)_{u+1}$, satisfying $\cod
I(1)_u=\cod I'_u-1$.  Continuing this process, we may choose
an ideal $I(1)\supset I(0)=I'$ of  Hilbert function $H(R/I(1))=T(1)$, as
claimed. Continuing in this manner, we eventually construct $I(s)$ of
Hilbert function
$H(R/I(s))=T(s)=T$, and satisfying
$I(s)\supset I'$, as claimed. This completes the proof of the Lemma for
the pair $(T,T')$ when 
$c_T=c_{T'}=0$. \par
When $c_T\not=0$, by Corollary \ref{commonfactor} any ideal $I$ with
$H(R/I)=T$ must have a common factor $f=\GCD (I)$ of degree $c_T$. We
have 
$T\le T'\Rightarrow c(T)\le c(T')$. Suppose the pair of ideals $I,I'$
satisfies
$I\supset I', H(R/I)=T,H(R/I')=T'$, then
 $f=\GCD (I)$ divides any common factor
$f'=\GCD (I')$ of
$I'$. Given $I'$, we now refine the choice of $I$ by choosing in advance a
degree
$c(T)$ factor
$f$ of
$\GCD (I')$ to be the common factor of $I$. Now it will suffice to choose
$J=I:f$ of Hilbert function $T:c(T')$ containing
$I':f$, of Hilbert function $ T':c(T')$, and then set $I=fJ$. Thus
we have reduced to showing the Lemma when $T$ is eventually zero, but
$c_{T'}>0$.\par
Suppose now that $c_T=0, c'=c_{T'}\not= 0$, and define $s'$ by $
T'_{s'-1}>T'_s=c_{T'}>0$. (When no such integer $s'$ exists, then
$I'=(f')$ and choosing $I\supset (f')$ poses no difficulty). Let
$f'$ be the degree $c'$ common factor of
$I'$. When $e'_i$ of \eqref{ee'T} satisfies $e'_i>0$ we choose
$T(1) $ as in the case $c_T=c_{T'}=0$, however to construct
$I(1)$, we first construct
$I(1):f'$ of Hilbert function $T(1):c'$ such that $I(1):f'\supset I':f'$,
as above, then we let $I(1)=f'\cdot (I(1):f')$. When $i=s'+1$ and $e'_i=0$
in
\eqref{ee'T}, then $a=+\infty$ in \eqref{ee'T}. We choose $I(1)\cap
M^{s+1}=(f'_1)\cap M^{s+1}$ with $f'_1$ a degree $c'-1$ divisor of $f'$.
Continuing in this way, we obtain finally an ideal $I\supset I'$ of
Hilbert function $H(R/I)=T$. This completes the proof of the statements
involving $T,T'$ of the Lemma in all cases.\par We now turn to the case
of a pair $H,H'$ of acceptable Hilbert functions. When
$H$ is eventually zero, one uses the above methods to first construct
$I+M^{j+1}$ and then construct
$I\cap M^j$, which together determine the ideal $I$ (since $I_j=I'_j$ is
given). When
$H$ is eventually
$c$, then one chooses $f$ of degree $c$ dividing the common factor $f'$ of
$I'$ of degree $c(T')\ge c$. Then one chooses $I:f$ of Hilbert function 
$T:c$, as above from $I':f$ of Hilbert function $T':c$, then sets
$I=f\cdot (I:f)$. Since $H=H(N,T)$ is acceptable (Definition
\ref{defaccept}) if and only if $N, T$ have the same $\tau$ and are
both permissible (satisfy \eqref{eN} or \eqref{eT}, respectively), this
completes the proof of the Lemma.
\end{proof}
\begin{example}We illustrate the process of choosing $N(1)$ in the proof
above. Suppose that the two sequences $N',N$ are $N'=(1,2,\ldots
,13,11,9,7,4,0)$ with $N'_{16}=4$, and $N=$ \linebreak $(1,2,\dots
,13,12,11,8,4,0)$. We choose $N(1)$: here $t=15$, and one chooses
$N(1)_{15,16}=(8,4).$ However, if this were the only change, the
intermediate sequence
$(1,\ldots ,13, 11,9,8,4,0)$ would violate the condition on first
differences, as it has first differences  $(\ldots 2,1,4,4)$ , which has
a decrease from 2 to 1. Instead, we must choose $N(1)=(1, \ldots
,13,12,10,8,4)$, which is also next to
$N'$ in the partial order among the subset of sequences possible for level
algebras
$\LA (13,16)$ and having $N(1)_{15}>7$. Then 
$N(2)=N$.  Note that $N(0)=(1,\ldots ,13,12,10,7,4,0)$ is next to $N'$ in
the partial order, but we have chosen to step to $N(1)$, which is the
closest to
$N'$ among those between $N'$ and $N$ and differing from $N'$ in the
highest possible degree.
Note that in the proof of Lemma~\ref{build}, the occuring Hilbert
functions $N(i), T(i)$ must be permissible for a level algebra, graded
ideal, respectively of a vector space of forms. But the intermediate
ideals
$I(1),
\ldots $ that we construct are not themselves level ideals, nor ideals
generated by
$I_j$, respectively.
\end{example}
 Recall from Definition \ref{defpo} that we denote by $\mathcal P=\mathcal P(d,j)$
the partial order on the set $\mathcal H(d,j)$ of acceptable Hilbert functions.
The acceptable Hilbert functions are described in Definition \ref{defaccept}, and
further in Lemma \ref{accept}. Recall that we showed in Theorem \ref{allH} that these
$H\in \mathcal H(d,j)$ are exactly the sequences occurring as Hilbert functions of
ancestor algebras.
\begin{theorem}\label{closureofstrata} Let $d,j$ be positive integers satisfying
$d\le j$, assume that the field $k$ is algebraically closed, and suppose
that
$H$ is an acceptable
$O$-sequence (Definition
\ref{defaccept}).
\begin{enumerate}[A.]
\item{\sc {Frontier property.}}\label{closureA} The Zariski
closure $\overline {\Grass _H(d,j)}$ satisfies 
\begin{equation}\label{eclosure}
\overline {\Grass
_H(d,j)}=
\bigcup_{H'\ge_{\mathcal P} H} \Grass_{H'}(d,j).
\end{equation}
The analogous equality holds for $\overline{\LA_N (d,j)}$ and for
$\overline{\GA_T (d,j)}$.
\item{\sc $G(H)$ is a desingularization of $\overline{\Grass
_H(d,j)}$.}\label{closureB} There is a surjective morphism $\pi :G(H)\to
\overline{Grass(H)}$ from the nonsingular variety $G(H)$, given by
$I\to I_j$. The inclusion $\iota :\Grass_H(d,j)\subset G(H), \iota :
V\to \overline{V}$ is a dense open immersion. For $H'\in \mathcal
H(d,j), H'\ge_\mathcal P H$, the fibre of $\pi$ over
$V'\in
\overline{\Grass_H(d,j)}\cap
\Grass_{H'}(d,j)$ parametrizes the family of graded ideals 
\begin{equation}\label{IandV'}
\{ I\mid H(R/I)=H \text { and } I_j=V'\}.
\end{equation}
The schemes $\overline{\LA_N(d,j)}$ and $\overline{\GA_T(d,j)}$ have
desingularizations $G(N)$ and $G(T)$, respectively, with analogous
properties.
\end{enumerate}
\end{theorem}
\begin{proof} By Theorem \ref{basic} \eqref{basici},\eqref{basiciii} $G(H)$ is
nonsingular and has as open dense subset the subfamily of ideals with minimum
number of generators; by Proposition
\ref{ancideals} \eqref{ancidealsv}, this subfamily is $\iota (\Grass_H(d,j))$ (see
also Theorem \ref{dimstrata}\eqref{dimstrataii}). By definition of $\pi$ the fibre of
$\pi$ is the family specified in
\eqref{IandV'}. That $\pi$ is surjective we will show next, thus completing the
proof of (B). \par
We now show \eqref{eclosure}. Suppose that $H'\ge H\in
\mathcal H(d,j)$: so
$H,H'$ satisfy the condition of Proposition \ref{etau} and each occurs as the Hilbert
function of an ancestor ideal, and let $V'\in \Grass_{H'}(d,j)$. By Lemma
\ref{build} there is an ideal $I$ of Hilbert function $H$ satisfying
$I_j=V'$. Since $G(H)$ is irreducible with open dense subscheme
$\Grass_H(d,j)$ we have that there is a family $I(t), t\in \Z$ of ideals parametrized
by  a curve $\Z \subset G(H)$ such that for $t\not= t_0$, $I(t)\in 
\iota (\Grass_H(d,j))$, with
$I=\lim_{t\to t_0} I(t)$; it follows that $V'=\lim_{t\to t_0} V(t)=(I(t)_j$ is in
the closure of $\Grass_H(d,j)$. This shows that the closure
$\overline{\Grass_H(d,j)}$ includes the union of lower strata in \eqref{eclosure}. By
Theorem
\ref{loclostrata} the closure 
$\overline {\Grass_H(d,j)}$ is a subset of $
\bigcup_{H'\ge_{\mathcal P}H }\Grass_{H'}(d,j)$.
This completes the proof of \eqref{eclosure} and (A), as well as (B) for
$\overline{\Grass_H(d,j)}$.
 An analogous argument proves
the results in (A) concerning the closures 
 $\overline{\LA_N (d,j)}$ and
$\overline{\GA_T (d,j)}$.  This
completes the proof.
\end{proof}
\begin{corollary}\label{grasstauirred}
The scheme
$\Grass_\tau (d, R_j)$ is irreducible and 
$\Grass_{H_\tau}(d,j)$ (see \eqref{eHtaub}) is a dense open subscheme. The
Zariski closure of
$\Grass_\tau (d,j)$ satisfies
$\overline{\Grass_\tau )(d,j)}=\bigcup _{\tau '\le \tau} \Grass_{\tau
'}(d,j)$.
\end{corollary}
\begin{proof}
We fix $(d,j,\tau )$. Evidently, by Lemma \ref{taubasic}\eqref{taubasicii}
and equation \eqref{HAdual},  the Hilbert function $N(H_\tau)$ is
maximum, among the Hilbert functions $N(H)$ for $H$
satisfying
$\tau (H)\le
\tau$. Similarly
\eqref{BandH} and \eqref{BandT} show that $T(H_\tau)$ has the
minimum values among such $H$. Thus, Theorem \ref{closureofstrata} implies
the Corollary.
\end{proof}\par
\begin{definition}\label{defpos}
We denote by $\mathcal PA(d,j)$ the partially ordered set of pairs of
partitions $(P,Q)$ such that $P $ partitions $d$, $Q $
partitions an integer no greater than $j+1-d$, and the largest part $p_1$ of
$P$ and the largest part $q_1$ of $Q$ satisfy $p_1=q_1+1$. We let
$(P,Q)\le (P',Q')$ if both
$P\le P'$ and $Q\le Q' $ in the respective
majorization partial orders.
\end{definition}
\begin{theorem}\label{equivpos}  There is an isomorphism of partially
ordered sets $\mathcal H(d,j)$ under the partial order $\mathcal P(d,j)$
and the partially ordered set $\mathcal PA(d,j)$, under the
product of the majorization partial orders (see Definition \ref{defpos})
given by
$H\to (P,Q),P=P(H)=A(H)^\ast, Q=Q(H)=B(H)^\ast$ (see Definitions \ref{defP} and
\ref{defA-D}). This is the same order as is
induced by specialization (closure) of the strata $\Grass(H)$.\end{theorem}
\begin{proof}
This is immediate from \eqref{eclosure}, Theorem \ref{allH} \eqref{allHiii}, and
Lemma \ref{comparepart}.
\end{proof}
\begin{example}\label{exsimple}
 We consider the partial order on all sequences $H$ for
$(d,j)=(4,5)$ (see Table \ref{table45}). Thus, $A$ partitions the
dimension
$d=4$ into
$\tau
\le 3$ parts, and $B$ partitions the integer $\cod (V)-c=2-c$ into $\tau
-1$ parts. $\Grass(4,R_5)$ has dimension 8; the open cell is given by the
pair
$A=(2,1,1), B=(1,1)$. When $\tau=2$ there are two
sequences, and for $\tau=1$ a single sequence. They are here linearly ordered
by $\ge_{\mathcal P(4,5)}$, so by Theorem \ref{closure} the closure of
each stratum listed in Table~\ref{table45} is the union of the stratum
itself with the strata below it. Note that the $A, P$ and $Q$ columns of
partitions in Table~\ref{table45} are simply ordered in the majorization partial
order, but the
$B$ column is not. The order on $\mathcal H(d,j)$ is equivalent to the
product of majorization orders on the pairs $(P,Q)$.
\end{example}

\begin{table}[tbh]
\begin{center}
\begin{tabular}{|l| r r r  r r r l l |}
\hline
Stratum &   $\tau$&   $A$&  $B$& $P=A^\ast$ &$ Q=B^\ast$ & c&$\cod $ &   $H$
\\
\hline 
   H(0) &    3  & (2,1,1)& (1,1) &(3,1)&(2)  & 0  & 0  &
(1,2,3,4,4,2,0,\underline{0})\\
   H(1) &    2&   (2,2) & (2) & (2,2) &(1,1)  &  0 &
1&(1,2,3,4,3,2,1,\underline{0}) \\ 
   H(2) &    2&   (2,2) &(1) & (2,2)&(1) & 1 &  
3&$(1,2,3,4,3,2,\underline{1})$
\\
   H(3)&     1 &  (4) & - &(1,1,1,1)& -&   2 & 6  &$(1,\underline
{2})$
\\
\hline
\end{tabular}
\caption{Hilbert functions $H $ for $(d,j)=(4,5)$}\label{table45}
\end{center}
\end{table}
\begin{remark} Possibly relevant to the frontier property, given Theorem
\ref{closureofstrata}\eqref{closureA} and Theorem~\ref{equivpos},
C.~Greene and D.~J.~Kleitman
have studied the longest simple chains in the lattice of
partitions of an integer \cite{GrK}. \par
Relevant to the desingularization of Theorem \ref{closureofstrata}
\eqref{closureB}, a basis for the homology of
$G(H)$ is given in
\cite{IY}, in terms of the classes $\pi_\ast (E(J))$ determined by the
monomial ideals
$J$ of Hilbert function
$H(R/J)=H$: here $E(J)$ is the affine cell parametrizing graded ideals
having initial ideal
$J$, and it the set $\{ E(J)\}$ form a cell decomposition of $G(H)$. 
A natural cobasis of a monomial ideal of colength
$n, H(R/J)=H$ is a vector space  $E^c (J)$ of monomials whose graph is
the Ferrers graph of a partition $P(E^c)$ of $n$ with
diagonal lengths
$H$. The dimension of the cell $E(J)$ is the number of
difference one hooks (arm-leg=1) in the partition $P(E^c)$ When
$\mid H\mid=\sum H_i=n$ a basis for the degree-$i$ homology corresponds
one-to-one with the  partitions of n having the given diagonal lengths
$H$; and having the given number $i$ of hooks of difference one. In a few
cases the homology ring structure of
$G(H)$ is known, but in general the homology ring structure is not known
(see
\cite{IY}).
\end{remark}
\section{Waring problem, related vector spaces}\label{relatedapply}

 In 
Section
\ref{subsimultWaring} we apply the previous results
to a refinement of the simultaneous Waring problem for a vector space of
forms. In Section \ref{related} we first return to polynomial rings $R$ of
arbitrary dimension $r$, to develop the notion of a
space $W\subset R_i$ related to a vector space $V\subset R_j$ if
$W$ is obtained by a chain whose elements are each a homogeneous
component of the ancestor ideal of the predecessor space. When $r=2$ we
bound the number of classes $\overline{W}$ related to $\overline{V}$
in terms of the $\tau$ invariant $\tau(V)$. Finally, we state some open
problems. 
\subsection{The simultaneous Waring problem for degree-$j$ binary
forms}\label{subsimultWaring}
 We let $r=2$ and denote by $\mathcal
R=k[X,Y]$ the \emph{dual} polynomial ring to $R$. We suppose that $\cha
k=0$ or $\cha k =p>j$ throughout this section. The \emph{simultaneous
Waring problem} is to find the minimum
number
$\mu (c,j)$ of linear forms, needed to write each element of a general
dimension-$c$ vector space $\mathcal W\subset \mathcal R_j$ as a sum
of $j$-th powers of the linear forms; here the choice of the linear forms
depends on
$\mathcal W$. Our refinement is to fix also the differential $\tau$
invariant of $\mathcal W$. \par  The case
$c=1$ of a single binary form
$F$ is quite classical: it is related to the secant varieties of rational
normal curves, and is resumed along with this connection in
\cite[\S 1.3]{IK}. Note that in this section \ $c=\dim \mathcal W$
satisfies $c=\cod (V)=j+1-\dim V$ where $ V=(\Ann \mathcal W)_j$ (see
\eqref{annW} below). Letting $\mu(W)$ denote the minimal length of a
simultaneous (generalized) additive decomposition of $W$, our results rest
on the identity $\mu(W)=\mu (L(V))$, the order of the level ideal
$L(V)$ determined by $V$ (Lemma \ref{lWaring}), valid for $r=2$ only.
 For
$u\le c$ we let
$c_a=c(c-1)\cdots (c+1-a)$.
\begin{definition}
The ring $R=k[x,y]$ acts on $\mathcal R$ by differentiation
\begin{equation}
x^ay^b\circ X^cY^d=
\begin{cases}& (c_a\cdot d_b) X^{c-a}Y^{d-b}  \text{ if }  c\le a \text {
and } b\le d\\
&0 \quad \text { otherwise }.
\end{cases}
\end{equation}
Let $V\subset R_j$ be a vector subspace. We denote by $V^\perp\subset
\mathcal R_j$ the subspace 
\begin{equation}
V^\perp=\{ F\in \mathcal R_j\mid v\circ
F=0\,\, \forall v\in V\} .
\end{equation}\label{Vperp}
Given $\mathcal W\subset \mathcal R_j$ we denote by $\Ann (\mathcal
W)\subset R$ the ideal
\begin{equation}\label{annW}
\Ann\mathcal W =\{ f\in R\mid f\circ w=0\,\, \forall w\in \mathcal W\} .
\end{equation}
Let $V=(\Ann(\mathcal W))_j\subset R_j$. We define the \emph{differential
$\tau$-invariant}
$\tau_\delta (\mathcal W)$ as 
\begin{equation}\label{tauW}
\tau_\delta (W)=\tau (V)=\dim R_1\cdot V-\dim V.
\end{equation}
 We need also the following notions of
additive decomposition: let $F\in \mathcal W$ then $F=\sum_{i=1}^s \alpha
_i L_i^j$ is an additive decomposition of length $\mu$ of $F$, assuming that
the $\{ L_i\}$ are pairwise linearly independent. The form $F\in \mathcal
R_j$ has a \emph{generalized additive decomposition} (GAD) of
 length $\mu$ and weights $\beta_1,\ldots ,\beta_t$ into powers of the
linear forms
$L_1,\ldots ,L_t\in \mathcal R_1$ if 
\begin{equation}\label{GAD}
 F=\sum_{i=1}^t G_iL_i^{j+1-\beta_i} \text { where } \deg G_i= \beta_i-1
\text { and }
\sum
\beta_i=\mu. 
\end{equation}
The vector space $\mathcal W\subset \mathcal R_j$ has a
\emph{simultaneous decomposition} of length $\mu$ if there is a single
ordered set $L=(L_1,\ldots ,L_t)$ of linear forms $L_i\in R_1$
(which may depend on W) and weights
$\beta=(\beta_1,\ldots ,\beta_t)$ such that each $F\in W$ has a GAD of
length $\mu$ and weights
$\beta$ into the  forms $L$. We denote by $\mu(\mathcal W)$ the shortest
length of a simultaneous additive decomposition of $W$.\par
We define $\mu (c,j), \mu (\tau ,c,j)$, respectively, as the common value
of
$\mu (W)$ for
$\mathcal W$ in a suitable open dense subset of $\Grass (c, \mathcal
R_j)$, or of $\Grass_{\tau_\delta} (c,R_j)$ (where $\tau_\delta (\mathcal
W) =\tau$), respectively.
\end{definition}
Note that we defined $\tau_\delta (\mathcal W )$ for $ \mathcal W\subset
\mathcal R_j$ using the annihilating degree-$j$ space
$V=(\Ann ({\mathcal W}))_j$. Here is a direct definition. Let $R_1\circ
\mathcal W\subset
\mathcal R_{j-1}$ be
$R_1\circ \mathcal W=\{
\ell\circ w,
\ell\in R_1, w\in \mathcal W\}$. Letting
$N=(n_0,n_1,\ldots )=H(R/\Ann({\mathcal W}))$, we have from
 $(\Ann({\mathcal W})_{j-1})^\perp=R_1\circ
\mathcal W$ and \eqref{tauseq2}
\begin{align}
\tau_\delta (\mathcal W)&=1+e_j(N)=1+n_{j-1}-n_j\notag\\
&=1+\dim R_1\circ \mathcal W - \dim \mathcal W.
\end{align}
\par\noindent  For $L_i=a_iX+b_iY\in \mathcal R_1$ we let
$\ell_i=b_ix-a_iy\in
\mathcal R_1$: then $\ell_i\circ L_i=0$. We have the
following well-known result. Recall that $\mu (L(V))$ is the \emph{order}
of the level ideal 
$L(V)$.
\begin{lemma}\label{lWaring} Let $V\subset R_j$ and set $\mathcal
W=V^\perp$. The level ideal
$L(V)$ satisfies
\begin{equation}\label{levelann}
L(V)=\Ann (\mathcal W), \,\mathcal W=V^{\perp}.
\end{equation}
Let $F\in \mathcal R_j$. Then $F$ has a GAD of length $\mu$ as in
\eqref{GAD} if and only if \begin{equation}\label{annF}
\exists f\in\Ann (F) \text{ such that } \deg f=\mu \text { and } f=\prod
\ell_i^{\beta_i}, \ell_i\in R_1 
\end{equation} 
Let $\mathcal W\subset \mathcal R_j$ and $\dim \mathcal W =c$. Then
$\mu(\mathcal W)=\mu(L(V))$ for $ V=(\Ann (\mathcal W )_j$. Also $1\le
{\tau_\delta}$ and
\begin{equation}\label{taumax}
{\tau_\delta}(\mathcal W)\le \min\{
c+1,j+1-c\},
\end{equation}
with equality in \eqref{taumax} for a generic choice of $\mathcal
W\subset R_j$ of dimension $c$.
\end{lemma}
\begin{proof}
The identity \eqref{levelann} is a basic property of inverse systems -
see in general
\cite[\S 60ff]{Mac1}\cite{EmI,Ge} or for a modern proof, \cite[Lemma
2.17]{IK}. Equation \eqref{annF} is \cite[Lemma 1.33]{IK}; that
$\mu(W)=\mu(L(V))$ is a straightforward consequence. The last statement is
a consequence of the upper bound on
$\tau (V), V= (\Ann W)_j$ from Lemma \ref{taugen}, rewritten in terms of
$c,j$, since $\tau_\delta (\mathcal W)= \tau (V)$.
\end{proof}\par
We let $c=j+1-d$ and define $
\mu (\tau ,d,j)= j+1-\lceil d/\tau\rceil$.
When  $\mu\le \mu(\tau ,d,j)$, we define the 
Hilbert function sequence
$N(\mu,\tau ,d,j )$ by
\begin{equation}\label{eNs}
N(\mu,\tau,d,j)_i=
\begin{cases}
&\min\{i+1,\mu ,c+(\tau -1)(j-i)\} \text { for }  i\le j\\
& 0 \text { for } i>j.
\end{cases}
\end{equation}
We define $N(\tau,d,j)=N(H_\tau (d,j) )$ with $H_{\tau} (d,j)$ from
equation
\eqref{eHtaub}:  thus we have \linebreak $N(\tau,d,j)_i=\min\{i+1 ,c+(\tau
-1)(j-i)\}$ for
$ i\le j$. We define $a,\kappa\in \mathbb N$ by $\mu-c=a(\tau-1)+\kappa$
with
$0\le \kappa =\rem (\tau-1,\mu-c)<\tau -1$,
\begin{lemma}\label{Ndj} $N(\tau,d,j)$ is the maximum level algebra
Hilbert function for a $d$ dimensional vector space $V\subset R_j$ with
$\tau(V)=\tau$; it has order 
$\mu(\tau,d,j)$ and
partition $P(\tau,d,j)=(\tau^{\lfloor d/\tau\rfloor},\rem (\tau,j))$
from
\eqref{eHtau}.
$N(\mu,\tau,d,j)$ is the maximum level algebra Hilbert function that is
both bounded above by $\mu$ and possible for a vector space $V\subset R_j$
with
$\tau(V)=\tau$. It has order $\mu$ and partitions $P,A$ of $d$
\begin{align}
P=P(\mu,\tau,d,j)&=(\tau^a,\kappa  +1,1^{j-\mu -a}),\label{eqPmutau}\\
A=A(\mu,\tau,d,j)&=P^\ast=(j+1-\mu ,\lceil (\mu -c)/(\tau
-1)\rceil^{(\kappa  -1)^+},a^{\tau-\kappa }).\label{eqAmutau}
\end{align}
 The dimension of $\LA_N(d,j),
N=N(\tau,d,j)$ is $\tau (j+2-\tau)-d$. 
\end{lemma}
\begin{proof}  The order $\mu=\mu (\tau,d,j)$ of $N(\tau ,d,j)$ satisfies
\begin{equation*}
\mu =\max\{ i\bigm| N(\tau ,d,j)_{i-1}\ge i\} =\max\{ i\bigm|
c+(j-(i-1))(\tau -1)\ge i\},
\end{equation*}
which leads to $\mu = \mu(\tau ,d,j)$. The calculation of $P(\mu ,\tau
,d,j), A(\mu,\tau,d,j)$ is routine, and the dimension formula for $\LA_N(d,j),
$ is 
\eqref{edimtau}.
\end{proof}\par
 One part \eqref{sWaringii} of the following
Theorem may be classical; it was shown by J. Emsalem and the author in an
unpublished preprint, and also in
\cite{Ca,CaCh}. 
\begin{theorem}\label{sWaring}
We will suppose that $\mathcal W\subset \mathcal R_j, \mathcal R=k[X,Y],\,
\dim \mathcal W =c$, and $d=j+1-c$.
\begin{enumerate}[i.]
\item\label{sWaringi}
Each dimension $c$ subspace $\mathcal
W\subset \mathcal R_j$ with $\tau_\delta (\mathcal W)=\tau $ satisfies
$c\le \mu (W)\le \mu (\tau,d,j)$, with equality $\mu (W)=\mu (\tau,d,j)$
for a generic choice of such $\mathcal W$.
\item\label{sWaringii} For general $\mathcal W$ the value of $\mu
(\mathcal W)$ is $\lfloor c(j+2)/(c+1)\rfloor$ if $c<j/2$, and $j$
otherwise.
\item\label{sWaringiii} Let $c\le \mu\le \mu(\tau ,d,j)$. When $k$ is
algebraically closed, the subfamily
$\GAD_\mu (\tau ,c,j)$ of\linebreak
$\Grass_{\tau_\delta}(c,\mathcal R_j)$ parametrizing $\mathcal W$
satisfying $\tau_\delta (\mathcal W)=\tau$ and $\mu (\mathcal W)\le \mu$
is isomorphic under
$\mathcal W\to (\Ann\mathcal W)_j$ to $\overline{\LA_{N}(d,j)}$, 
where $N=N(\mu,\tau ,d,j)$. The
codimension of
$\LA_N(d,j)$ in
$\Grass_{\tau}(d,j)$ satisfies, for $1\le \mu < \mu (\tau,d,j)$
\begin{equation}\label{codmutau}
\cod_{\tau_\delta} \GAD_\mu (\tau ,c,j)=\ell (A)=(j-\mu )\tau -(d+1).
\end{equation} 
\end{enumerate}
\end{theorem}\noindent
\begin{proof} By Lemmas \ref{lWaring} and \ref{Ndj} each of the statements
\eqref{sWaringi},\eqref{sWaringii}, and the first part of
\eqref{sWaringiii} translates into one about the order of $N(\tau,d,j)$,
or the dimension of $N(\mu,\tau,d,j)$.
 Corollary \ref{grasstauirred} implies that for an
open dense set of
$V\in\Grass_\tau(d,j)$, the Hilbert function of $\LA (V)$ is $N(\tau
,d,j)$, derived from $H(\tau,d,j)$ of
\eqref{eHtaub}. Thus, the order $\mu (\tau,d,j)$ of $ N(\tau ,d,j)$, is
the generic value for $\mu (W), W, \tau_\delta (W)=\tau$. This gives
\eqref{sWaringi}, and \eqref{sWaringii} follows from substituting
$\tau=c+1$ or $j+1-c$ from \eqref{taumax} into the formula of
\eqref{sWaringi}.  The codimension of
$\LA_N(d,j)$ in $\Grass_\tau (d,j)$ of \eqref{sWaringiii} is by
\eqref{codb} the invariant
$\ell (A)$ of
\eqref{ellA} for the partition $A=A(\mu,\tau,d,j)$ from \eqref{eqAmutau};
however a routine calculation using $\dim N(\tau ,d,j)$ from Lemma
\ref{Ndj} and \eqref{edimN} --- assuming $e_\mu=0$ for $N=N(\mu,\tau
,d,j)$ --- gives 
\eqref{codmutau} for $\mu < \mu (\tau ,d,j)$ (when $\mu=\mu (\tau d,j)$
the assumption $e_\mu=0$ for \eqref{codmutau} may not hold).
Theorem~\ref{closureofstrata} completes the proof of
\eqref{sWaringiii}.
\end{proof}
\begin{remark} Theorem \ref{sWaring} states that vector spaces
$\mathcal W$ with higher $\tau$ in general require a larger number of
linear forms $L_1,\ldots ,L_\mu $ so that 
\begin{equation}\label{edistinct}
\mathcal W\subset \langle L_1^j,\ldots
,L_\mu ^j\rangle.
\end{equation}
 Thus, letting $V=(\Ann (W))_j$ when $\tau (V)=1$ so
$V=f_c R_{j-c}$, we have $\mu (\mathcal W)=c$. When $c\ge j/2$ and $\tau
(V)=j+1-c$, the maximum value, then $\mu (\mathcal W)=j$ in general. Note
that, given $(\mu, \tau, d,j)$ satisfying $c\le\mu\le \mu (\tau,d,j)$, the
proof of Theorem \ref{basic} in \cite{I2} shows that one can choose a
vector space
$V\in \LA_N(d,j), N=N(\mu,\tau,d,j)$ such that there is a form $f\in
L(V)_\mu$ with distinct roots, thus one may suppose that a general
$\mathcal W\in \GAD_\mu (\tau ,c,j)$ satisfies \eqref{edistinct}
\end{remark}
\subsection{Vector spaces related to $V$; open problems}\label{related}
 In section \ref{related} the dimension $r$ of $R$ is arbitrary
unless otherwise specified. We say that $W\subset R_i$ is related to
$V\subset R_j$ if there is a sequence $(i_1,\ldots ,i_k)\in \mathbb Z^k$
such that 
\begin{equation}\label{erelationexp}
W=R_{i_k}\cdot R_{i_{k-1}}
\cdots  R_{i_1}V= R_{i_k}\cdot (R_{i_{k-1}}\cdot ( \cdots
R_{i_1}V)\ldots  ).
\end{equation}
We give some basic identities, valid for $R=k[x_1,\ldots ,x_r]$. 
\begin{lemma}\label{compare}
We have for arbitrary vector spaces $V\subset R_j$,
\begin{align}\label{etwosame}
R_sR_t V&=R_{s+t}V \text { if } s,t\le 0 \text { or } s,t\ge 0;\\
R_sR_tV&\subset R_{s+t}V \text { if } s\ge 0 \text { or } t\le
0;\label{etwosub}\\ R_sR_tV&\supset R_{s+t}V \text { if } s\le 0 \text {
or } t\ge 0.\label{etwosup}
\end{align}
Also, 
\begin{align}\label{e3equal}
R_sR_tR_uV&=R_{s+t+u}V \text { if } s,t,v \text { have the
same sign,  or if  sign $s=$ sign $u$ and } \mid
t\mid \le \mid s\mid,\mid u\mid.\\
R_sR_tR_uV&\subset R_{s+t+u}V \text { if } s,s+t\ge 0 \text { or } u,
t+u\le 0\label{e3sub}\\
R_sR_tR_uV&\supset R_{s+t+u}V \text { if } s,s+t\le 0 \text { or }
u,t+u\ge 0.\label{e3sup}
\end{align}
\end{lemma}\noindent
The proofs are immediate from the definitions. The following Lemma gives
a normal form for relations, that need not be unique. 
\begin{lemma}\label{normalformrel} Let $W$ be related to $V$. Then there
is an expression
$W=R_{i_k}\cdot R_{i_{k-1}}
\cdots  R_{i_1}V$ satisfying
\begin{enumerate}[i.]
\item\label{normalformalt}
The sequence $i_1,\ldots ,i_k$ is alternating in sign.
\item\label{normalunimod} $\exists t, 1\le t \le k$ such that $\mid
i_1\mid <
\cdots <\mid i_t\mid$, and if $k>t, \mid i_t\mid\ge\mid I_{t+1}\mid \ge
\cdots
\ge
\mid i_s\mid$.
\end{enumerate} 
\end{lemma}
\begin{proof} First, using \eqref{etwosame} to collect $R_a\cdot R_b$ for
which sign $a=$ sign $b$,  we may assume the expression is alternating in
sign and is no longer than the original expression. Then using
\eqref{e3equal} we collect adjacent triples $R_a\cdot R_b\cdot R_c$ in
the expression for
$W$, for which 
$\mid b\mid\le \mid a\mid,\mid c\mid. $ Since collecting terms shortens
the length of the relation, after a finite number of steps of
collecting such triples and assuring that the signs alternate, we will
arrive at an expression where the indices alternate in sign, and for which
each adjacent triple $R_a\cdot R_b\cdot R_c$ we have $\mid b\mid >\mid
a\mid,\mid c\mid$. This is possible only if the indices satisfy the
condition \eqref{normalunimod}. 
\end{proof}\par
One might ask whether $W$ related to $V$ and $V$ related to $W$ imply
equality
$\overline{V}=\overline{W}$. We will shortly show that this holds when
$r=2$ (Corollary \ref{reltauequal}). The following counterexample when
$r=3$ is due to David Berman
\cite{Be}. 
\begin{example}\label{berman} (D. Berman: loops in the natural partial
order). Let
$V=\langle x^2y^3,y^2z^3,x^3z^2\rangle\subset R_5, R=k[x,y,z]$, and let
$W=R_2V$. Then $V=R_{-2}W$ but $R_{-1}W$ contains $x^2y^2z^2$, which is
not in
$R_1V$, hence $\overline{V}\ne\overline{W}$. 

\end{example}
 We now restrict to $r=2$.
\begin{proposition}\label{numbrelated} Supppose that $r=2$ and $V\subset
R_j$ satisfies
$\tau (V)=\tau$. Then there are at most
$2^\tau -1$ non-zero equivalence classes $\overline{W}$ of vector spaces
related to $V$. Any nonzero $W$ related to $V$ has an expression of
length $k\le \tau (V)-\tau(W)+1$.
\end{proposition}
\begin{proof}
When $\tau (V)=1$, Lemma \ref{taugen} implies that the vector space $V$
satisfies
$V=f\cdot R_{j-d}$, and $\overline{V}=(f)$. Evidently, any nonzero
$W$ related to $V$ must satisfy $\overline {W}=(f)$. Let
$n>1$ and assume inductively that the statement is true for all $j$, for
vector spaces $V$ satisfying
$\tau (V)\le n -1$. Let $V\subset R_j$ satisfies $\tau
(V)=n$, and let $u,v$ be the minimum positive integers such that
$\overline{R_{-u}V}$ and $\overline{R_vV}$ are each not equivalent to $V$.
Since both $\tau (R_{-u}V)\le n-1 $ and $\tau (R_v(V))\le n-1$, the
induction step would follow from the following claim, as we would then
have that the number of classes $\overline{W}$ related to $V$ would
satisfy
\begin{align*}
\# \{ \,\overline{W}\text { related to } V\}&=\# \{\,\overline{W}\text {
related to } R_{-u}V\}+\#\{\,\overline{W}
\text {related to} R_vV \}+
\text {  one for } \overline{V}\\
&\le 2(2^{n -1}-1)+1=2^n -1.
\end{align*} 
\par\noindent
{\it Claim}:
Let $W\ne 0$ be related to $V$, and assume $\overline{W}\ne\overline{V}$.
Then $W$ is related to $R_{-u}V$ or to $R_vV$, where $u,v$ are defined
above.
\smallskip\par\noindent
{\it Proof of claim}.
We first observe that
\begin{equation}\label{erelreduction}
\overline{R_wV}=\overline{V} \Rightarrow R_aR_wV=R_{a+w}V\text { for }
a\in \mathbb Z. 
\end{equation}
When sign $a=$ sign $w$, this is just \eqref{etwosame}; when sign $a \ne
$ sign $w$ and $\mid a\mid \ge \mid w\mid$ then\
\begin{align*}
R_a\cdot R_w &= R_{a+w}R_{-w}\cdot R_wV  \text { by } \eqref{etwosame}
\text { as sign $a+w= $ sign $-w$ }\\
&=R_{a+w}V  \text { since } \overline{V}=\overline {R_wV}.
\end{align*}
Suppose now that $W$ is related to $V$. Unless
$\overline{V}=\overline{W}$, by
\eqref{erelreduction} we may assume that in the expression
$ W=R_{i_k}\cdot R_{i_{k-1}} \cdots 
R_{i_1}V$ for $W$ we have $i_1\le -u$ or $i_1\ge v$. Then by
\eqref{etwosame} $R_{i_1}V=R_{i_1+u}\cdot R_{-u}V$ in the first case, or
$R_{i_1}V=R_{i_1-v}R_vV$ in the second case. This completes the proof of
the Claim, and of the first statement of the Proposition.\par
The Claim and above proof shows that we need only allow at most one factor
of the form
$R_{i_t}$ in the expression for $W$ for each reduction by one in $\tau$,
and one more for the last step, giving us $k\le \tau (V)-\tau(W)+1$ as
claimed.
\end{proof}
\begin{corollary}\label{reltauequal} Let $r=2$, and suppose that $V\subset
R_j$ and
$ W\subset R_w$ satisfy $W$ is related to $V$ in the sense of
\eqref{erelationexp}, and also $V$ is related to $W$. Then
$\overline{V}=\overline{W}$.
\end{corollary}
\begin{proof}By repeated application of Proposition
\ref{taubasic}\eqref{taubasici}, we have $\tau ( W)\le \tau (V)$,
and vice-versa, hence $\tau (W)=\tau(V)$. Then there is an expression
$W=R_aV$ by the second part of Proposition \ref{numbrelated}.  Proposition
\ref{taubasic}\eqref{taubasiciii} now implies that
$\overline{V}=\overline{W}$.
\end{proof}

\subsubsection{Open problems}
\begin{enumerate}[A.]
\item The dimension
and closure results of Theorems \ref{dimstrata},
\ref{codpartition}, and \ref{closureofstrata} have a naturality that
suggest they might extend to strata not only by the Hilbert function and
partial Hilbert functions (analagous to
\cite[\S 4B]{I2}), but also to more refined strata closer to the complete
Hilbert function where the dimension of each vector space $W$
related to
$V$ is specified (see Section \ref{related} and \cite{Be}). For example
suppose that
$D(u,v)(V)=\dim R_uR_vV$ is specified for all $u,v$: what is the
dimension and closure of the stratum of
$\Grass(d,R_j)$ determined by $D=\{ D(u,v)\}$?
\item The desingularization morphism $G(H)\to \overline{\Grass_H(d,j)}$
is a semi-small resolution. What can be said about
the singularities of $\overline{\Grass_H(d,j)}$? What is the class of
$\overline{\Grass_{H'}(d,j)}$ in  the homology ring $H_\ast (G(H))$? Is
$\overline{\Grass_{H}(d,j)}$ Cohen-Macaulay? A.
King and C. Walter have shown that the homomorphism
$i_\ast :H_\ast (G(H))\hookrightarrow \prod_{\mu\le i \le s}H_\ast
(\Grass(i+1-H_i,R_i))$  is an inclusion \cite{KW}.
\item In Corollary \ref{NTintprop} we showed that
$\Grass_H(d,j)=\LA_N(d,j)\cap
\GA_T(d,j)$, is a proper intersection in $\Grass_\tau (d,j)$. Thus, the
only condition tying $\LA_N(d,j)$ and $\GA_T(d,j)$, with $N=N_H$ and
$T=T_H$ is that
$\tau(N)=\tau(T)$. Do these subvarieties intersect transversely? 
\item Is there a relation between the cohomology rings 
$H^\ast (\overline{\LA _N(d,j))}$ and $H^\ast (\overline{\GA _T(d,j)}$,
when the related partitions $A,B$ correspond? Or a relation between
$H^\ast (\overline{\LA _N(d,j))}$ and $H^\ast (\overline{\LA _N'(d,j))}$
when the partition $A'$ determining $N'$  has one more part than the
partition
$A$ determining $N$?
\item There is a well-known geometric interpretation of the Hilbert
function stratum
$\GA_T(d,j)$. The vector space $V$ determines a rational curve
$X\subset{\mathbb P}^{d-1}$; the restriction $\mathcal T$ to $X$ of the
tangent bundle to $\mathbb P^{d-1}$ decomposes into a direct sum of the
line bundles
$\mathcal T\cong\oplus \mathcal O(-j-d_i)$ where $D$ is the partition
we defined in Definition \ref{defA-D}
\cite{GhISa}. 
Also, the partition $C$ corresponds to the generator degrees of the
ancestor ideal $\overline{V}$, and these are related to the minimum
dimension rational scroll containing the rational curve determined by (a
basis of) $V$ \cite{I5}. Is there a natural geometric interpretation
of the pair $C,D$, that could generalize to other curves in
${\mathbb P}^{d-1}$? 
\end{enumerate}
\begin{ack}
We acknowledge gratefully the many conversations with J. Emsalem
that have influenced the development of this article since the original
preprint. We are grateful also for a collaboration with
F. Ghione and G. Sacchiero in \cite{GhISa}, the results of which have
influenced Section \ref{minressub} of this work, and as well the
collaboration with V. Kanev on \cite{IK} which stands as a reference.
We have benefited from the interest of several colleagues in the
Macaulay inverse systems and level algebras, in particular A. Geramita
and his collaborators, also M. Boij, and Y. Cho, our collaborator on
\cite{ChoI}. We are grateful to the referee for many  
helpful suggestions to improve clarity.
\end{ack}
\bibliographystyle{amsalpha}

\begin{thebibliography}{ACGHOM}
\addcontentsline{toc}{section}{References}

\bibitem[Be]{Be}
Berman D.: \emph{Simplicity of a vector space of forms: finiteness of the
number of complete Hilbert functions}, J. Algebra {\bf 45} (1977), 52--57.

\bibitem[BiGe]{BiGe}
Bigatti A., Geramita A.: \emph{Level algebras, lex segments, and minimal
Hilbert functions}, Communications in Algebra {\bf 31} no. 3 (2003),
1427--1451. 

\bibitem[Bj]{Bj1}
Boij M.: \emph{Betti numbers of compressed level algebras}, J. Pure Appl.
Algebra {\bf 134} (1999), no. 1, 11--16.

\bibitem[BrPV]{BPV}
Brugui\`{e}res A.: \emph{Fibr\'{e}s de Harder-Narasimham et Stratification de
Shatz}, p. 81--204 in
Le Potier J., Verdier J.-L., Module des Fibr\'{e}s Stables sur les Courbes
Alg\`{e}briques, Progress in Math. Vol. 54 (1985), Birkhauser-Boston, Boston,
MA.

\bibitem[BrH]{BH}
Bruns W. , Herzog J.: \emph{Cohen-Macaulay Rings}, Cambridge Studies in
Advanced Mathematics \# 39,
 Cambridge University Press, Cambridge, U.K., 1993; revised paperback
edition, 1998.

\bibitem[BuEi]{BE}
Buchsbaum D., Eisenbud D.: \emph{Algebra structures for finite free
 resolutions, and some structure theorems for codimension three},
 Amer. J. Math. {\bf 99} (1977), 447--485.

\bibitem[Ca]{Ca}
Carlini E.: \emph{Varieties of simultaneous sums of powers for binary
forms}, preprint, 2002, AG/0202050.

\bibitem[CaCh]{CaCh}
\bysame , Chipalkatti J.: \emph{On Waring's problem for several algebraic
forms}, preprint, 2001, AG/0112110.

\bibitem[ChGe]{ChGe}
Chipalkatti J., Geramita A.: \emph{On  parameter spaces for Artin level
algebras}, preprint, 2002, 25 p., AG/0204017,to appear, Michigan Math. J.

\bibitem[ChoI]{ChoI}
Cho Y., Iarrobino, A.: \emph{Hilbert functions of level algebras}, Journal
of Algebra {\bf 241} (2001), 745--758.

\bibitem[Di]{D}
 Diesel S. J.: \emph{Some irreducibility and dimension theorems for
families
 of height 3 Gorenstein algebras}, Pacific J. Math. {\bf 172} (1996),
365--397.

\bibitem[DF]{DiFo}
Dionosi C., Fontanari C.: \emph{Grassmann defectivity \`{a} la Terracini},
preprint, 2001.

\bibitem[EmI1]{EmI}
Emsalem J., Iarrobino A.: \emph{Inverse system of a symbolic power I}, J.
Algebra {\bf 174} (1995), 1080-1090.

\bibitem[FL]{FL}
\bysame ,Laksov D.: \emph{Compressed algebras}, Conf. on
Complete Intersections
in Acireale, (S.Greco and R. Strano, eds), Lecture Notes in  Math.
\# 1092, Springer-Verlag, Berlin and New York, 1984, pp. 121--151.

\bibitem[G]{Ge}
Geramita, A.~V.: \emph{Inverse systems of fat points}, Queen's Papers in
Pure and Applied Mathematics, Vol X, Queens University, 1995.

\bibitem[GHS1]{GHS1}
\bysame , Harima T., Shin Y.S.: \emph{Some special configurations of
points in $\mathbb P^n$}, preprint, 32p., 2002, to appear, J. Algebra.

\bibitem[GHMS1]{GHMS}
\bysame ,\bysame , Migliore J., Shin Y.~S.: \emph{Some remarks on the
Hilbert function of a level algebra};  Appendix: \emph{Construction of
all codimension three level Artinian Hilbert functions of
socle degree five, or of type two and socle degree six}, preprint, 2003.

\bibitem[GhISa]{GhISa}
Ghione F., Iarrobino A., Sacchiero G.: \emph{Restricted Tangent Bundles of
rational curves in ${\mathbb P}^n$}, preprint.

\bibitem[Go1]{Go}
Gotzmann G.: \emph{Eine Bedingung f\"{u}r die Flachheit und das
Hilbertpolynom eines graduierten Ringes}, Math. Z. {\bf 158} (1978),
no. 1, 61--70.

\bibitem[GreK]{GrK}
 Greene C., Kleitman D.J.: \emph{Longest chains in the
lattice of integer partitions ordered by majorization}, European J.
Combinatorics {\bf 7} (1986), 1--10.

\bibitem[Gro]{Gro}
Grothendieck A.: \emph{Techniques de construction et th\'{e}or\`{e}mes
d'existence en g\'{e}ometrie Alg\`{e}brique}, Sem. Bourbaki \# 221
(1961), or Fondements de la G\'{e}ometrie Alg\`{e}brique, Sem.
Bourbaki, 1957-1962, Secretariat Math. Paris (1962).

\bibitem[HN]{HN}
Harder G., Narasimham M.: \emph{On the cohomology groups of moduli spaces},
Math. Ann. {\bf 212} (1975), 215--248.

\bibitem[I1]{I1}
Iarrobino A.: \emph{Vector spaces of forms I, Ancestor ideals of a
vector space of forms}, preprint, 56 p., 1975.

\bibitem[I2]{I2} 
\bysame : \emph{Punctual Hilbert schemes}, A.M.S. Memoirs Vol. 10,
no. 188, Amer. Math. Soc. 1977.

\bibitem[I3]{I3}
\bysame : \emph{Deforming complete intersection Artin algebras. Appendix:
Hilbert functions of $\mathbb C [x,y]/I$}, Proc. of Symposium in Pure
Math. Vol. 40 Part I (1983), 593-608.

\bibitem[I4]{I4}
\bysame : \emph{Compressed algebras: Artin algebras having given socle
degrees and maximal length}, Transactions of the A.M.S., vol. 285, no. 1
pp. 337-378, 1984.

\bibitem[I5]{I5}
\bysame : \emph{Rational curves on scrolls and the restricted tangent
bundle: the ancestor ideal of a vector space of forms in $k[x,y]$},
preprint.

\bibitem[I6]{I6}
\bysame : \emph{Betti strata of height two ideals}, preprint.

\bibitem[IK]{IK}
\bysame , Kanev V.: \emph{Power Sums, Gorenstein Algebras, and
Determinantal Loci}, 345+xxvii p. (1999)  Springer Lecture Notes in
Mathematics \# 1721, Springer, Heidelberg.

\bibitem[IKl]{IKl}
\bysame, Kleiman S.: \emph{The Gotzmann
theorems and the Hilbert scheme}, Appendix C., p. 289--312, in A. Iarrobino and V. Kanev,
Power Sums, Gorenstein Algebras, and Determinantal Loci,  (1999), 
Springer Lecture Notes in Mathematics \# 1721, Springer, Heidelberg.

\bibitem[IY]{IY}
\bysame , Yam\'{e}ogo J.: \emph{The family $\mathrm{G_T}$ of graded
Artinian quotients of $k[x,y]$ of given Hilbert function}, preprint,
37p. 2002. (revision of \emph{Graded ideals and partitions, II:
ramificaton and a generalization of Schubert calculus}, preprint, 44p.,
alg-geom/9709021.)

\bibitem[KW]{KW}
King A., Walter C.: \emph{On Chow rings of fine moduli spaces}, J. reine
angew. Math. {\bf 461} (1995), 179--187.

\bibitem [Klp]{Klp}
Kleppe J.~O. : \emph{The smoothness and the dimension
of PGOR(H) and of other strata of the punctual Hilbert scheme}, J.
Algebra {\bf 200} (1998), 606--628.

\bibitem[Mac1]{Mac1}
Macaulay F. H. S.: \emph{The Algebra of Modular Systems},
Cambridge Univ. Press, Cambridge, U. K. (1916);
reprinted with a foreword by P. Roberts, Cambridge Univ. Pres4).

\bibitem[Mac2]{Mac2}
\bysame : \emph{Some properties of enumeration in the theory of
 modular systems}, Proc. London Math. Soc. {\bf 26} (1927), 531--555.

\bibitem[Mall]{Mall}
Mall D.: \emph{Connectedness of Hilbert function strata and other
connectedness results}, J. Pure Appl. Algebra
{\bf 150} (2000), no. 2, 175--205.

\bibitem[Par]{Par}
Pardue K.: \emph{Deformations of graded modules and connected loci on the
Hilbert scheme}, Queen's Papers Pure Appl. Math.  {\bf 105} (1997),
131--149. 

\bibitem[Ra]{Ra}
Ramella L.: \emph{La stratification du sch\'{e}ma de Hilbert des courbes
rationelles de ${\mathbb P}^n$ par le fibr\'{e} tangent restreint}, C. R. Acad.
Sci. Paris S\'{e}r. I Math. {\bf 311} (1990), no. 3, 181--184.

\bibitem[St1]{St}
Stanley R.: \emph{Hilbert functions of graded algebras}, Advances in
Math. 28, 1978, {\bf 78}, 57--83.

\bibitem[St2]{St2}
\bysame : \emph{Enumerative Combinatorics, Vol 1.}, Wadsworth \&
Brooks/Cole, Belmont, California (1986).

\bibitem[Ve]{Ve}
Verdier, J.-L. : \emph{Two-dimensional $\sigma$-models and harmonic maps
from $S^2$ to $S^n$}, in Group theoretical methods in physics (Istanbul,
1982), 136--141, Lecture Notes Phys. 180, Springer, Berlin, 1983.

\end{thebibliography}

\end{document}